\documentclass[a4paper,10pt]{scrartcl}
\usepackage[utf8]{inputenc}
\usepackage{amssymb,amsmath,amsthm}
\usepackage{array, hhline}
\usepackage{tikz}
\usepackage{authblk}
\usepackage{graphicx}
\usepackage{subfig}

\newcommand{\MeshRefinement}{3}

\renewcommand{\vec}{\boldsymbol}

\newcommand{\ud}{\,\mathrm{d}}
\newcommand{\R}{\mathbb{R}}
\renewcommand{\P}{\mathbb{P}}

\newtheorem{defi}{Definition}[section]
\newtheorem{theorem}[defi]{Theorem}
\newtheorem{lemm}[defi]{Lemma}
\newtheorem{rem}[defi]{Remark}
\newtheorem{cor}[defi]{Corollary}
\newenvironment{mproof}{\paragraph{Proof.}}{\hfill$\blacksquare$}

\numberwithin{equation}{section}
\numberwithin{table}{section}
\numberwithin{figure}{section}

\newcommand{\N}{\mathbb{N}}
\newcommand{\ds}{\displaystyle}

\newcommand{\braum}{L^2(I ; L^2(\Omega))}

\newcommand{\la}{\langle}
\newcommand{\ra}{\rangle}

\begin{document}

%\mathtoolsset{showonlyrefs,showmanualtags}

\title{Error analysis for discretizations of parabolic problems using
continuous finite elements in time and mixed finite elements in space}
\author[M.\ Bause et al.]
       {M.\ Bause\thanks{bause@hsu-hh.de (corresponding author),
$^\dag$Florin.Radu@uib.no, $^\ddag$koecheru@hsu-hh.de}\mbox{~}, F.\ A.\ Radu$^\dag$, U.\
K\"ocher$^\ddag$
       \\{\small
       $^\ast \, \ddag$ Helmut Schmidt University, Faculty of Mechanical Engineering,
Holstenhofweg
85, 220433 Hamburg, Germany\\
       $^\dag$ University of Bergen, Department of Mathematics, All\'{e}gaten 41, 50520
Bergen, Norway}
}
\maketitle

\begin{abstract}
{\bfseries Abstract.}
Variational time discretization schemes are getting of increasing importance for the 
accurate numerical approximation of transient phenomena. The applicability and value of 
mixed finite element methods (MFEM) in space for simulating transport processes have been demonstrated in a wide class of 
works. We consider a family of continuous Galerkin-Petrov time discretization schemes 
that is combined with a mixed finite element (MFE) approximation of the spatial 
variables. The existence and uniqueness of the semidiscrete approximation and of the 
fully discrete solution are established. For this, the Banach-Ne\v{c}as-Babu\v{s}ka 
theorem is applied in a non-standard way. Error estimates with explicit rates of 
convergence are proved for the scalar and vector-valued variable. An optimal order 
estimate in space and time is proved by duality techniques for the scalar variable. The 
convergence rates are analyzed and illustrated by numerical experiments, also on 
stochastically perturbed meshes. 
\end{abstract}

\clearpage

\section{Introduction}
\label{Sec:Intro}

Numerical simulations of time dependent single and multiphase phase flow and
multicomponent transport processes in complex and porous media with strong
heterogeneities and anisotropies are desirable in several fields of natural sciences and
civil engineering as well as in a large number of branches of technology; cf.\ e.g.\
\cite{Chen2005,Helmig1997}. Typically, the discretization in space involves a significant
set of complexities and challenges. MFEM (cf.\ \cite{Brezzi1991,Chen2010}) have proved 
their potential and  capability to approximate solutions with high accuracy and physical 
consistency; cf.\ e.g.\ \cite{Bause2010,Celia1990}. So far, the temporal approximation of flows and transport phenomena in 
porous media have received relatively little interest (cf., 
e.g., \cite{Farthing2002,Arbogast96,Brunner2012,Radu2004,Radu2008,Radu2010,Woodward00}
and the references therein) and have been limited to traditional non-adaptive first and
second order methods, even if strong chemical reactions with high temporal variations in
profiles are present. Rigorous studies of higher order time discretizations are still
missing. The low-order implicit time discretization is of particular concern with respect
to numerical diffusion for smooth solutions of transport problems (cf.\ \cite{RaduSuciu}
for a study on numerical diffusion for different temporal and spatial discretizations of
a transport equation).

The Galerkin method is a well-recognised approach to solve time dependent problems;
cf., e.g., \cite{Aziz1989,Thomee2006}. However, until now it has rarely been used
in practice for discretizing the time variable in approximations of initial-boundary
value problems. Since recently, variational time discretization schemes based on
continuous or discontinuous finite element techniques have been developed to the point
that they can be put into use (cf.\ \cite{Hussain2013,Hussain2012}) and demonstrate their 
significant advantages. Higher order methods are naturally embedded in these schemes and 
the uniform variational approach simplifies stability and error analyses. Further, 
goal-oriented error control \cite{Bangerth2003} based on the dual weighted residual 
approach relies on variational space-time formulations and the concepts of adaptive finite 
element techniques for changing the polynomial degree as well as the length of the time 
intervals become applicable. Variational time discretization schemes that are combined 
with continuous or
discontinuous finite element methods for the spatial variables are studied for flow and 
parabolic problems in, for 
instance, \cite{Ahmed2015,Ahmed2012,Ahmed2012_2,Andreev2014,Basting2013,Bause2015,
Hussain2013, Hussain2012,Hussain2011,Matthies2011,Sudirham2006} and for wave problems in, 
for instance, \cite{Bangerth2010,Koecher2014,Koecher2015}. In these works algebraic
formulations of the variational time discretizations are developed
\cite{Andreev2014,Hussain2012,Hussain2013,Koecher2014,Koecher2015,Sudirham2006},
preconditioning techniques for the arising block matrix systems are addressed
\cite{Andreev2014,Basting2013,Hussain2011,Koecher2015} and, finally, computational
studies are performed. 

Numerical analyses of semidiscretizations in time by variational methods and of  variational space-time approaches can be found in, for instance, 
\cite{Cesenek2012,Karakashian2004,Karakashian1999,Schieweck2010,Thomee2006}.
In \cite{Thomee2006} discontinuous variational approximations of the time variable are
studied for abstract parabolic problems whereas in \cite{Schieweck2010} their continuous
counterparts are analyzed. In \cite{Cesenek2012,Sudirham2006} discontinuous variational 
approximations in time and space are studied and error estimates are proved. In 
\cite{Sudirham2006} time-dependent domains are considered in an arbitrary Lagrangian
Eulerian (ALE) framework and the advection-diffusion equation is written in mixed form as 
a system of first order equations in space. In \cite{Ern2016} a discontinuous Galerkin 
method in time combined with a stabilized finite element approach in space for first 
order partial differential equations is investigated for static and dynamically changing 
meshes. Error estimates in the $L^\infty(L^2)$ and $L^2(L^2)$ norm are derived. In \cite{Karakashian2004,Karakashian1999} continuous 
space-time approximations for nonlinear wave equations with mesh modifications and for 
the Schr\"odinger equation are considered. Existence and uniqueness of the discrete 
solutions are discussed and error estimates are proved for the schemes. 

As far as the MFE approximation of parabolic problems is concerned, in 
\cite{Thomee2006} an error estimate for the semidiscretization in space is given. However, 
for the flux variable an error estimate is proved only for the $\vec L^2$ norm. No 
estimate is provided for the error in divergence of the flux, that is part of the natural 
norm of the underlying function space $\vec H(\mathrm{div};\Omega)$. In 
\cite{Cristina1987} and \cite{Johnson1981} similar error estimates, also in negative 
norms, are presented. In particular, estimates similar to the error estimates for 
conventional finite element approximations are established. The singular behavior of the 
error estimates as $t\rightarrow 0$ for initial data in $L^2(\Omega)$ is further 
included. 

In this work a continuous Galerkin-Petrov (cGP) method is used for the discretization in 
time, whereas the MFEM \cite{Brezzi1991,Chen2010} is applied for the spatial 
discretization. Appreciable advantages of the MFEM are its local mass conservation 
property and the inherent approximation of the flux field as part of the formation 
itself. In simulating coupled flow and transport processes in porous media the flux 
approximation of the flow problem is usually of higher practical interest than the 
approximation of the scalar variable itself. To the best of our knowledge, rigorous error 
estimates for fully discrete variational space-time discretization schemes that are based 
on MFE approximations are still missing. In our numerical analysis we split the temporal 
discretization error from the spatial one by introducing an auxiliary problem based on 
the semidiscretization in time. We firstly estimate the temporal discretization error 
and secondly the error between the semidiscrete and the fully discrete solution. The order of convergence estimates are derived in the natural norms of the variational space-time approach. They are summarized in Thm.\ \ref{theorem_error_cont2fully}. {\color{black} For the scalar variable of the MFE approach one of the given error estimates, measured in the norm of $L^2(0,T;L^2(\Omega))$, is optimal in space and time if a certain regularity assumption is supposed to be satisfied. For constant scalar-valued diffusion coefficients an error estimate for the flux variable in the norm of $L^2(0,T;\vec L^2(\Omega))$ is further provided. It is optimal in space and suboptimal in time. In the Gaussian quadrature points of the temporal discretization optimal order error estimates for the flux variable in $\vec L^2(\Omega)$ are even obtained for heterogeneous diffusion matrices.} The existence and uniqueness of the semidiscrete and fully discrete solution is further established. Even though a prototype
model problem is studied here only, we believe that the techniques for analyzing mixed variational space-time approximation schemes can be applied similarly to more complex flow and transport problems in porous media. 

This work is organized as follows. In Sec.\ \ref{Sec:FullDisSchm} our fully discrete
variational space-time method is developed. In Sec.\ \ref{Sec:ExUniErr} we address the
semidiscrete problem by proving existence and uniqueness of its solution and error
estimates for the semidiscretization in time. In Sec.\ \ref{Sec:MFEMcGDis} we study the
fully discrete problem and show the existence and uniqueness of its solution. The error
between the semidiscrete and fully discrete problem is estimated. In Thm.\ 
\ref{theorem_error_cont2fully} an error estimate for the simultaneous space-time
discretization is provided by combining the before-given estimates of the temporal and
spatial discretization. In Sec.\ \ref{Sec:NumStudies} we illustrate and validate our 
derived error estimates by numerical experiments. We end our work with some conclusions 
in Sec.\ \ref{Sec:Conclusions}.

\section{The fully discrete variational scheme}
\label{Sec:FullDisSchm}

\subsection{Notation and preliminaries}
\label{Sec:NotPrem_1}

Throughout this paper, standard notations are used. Let $\Omega\subset  \R^d$, with $d=2$ or $d=3$, be an polygonal 
or polyhedral bounded domain. We denote by $H^p(\Omega)$ the Sobolev space of $L^2$ functions with derivatives up to order $m$ in $L^2(\Omega)$ and by $\langle \cdot,\cdot \rangle$ the inner product in $L^2(\Omega)$. Sobolev spaces of vector-valued functions are written in bold letters. Further, let 
$H^1_0(\Omega)=\{u\in H^1(\Omega) \mid u=0 \mbox{ on } \partial \Omega\}$ and  $H^{-1}(\Omega)$ denote its dual space. For the norms of the Sobolev 
spaces the notation is 
\begin{align*}
\| \cdot \| := \| \cdot\|_{L^2(\Omega)}\,,\qquad 
\| \cdot \|_p := \| \cdot\|_{H^p(\Omega)}, \,\, \mbox{ for } p \in \N, p \ge 1\,.
\end{align*} 
For the mixed problem formulation we use the abbreviations
\begin{equation*}
 \vec V = \vec H(\mathrm{div};\Omega)= \{\vec q \in \vec L^2(\Omega) \mid \nabla \cdot \vec q \in L^2(\Omega)\} \,, \quad W = L^2(\Omega)\,,
\end{equation*}
and 
\[
\|\vec v\|_{\vec V}:= (\|\vec v\|^2 + \|\nabla \cdot \vec v \|^2)^{1/2}\,.
\] 

Let $X_0 \subset X \subset X_1$ be three reflexive Banach spaces with continuous embeddings. Then we consider
the following set of Banach space valued function spaces,
\begin{align*}
C(\overline{I};X) & =
\{ w : [0,T] \rightarrow X \mid \mbox{$w$ is continuous} \}\,,\\
L^2(I;X) & =
\bigg\{w: (0,T) \rightarrow X \;\; \bigg|\;\;
\int_I \| w(t) \|_X^2 \; \ud t < \infty \bigg\}\,,\\
H^1(I;X_0,X_1) & =
\{w \in L^2(I;X_0) \mid \partial_t w \in L^2(I;X_1)\}\,,
\end{align*}
that are equipped with their naturals norms (cf.\ \cite{Ern2010}) and where the
time derivative $\partial_t$ is understood in the sense of distributions on
$(0,T)$. In particular, every function in $H^1(I;X_0,X_1)$ is continuous on
$[0,T]$ with values in $X$; cf.\ \cite{Ern2010}. For $X_0=X=X_1$ we simply
write $H^1(I;X)$. Moreover, we put $H_0^1(I;X)=\{u\in H^1(I; X)\mid u(0)=0\}$.

For $u\in H^1_0(\Omega)$ let $A:H^1_0(\Omega)\mapsto H^{-1}(\Omega)$ be defined uniquely 
by
\begin{equation}
\label{Eq:Def_A}
\langle Au, v\rangle = a(u,v) \quad \mbox{for all}\;  u,v\in H^1_0(\Omega)
\end{equation}
with
\[
a(u,v):= \langle \vec D \nabla u, \nabla v\rangle\,,
\]
where the matrix  $\vec D=\vec D(\vec x)=(d_{ij}(\vec
x))_{i,j=1}^d$ satisfies $d_{ij}\in L^\infty(\Omega)$ and is elliptic with
\begin{equation}
\label{Eq:PosDefD}
D_M | \vec \xi |^2  \geq \vec \xi^\top \vec D(\vec x) \vec \xi  \geq D_m |\vec \xi
|^2\,, \qquad \theta_M | \vec \xi |^2  \geq \vec \xi^\top \vec D(\vec x)^{-1} \vec \xi  
\geq \theta_m |\vec \xi |^2 \,,
\end{equation}
for almost every $\vec x \in \Omega$, all $\vec \xi \in \R^d$ and some constants $0 < D_m 
\leq D_M < \infty$. In \eqref{Eq:PosDefD} we put $\theta_m := {D_M}^{-1}$ and $\theta_M := {D_m}^{-1}$. Under the previous assumptions it holds that
\begin{align}
 a(v,v) & \geq \alpha \|v\|^2_{1} \quad \mbox{for all}\; v \in
H^1_0(\Omega)\,,
\label{Eq:Coerc_1} \\[1ex]
 |a(u,v)| & \leq \beta \|u\|_{1}\|v\|_{1} \quad \mbox{for all}\;
u,v\in H^1_0(\Omega)\,.
\label{Eq:Coerc_2}
\end{align}
Thus, $A: H_0^1(\Omega) \mapsto H^{-1}(\Omega)$ is a linear and continuous operator. For a subspace $D(A) \subset H^1_0(\Omega)$ let $A:D(A)\mapsto H^{-1}(\Omega)$ be a bijective linear continuous operator. For instance, if $\Omega$ is a convex polygonal or
polyhedral bounded domain and $d_{ij}\in W^{1,\infty}(\Omega)$, for $i,j=1,\ldots
d$, is satisfied, then the operator $A$ is a bijective linear continuous
operator from $D(A) = H^2(\Omega)\cap H^1_0(\Omega)$ to $L^2(\Omega)$; cf.\ 
\cite{Grisvard85}.

Due to the properties \eqref{Eq:Coerc_1}, \eqref{Eq:Coerc_2} of the
bilinear form $a(\cdot,\cdot)$ the lemma of Lax--Milgram ensures that the
operator $A: H_0^1(\Omega)\mapsto H^{-1}(\Omega)$ defined in \eqref{Eq:Def_A} is
invertible and satisfies in the corresponding operator norm the stability estimates
\begin{equation*}
\label{Eq:Stab_A}
\| A\|_{} \leq \beta \quad \mbox{and} \quad \| A^{-1}\| \leq \alpha\,.
\end{equation*}
Moreover, for all $ g\in H^{-1}(\Omega)$ it holds that
\begin{equation}
\label{Eq:Prop_A}
\langle g, A^{-1} g \rangle = \langle A A^{-1} g, A^{-1}g\rangle \geq \alpha \|
A^{-1} g\|^2_{1} \geq \frac{\alpha}{\beta^2}\, \|g\|^2_{H^{-1}(\Omega)}\,.
\end{equation}

As usual, by $c>0$ we denote a generic constant throughout the paper. 

\subsection{Problem formulation}
\label{Sec:NotPrem_2}

As a prototype model for more sophisticated multiphase flow and multicomponent
reactive transport systems in porous media (cf.\ e.g.\ \cite{Chen2005,Helmig1997}) we
study in this work 
\begin{align}
\partial_t u - \nabla \cdot (\vec D \nabla u ) & = f  && \mathrm{in}\;
\Omega \times I\,, \label{Eq:Diff}\\
u & = 0 && \mathrm{on}\;  \partial \Omega \times I\,, \label{Eq:HomDir}\\
u(\cdot, 0) & = u_0 && \mathrm{in}\; \Omega\,,\label{Eq:InVal}
\end{align}
equipped with homogeneous Dirichlet boundary conditions for simplicity only, where $I=(0,T]$ with final time $T>0$ and the diffusion matrix $\vec D$ satisfies the assumptions made in the previous subsection.

Let $f\in L^2(I;W)$ and $u_0\in H^1_0(\Omega)$ be given. Then the existence of a unique weak solution
\begin{equation}
\label{Eq:RegSol}
u \in L^2(I;H^1_0(\Omega)) \cap H^1(I;W) \cap C(\overline{I};W)
\end{equation}
to \eqref{Eq:Diff}--\eqref{Eq:InVal} is ensured; {\color{black} cf.\ \cite[p.\ 382, Thm.\ 5]{Evans2010}}. {\color{black} We note that \eqref{Eq:RegSol} already provides an improved regularity for the weak solution of \eqref{Eq:Diff}--\eqref{Eq:InVal}; cf.\ \cite[p.\ 378, Thm.\ 3]{Evans2010}}. 

In order to derive our family of discretization schemes, we first define the auxiliary flux variable $\vec q := - \vec D \nabla u$ for the weak solution $u$ of  \eqref{Eq:Diff}--\eqref{Eq:InVal} that is given by \eqref{Eq:RegSol}.  Since $\partial_t u \in L^2(I;W)$ is satisfied by \eqref{Eq:RegSol} and $f\in L^2(I;W)$ holds by assumption, it directly follows that $\vec q\in L^2(I;\vec V)$. {\color{black} The pair $\{u,\vec q\} \in H^1(I;W)\cap C(\overline{I};W) \times L^2(I;\vec V)$ is then also the unique solution to the set of variational equations}
\begin{align}
\int_0^T \langle  \partial_t u, w \rangle \ud t + \int_0^T \langle
\nabla \cdot \vec q , w\rangle \ud t  & = \int_0^T \langle f, w \rangle \ud t \,,
\label{Eq:IntMix_1}\\[1ex]
\int_0^T \langle \vec D^{-1} \vec q, \vec v\rangle \ud t  - \int_0^T
\langle u, \nabla \cdot \vec v \rangle \ud t  & = 0  \label{Eq:IntMix_2}
\end{align}
for all $w \in L^2(I;W)$ and $\vec v \in L^2(I;\vec V)$ and satisfies the initial condition $u(0)=u_0$. {\color{black} To find \eqref{Eq:IntMix_2} integration by parts was used.} The global problem formulation \eqref{Eq:IntMix_1}, \eqref{Eq:IntMix_2} motivates our semidiscretization in time.

{\color{black}
\begin{rem}
\label{Rem:Assump}
\begin{itemize}
\item Below, in order apply Lagrange interpolation in time to the function $f$, we need the stronger assumption that $f\in C([0,T];W)$ is satisfied.
\item Below, we introduce a semidiscrete approximation in time of the flux $\vec q$ in a subspace of $C([0,T];\vec V)$. For this we need to assume that $\vec D\nabla u_0\in \vec V$ holds.
\item Higher order regularity of weak solutions to \eqref{Eq:Diff}--\eqref{Eq:InVal}, that is need below for the proof of higher order convergence rates, can be obtained under further technical assumptions about the data, coefficients and the boundary of the domain $\Omega$. For the prototype model problem \eqref{Eq:Diff}--\eqref{Eq:InVal} such higher order regularity results are well-known; cf.\ \cite[p.\ 386, Thm.\ 6]{Evans2010}. For (elliptic) regularity results in domains with non-smooth boundaries we refer to, e.g., \cite{Grisvard85,Mazya2000}. Below, we tacitly assume that the required assumptions about the data and $\partial \Omega$ are satisfied such that the existence of a sufficiently regular solution can be assumed. Without such an assumption the application of higher order methods is not meaningful. 
\end{itemize}
\end{rem}
}

\subsection{Variational discretization in time by a continuous Galerkin method }
\label{Sec:SemiTimeCont}

For the discretization in time we decompose the time interval $(0,T]$ into $N$
subintervals $I_n=(t_{n-1},t_n]$, where $n\in \{1,\ldots ,N\}$ and $0=t_0<t_1<
\cdots < t_{n-1} < t_n = T$. Further $\tau$ denotes the discretization
parameter in time and is defined as the maximum time step size $\tau =
\max_{1\leq n \leq N} \tau_n $, where $\tau_n = t_n-t_{n-1}$. We introduce the
function spaces of piecewise polynomials of order $r$ in time,
\begin{align*}
 \mathcal{X}^{r}{(X)} & :=\left\{ u_\tau \in C({\bar{I};\,X}) \;\; \Big| \;\;
u_\tau{}_{|{\overline{I}_n}} \in \P_r(\overline I_n;\,X)\,,\;
\forall n \in \{1,\ldots ,N\} \right\} \,, \\[1ex]
\mathcal{Y}^{r}{(X)} & := \left\{ w_\tau \in L^2({I;\,X}) \;\; \Big| \;\;
w_\tau{}_{|{I_n}} \in \P_{r}(I_n;\,X)\,,\; \forall n \in \{1,\ldots ,N\} \right\}\,,
\end{align*}
where
\begin{equation*}
\P_r(J;\,X) = \bigg\{ p : J \to X \;\; \bigg| \;\;
p(t) = \sum\limits_{j=0}^{r}{\xi_n^j\, t^j}\,,\;
 \xi_n^j \in X\,,\; j=0,\ldots,r \bigg\}
\end{equation*}
and $\mathcal{X}^{r}{(X)}\subset H^1(0,T;W)$. We let 
\[
\mathcal{X}_0^{r}{(X)} = \left\{ u_\tau \in \mathcal{X}^{r}{(X)} \;\; \Big| \;\;
u_\tau(0) = 0 \right\}\,.
\]
Further, we put
\begin{displaymath}
\mathcal{W} = X_0^r (W) \times X^r (\vec V)  \qquad \mbox{and}
\qquad  \mathcal{V} = Y^{r-1} (W) \times Y^{r-1} (\vec V)  \,.
 \end{displaymath}
We equip the function spaces $\mathcal{W}$ and $\mathcal{V}$ with their natural
norms being defined by
\begin{align}
  \|\{u_\tau ,\vec q_\tau \}\|_{\mathcal{W}}^2 & = \|u_\tau \|_{L^2(I;W)}^2 +
\| \partial_t u_\tau
\|_{L^2(I;W)}^2  + \| \vec q_\tau \|^2_{L^2(I;\vec V)}\,,
\label{Eq:SpcTmNrm}
\\[1ex]
\|\{w_\tau ,\vec v_\tau \}\|_{\mathcal{V}}^2 & = \| w_\tau
\|_{L^2(I;W)}^2 +
\|\vec v_\tau  \|^2_{L^2(I;\vec V)}\,.
 \nonumber
 \end{align}
With respect to these norms the space $\mathcal{W}$ is a Banach space and the space $\mathcal{V}$ is a reflexive Banach space.
Further, we define the space-time bilinear form $a_\tau \in \mathcal{L}(\mathcal{W}\times
\mathcal{V};\R)$ by means of
\begin{align*}
a_\tau(\{u_\tau,\vec q_\tau \},\{w_\tau,\vec v_\tau\}) = &  \int_0^T
\Big(\langle \partial_t u_\tau, w_\tau \rangle + \langle \nabla \cdot
\vec q_\tau, w_\tau \rangle \big) \ud t \\[2ex]
& + \int_0^T \langle \vec D^{-1}
\vec q_\tau , \vec v_\tau \rangle \ud t - \int_0^T
\langle u_\tau , \nabla \cdot \vec v_\tau \rangle \ud t
\end{align*}
for $\{u_\tau,\vec q_\tau\}\in \mathcal{W}$ and $\{w_\tau,\vec v_\tau\}\in
\mathcal{V}$. Obviously, the mapping $a_\tau: \mathcal{W} \times
\mathcal{V} \mapsto \R$ is linear and continuous, i.e.
\begin{equation}
\label{Eq:ContA}
|a_\tau(\{u_\tau,\vec q_\tau \},\{w_\tau,\vec v_\tau\}) | \leq c \|
\{u_\tau^0,\vec q_\tau\}\|_{\mathcal{W}}\| \,
\{w_\tau,\vec v_\tau\}\|_{\mathcal{V}}
\end{equation}
with some constant  $c >0$ independent of $\tau$ and $T$.

For the family of continuous variational time discretization schemes the
spaces $\mathcal{X}^{r}{(X)}$ of continuous functions act as spaces for the
solution whereas the spaces $\mathcal{Y}^{r-1}{(X)}$ consisting of piecewise
polynomials that are discontinuous at the end points of the time intervals are
used as test spaces. Since the spaces of the trial and test functions differ
here, a discretization of Galerkin-Petrov type is thus obtained.

A semidiscrete variational approximation of the mixed form of problem
\eqref{Eq:Diff}--\eqref{Eq:InVal}, refered to as the exact form of cG$(r)$, is then defined by solving the variational
equations \eqref{Eq:IntMix_1}, \eqref{Eq:IntMix_2} in discrete
subspaces: \emph{Find $\{u_\tau, \vec q_\tau\}\in \mathcal{X}^r
(W)\times \mathcal{X}^r(\vec V)$ such that}
\begin{align}
\label{Eq:GtdP_1}
\int_0^T \langle \partial_t u_\tau  ,w_\tau  \rangle \ud t + \int_0^T \langle
\nabla \cdot \vec q_\tau, w_\tau \rangle \ud t & = \int_0^T \langle f,
w_\tau \rangle \ud t\,,\\[2ex]
\label{Eq:GtdP_2}
\int_0^T
\langle \vec D^{-1} \vec q_\tau , \vec v_\tau \rangle \ud t - \int_0^T
\langle u_\tau , \nabla \cdot \vec v_\tau \rangle \ud t & = 0\,,
\end{align}
\emph{for all $w_\tau \in \mathcal{Y}^{r-1}(W)$ and $\vec v_\tau \in
\mathcal{Y}^{r-1}(\vec V)$ with the initial conditions that $u_\tau (0) := u_0$ and {\color{black} $\vec q_\tau (0) := - \vec D\nabla u_0$ (cf.\ Rem.\ \ref{Rem:Assump})}.} 

%\label{Rem:IniFluxCond}  

We refer to the solution of Eqs.\ \eqref{Eq:GtdP_1}, \eqref{Eq:GtdP_2} as the
continuous Galerkin--Petrov method with piecewise polynomials of order
$r$ and use the notation cGP($r$). To ensure the existence and uniqueness of solutions to
\eqref{Eq:GtdP_1}, \eqref{Eq:GtdP_2}, it is sufficient to use the test spaces
$\mathcal{Y}^{r-1}(W)$ and $\mathcal{Y}^{r-1}(\vec V)$ with piecewise
polynomials of order $r-1$, since the continuity constraint at the discrete
time points $t_n$, $n=0,\ldots,N-1$, that is implied by the definition of the
solution spaces $\mathcal{X}^r (W)$ and $\mathcal{X}^r(\vec V)$, yields a
further condition. By using discontinuous test basis functions $w_\tau(t) = w
\psi_{n,i}(t)$ and $\vec v_\tau = \vec v \psi_{n,i}(t)$, for
$i=1,\ldots ,r$, with arbitrary time independent functions $w\in W$ and $\vec v
\in \vec V$, respectively, and piecewise polynomial
functions $\psi_{n,i}:I\mapsto \R$ that are of order $r-1$ on $I_n$
and vanish {\color{black} on} $I\backslash\overline{I}_n$, we can recast the variational equations
\eqref{Eq:GtdP_1}, \eqref{Eq:GtdP_2} as a time marching scheme: \emph{For $n=1,\ldots, N$
find $u_\tau{}_{|\overline{I}_n}\in \P_r(\overline{I}_n;W)$ and $\vec
q_\tau{}_{|\overline{I}_n}\in P_r(\overline{I}_n;\vec V)$ such that}
\begin{align}
\label{Eq:locP_1}
\int_{I_n} \langle \partial_t u_\tau  , w \rangle \, \psi_{n,i}(t) \ud t + \int_{I_n}
\langle
\nabla \cdot \vec q_\tau, w \rangle \, \psi_{n,i}(t) \ud t & = \int_{I_n} \langle f,
w\rangle \, \psi_{n,i}(t) \ud t\,,\\[2ex]
\label{Eq:locP_2}
\int_{I_n} \langle \vec D^{-1} \vec q_\tau , \vec v \rangle \, \psi_{n,i}(t) \ud t -
\int_{I_n}
\langle u_\tau , \nabla \cdot \vec v \rangle \, \psi_{n,i}(t) \, \ud t & = 0
\end{align}
\emph{for all $w\in W$ and $\vec v \in \vec V$ and $i=1,\ldots ,r$ with the continuity
constraints $u_\tau{}_{|I_n}(t_{n-1}) =
u_\tau{}_{|I_{n-1}}(t_{n-1})$ and $\vec q_\tau{}_{|I_n}(t_{n-1}) =
\vec q_\tau{}_{|I_{n-1}}(t_{n-1})$  for $n\geq 2$ and the initial 
conditions $u_\tau{}_{|I_n}(t_{n-1}) := u_0$, {\color{black} $\vec q_\tau{}_{|I_n}(t_{n-1}) := -\vec D\nabla u_0$} for $n = 1$.}

To determine $u_\tau{}_{|\overline{I}_n}$ and $\vec q_\tau{}_{|\overline{I}_n}$, we
represent them in terms of basis functions, with respect to the time variable, of the
spaces $\mathcal{X}^r(W)$ and $\mathcal{X}^r(\vec V)$ such that
\begin{equation}
 \label{Eq:RepBasis}
 u_\tau{}_{|\overline{I}_n} (t ) = \sum_{j=0}^r U_n^j\,  \varphi_{n,j}(t) \quad
\mathrm{and} \quad
 \vec q_\tau{}_{|\overline{I}_n} (t ) = \sum_{j=0}^r \vec Q_n^j \, \varphi_{n,j}(t)\,,
\quad
\mathrm{for} \; t \in I_n\,,
\end{equation}
with coefficient functions $U_n^j\in W$ and $\vec Q_n^j \in \vec V$ for $j=0,\ldots, r$
and polynomial basis functions $\varphi_{n,j}\in \P_r(\overline{I}_n;\R)$ that are
Lagrange functions with respect to $r+1$ nodal points $t_{n,j}\in I_n$
satisfying the conditions $\varphi_{n,j} (t_{n,i}) = \delta_{i,j}$ for $i,j=0,\ldots , 
r$. For the treatment of the continuity constraint in time we put $t_{n,0}=t_{n-1}$. The 
other points $t_{n,1},\ldots, t_{n,r}$ are chosen as the quadrature points of the
$r$-point Gaussian quadrature formula on $I_n$ which is exact if the function to be
integrated is a polynomial of degree less or equal to $2r-1$. The basis functions
$\varphi_{n,j}\in \P_r(\overline{I}_n;\R)$ of \eqref{Eq:RepBasis}, for $j=0,\ldots , r$,
are defined, as usual in the finite element framework, via the affine reference 
transformation onto $\hat I = [0,1]$. The test basis functions $\psi_{n,i}\in 
P_{r-1}(\overline{I}_n;\R)$ with $\psi_{n,i} (t_{n,l}) = \delta_{i,l}$ for 
$i,l=1,\ldots,r$ are defined similarly; cf.\ \cite{Bause2015,Koecher2015} for details. 
Now we transform all the time integrals in \eqref{Eq:locP_1}, \eqref{Eq:locP_2} to the 
reference interval $\hat I$. By a subsequent application of {\color{black} the $r$-point Gaussian quadrature formula with weights $\hat \omega_i$ and quadrature nodes $\hat t_i$ on $\hat I$} as well as the further notation
\[
\hat \alpha_{ij} := \hat \omega_i\cdot \dfrac{\ud}{\ud\hat t} \hat \varphi_j(\hat 
t_i) \quad \text{and} \quad \hat \beta_{ij}:= \hat \omega_i \cdot \delta_{i,j}
\]
for $i= 1,\ldots,r$, $j=0,\ldots, r$ (cf.\ \cite{Bause2015,Koecher2014,Schieweck2010}), 
we obtain the following system of variational problems for the coefficient functions 
$U_n^j\in W$ and $\vec Q_n^j\in \vec V$ of the representation
\eqref{Eq:RepBasis}: \emph{For $n=1,\ldots, N$ and $j=1,\ldots,r$ find coefficient 
functions $\{U_n^j,\vec Q_n^j\} \in W\times \vec V$ such that}
\begin{align}
\sum_{j=0}^r \hat \alpha_{ij} \langle U_n^j ,w\rangle + {\tau_n} \,
\hat \beta_{ii} \langle \nabla \cdot \vec Q_n^i ,w\rangle & = {\tau_n} \,
\hat \beta_{ii} \langle f(t_{n,i}) ,w\rangle\,,
\label{Eq:SemiDis1}
\\[1ex]
\langle \vec D^{-1} \vec Q_n^i,\vec v \rangle - \langle U_n^i ,\nabla \cdot
\vec v\rangle & = 0\,,
\label{Eq:SemiDis2}
\end{align}
\emph{for $i=1,\ldots, r$ and all $\{w,\vec v\}\in W\times \vec V$, and where due to continuity in time $U_n^0 = u_\tau{}_{|I_{n-1}}(t_{n-1})$, {\color{black} $\vec Q_n^0 = \vec q_\tau{}_{|I_{n-1}}(t_{n-1})$} for $n\geq 2$ and $U_n^0:= u_0$, {\color{black} $\vec Q_n^0:=-\vec D\nabla u_0$} for $n=1$.}

{\color{black}
\begin{rem}
\label{Rem:IniVal}
In the numerical scheme \eqref{Eq:SemiDis1}, \eqref{Eq:SemiDis2}, the flux coefficient functions $\vec Q_n^j$, for $j=1,\ldots, r$, arise only in the $r$ Gaussian quadrature points $t_{n,1},\ldots, t_{n,r}\in (t_{n-1},t_n)$ of the subinterval $I_n$. Nevertheless, the coefficient functions $\vec Q_n^0$ for $n\geq 1$, are needed for the unique determination of the semidiscrete flux function $\vec q_\tau\in \mathcal X^r(\vec V)$ and an explicit evaluation of $\vec q_\tau{}_{|I_n}$ by the representation \eqref{Eq:RepBasis} in other time points of $I_n$ than in the Gaussian quadrature nodes. The fact that the coefficient functions $\vec Q_n^0$ do not arise in  \eqref{Eq:SemiDis1}, \eqref{Eq:SemiDis2} is due to the definition of the Lagrange basis functions $\varphi_{n,j}$ in \eqref{Eq:RepBasis} and the fact that the time derivative of the flux variable $\vec q$ does not arise in the model equations.   
\end{rem}
}

For the derivation of \eqref{Eq:SemiDis1}, \eqref{Eq:SemiDis2} from \eqref{Eq:locP_1}, 
\eqref{Eq:locP_2} we tacitly replaced the integrand $f$ on the right-hand side by its 
Lagrange interpolate $\Pi_r f \in \P_r(I_n;L^2(\Omega))$ defined by
\begin{equation}
\label{Eq:DefLagIntRhs}
\Pi_r f(t)_{|I_n} = \sum_{j=0}^r f(t_{n,j}) \varphi_{n,j}(t) \quad \mbox{for}\; t\in
I_n\,.
\end{equation}

We note that the constants $\hat  \beta_{ii}$ are satisfying the following property.
%change coefficei
\begin{lemm}[{Coefficient property (C)}] 
\label{Lem:CoePropC}
There exist constants $\beta_m, 
\beta_M \in \R $ such that
\begin{equation}
\label{Eq:cond_on_beta}
0 < \beta_m \le  \hat  \beta_{ii} \le  \beta_M < \infty \,, \quad \mbox{for } \;
i=1,\ldots,r\,,
\end{equation}
is satisfied. The constants do not depend on the time step size, but only on the number 
$r$ of involved Gaussian quadrature points.
\end{lemm}

\begin{mproof}
Indeed, the coefficients $\beta_{ii}= \hat \omega_i$ are the Gauss-Legendre quadrature
weights
\begin{equation}\label{formula_weights}
\tilde w_i = \int_{-1}^1 \prod_{j=1 \atop j \neq
i}^r \left(\frac{x-x_j}{x_i-x_j}\right)^2 \ud x= \dfrac{1}{(1 - x_i^2)
(P_{r}^\prime(x_i))^2}\,, \quad i = 1,\ldots,r\,,
\end{equation}
scaled to the interval $[0,1]$, i.e.\ {\color{black} $\hat \omega_i = \tilde \omega_i/2$}. In
\eqref{formula_weights}, $P_r$ denotes the Legendre polynomial of degree $r$ and $x_i$,
for $i = 1,\ldots,r$, are its roots, cf., e.g., \cite[p.\ 436]{Quarteroni2007}. Since the
sum of the weights $\tilde \omega_i$ equals to two and the weights are all strictly
positive, we immediately conclude that an upper bound for $\hat \omega_i$ is given by one.
On the other hand, we know that $ | P_r^\prime(x) | \le {r(r+1)}/{2}$ for any $x \in
[-1,1]$; cf.\ \cite[p.\ 73]{Bernardi1992}. This gives us the lower bound
$\hat w_i \ge 2/(r(r+1))^2$.
\end{mproof}

%\begin{rem} \label{Rem:InitialFlux} The advantage of using the Gaussian quadrature rule for evaluating the time integrals in \eqref{Eq:locP_1}, \eqref{Eq:locP_2} is that in the resulting scheme \eqref{Eq:SemiDis1}, \eqref{Eq:SemiDis2} no flux approximation $\vec Q_n^0$ at the initial time point $t_{n-1}$ of the subinterval $I_n$ arises. Therefore, the initial flux at time $t=0$, that we usually do not have in practice, is not involved and needed for calculating the coefficient functions $U_n^j$ and $\vec Q_n^j$ for $j=1,\ldots, r$ at the intermediate Gaussian quadrature points $t_{n,l}\in (t_{n-1},t_n)$ for $l=1,\ldots ,r$. For evaluating \eqref{Eq:ContCond1} for $n\geq 2$, we calculate $u_\tau{}_{|I_{n-1}}(t_{n-1})$ explicitly by evaluating the fully discrete counterpart of the first of the representations in \eqref{Eq:RepBasis} on $I_{n-1}$. Approximations of the flux variable are computed similarly by evaluating the fully discrete counterpart of the second of the representations in \eqref{Eq:RepBasis} on $I_{n-1}$. For the initial value of the flux variable either a sufficiently smooth quantity $-\vec D \nabla \vec u_0$, if available, or an approximation computed by an extrapolation technique (cf.\ \cite{Hussain2012}) is used. \end{rem}

Below, we will also need the following auxiliary results.
\begin{lemm}
\label{lemma_timeintegration_basicfunc}
Let $ F(t, \vec x) = \sum_{i=0}^r F_n^i (\vec x) \varphi_{n, i} (t)$, for $t\in I_n$,
with coefficient functions $F_n^i \in W$ for $i=0,\ldots,r$. Then it holds that
\begin{equation}
\label{result_lemma_timeintegration_basicfunc}
\sum_{i = 1}^r  \sum_{j = 0}^r \hat \alpha_{ij} \langle F_{n}^j, F_{n}^i \rangle =
\int_{t_{n-1}}^{t_n} \langle \partial_t F, F \rangle dt = \dfrac{1}{2} \|  F (t_n)
\|^2  - \dfrac{1}{2} \|  F (t_{n-1}) \|^2
\end{equation}
and
\begin{equation}\label{result_lemma_L2_I_X_error}
\| F \|_{L^2(I_n; W)}^2 \le c\, \tau_n\sum_{j=0}^r \| F_{n}^j \|_W^2\,,
\end{equation}
for some $c >0$ independent of $\tau_n$. An analogous results holds for
coefficients $\vec F_n^i\in \vec V$.
\end{lemm}

\begin{mproof} Using the properties of the basis functions $\varphi_i$ and that the
$r$-point Gaussian quadrature formula is exact for polynomials of maximum degree $2r - 1$ there
holds that
\begin{align*}
\int_{t_{n-1}}^{t_n} \langle \partial_t F, F \rangle \ud t & = \int_{t_{n-1}}^{t_n}
\int_\Omega \sum_{j = 0}^r  \varphi_{n,j}^\prime (t) F_n^j (\vec x)   \sum_{i = 0}^r
\varphi_{n,i} (t) F_n^i (\vec x) \ud \vec x \ud t \\[1ex]
& =  \sum_{j = 0}^r \sum_{i = 0}^r \int_{0}^{1}  \frac{\ud}{\ud\hat t} \hat
\varphi_{j} (\hat t)  \cdot  \hat \varphi_i (\hat t) \ud \hat t \; \langle  F_n^i,  
F_n^j \rangle \\[1ex]
& =  \sum_{j = 0}^r \sum_{i = 1}^r \hat w_i \hat \varphi_j^\prime (\hat t_{i})  \langle
 F_n^i,  F_n^j \rangle =  \sum_{i = 1}^r  \sum_{j = 0}^r \hat \alpha_{ij} \langle
F_{n}^j, F_{n}^i \rangle.
\end{align*}
The second of the equalities in \eqref{result_lemma_timeintegration_basicfunc} follows
immediately from the first one. It remains to prove \eqref{result_lemma_L2_I_X_error}.
It holds that
\begin{align*}
\| F \|_{L^2(I_n; W)}^2  \le (r +1) \sum_{i=0}^r  
\int_{t_{n-1}}^{t_n} \varphi_{n, i}^2 (t) \ud t \; \|F_n^i\|^2 \le  c (r +1) \sum_{i=0}^r 
\tau_n \| F_n^i \|^2,
%\label{Eq:AuxIneq}
\end{align*}
with $c$ independent of $\tau_n$. Here we used that $\int_{t_{n-1}}^{t_n} \varphi_{n,
i}^2(t) \ud t \le c\, \tau_n$; cf.\ \cite[p.\ 1790]{Karakashian1999}.
\end{mproof}

\subsection{Discretization in space by the mixed finite element method}
\label{Sec:cGPrMFEM}

Now, we present the fully discrete approximation scheme that is obtained by discretizing \eqref{Eq:SemiDis1}, \eqref{Eq:SemiDis2} with respect to their spatial variables. For this we choose a pair of finite element spaces $W_h
\subset W$ and $\vec V_h\subset \vec V$ satisfying the inf-sup stability condition; cf.\ 
\cite{Brezzi1991,Chen2010}. Here, we denote by $\mathcal{T}_h=\{K\}$ a finite element 
decomposition of mesh
size $h$ of the polyhedral domain $\overline{\Omega}$ into closed subsets $K$,
quadrilaterals in two space dimensions and hexahedrals in three space dimensions. Since
the software library \texttt{deal.ii} \cite{DealIIReference} that we use for our
implementation of the schemes allows only quadrilateral and hexahedral elements, we
restrict ourselves to these types of elements in the following. Triangular and
tetrahedral elements can be treated in an analogous way. In our calculations (cf.\ Sec.\
\ref{Sec:NumStudies}) we use the Raviart--Thomas element on quadrilateral meshes for 
two space dimensions. For an application in three dimensions based on 
the Raviart--Thomas--N\'{e}d\'{e}lec element we refer to
\cite{Bause2015,Koecher2015}.

The construction of the discrete function spaces $W_h$ and $\vec V_h$ on quadrilateral 
and hexahedral finite elements is done by a transformation $\mathcal{T}_K: \hat K
\rightarrow K $ of the reference element $\hat K = [0,1]^d$, with $d=2$ or
$d=3$, to the element $K$ through a diffeomorphism $\mathcal{T}_K$
for all $K\in \mathcal{T}_h$. We sketch this briefly for $d=2$; cf.\ 
\cite{Chen2010,Koecher2015} for $d=3$. For this, let
\begin{align*}
\hat Q^{p_1,p_2} & : = \bigg\{\hat p: [0,1]^2 \rightarrow \R \;\Big|\; \hat p
(\vec{\hat x}) = \sum_{i=0}^{p_1} \sum_{j=0}^{p_2} p_{i,j} x_1^i x_2^j\,, \;
p_{i,j}\in \R \bigg\}\,.
\end{align*}
We then define the discrete subspaces $W_h^p \subset W$ and $\vec V_h^p \subset \vec V$
by
\begin{align}
W_h & = W_h^p := \left\{ w \in W \;\Big| \; w_K \circ \mathcal{T}_K^{-1} \in \hat
Q^{p,p}\,, \; \mathrm{for}\; K \in \mathcal{T}_h \right\}\,,
\label{Eq:DefWh}\\
\vec V_h & = \vec V_h^p := 
\left\{ \vec v \in \vec V \;\Big| \; \vec v_K \circ \mathcal{T}_K^{-1} \in \hat
Q^{p+1,p}\times \hat Q^{p,p+1}\,, \; \mathrm{for}\; K \in \mathcal{T}_h \right\}\,.
\label{Eq:DefVh}
\end{align}

The fully discrete continuous Galerkin-Petrov and MFE approximation
scheme, referred to as cGP($r$)--MFEM($p$), then defines fully discrete solutions $u_{\tau,h}\in \mathcal X^r(W_h)$ and $\vec q_{\tau,h}\in \mathcal X^r(\vec V_h)$ that are represent in terms of basis functions in time by 
\begin{equation*}
 u_{\tau,h}{}_{|\overline{I}_n} (t ) = \sum_{j=0}^r U_{n,h}^j\,  \varphi_{n,j}(t) \quad
\mathrm{and} \quad
 \vec q_{\tau,h}{}_{|\overline{I}_n} (t ) = \sum_{j=0}^r \vec Q_{n,h}^j \, \varphi_{n,j}(t)\,,
\quad
\mathrm{for} \; t \in I_n\,,
\end{equation*}
with coefficient functions $U_{n,h}^j\in W_h$ and $\vec Q_{n,h}^j \in \vec V_h$ for $j=0,\ldots, r$. The coefficient functions are obtained by solving the variational problem  \eqref{Eq:SemiDis1}, \eqref{Eq:SemiDis2} in the
 discrete subspaces $W_h\subset W$ and $\vec V_h\subset \vec V$:
\emph{For $n=1,\ldots, N$ and $j=1,\ldots , r$ find coefficient functions 
$\{U_{n,h}^j,\vec Q_{n,h}^j\} \in W_h\times \vec V_h$ such that}
\begin{align}
\sum_{j=0}^r \hat \alpha_{ij} \langle U_{n,h}^j ,w_h\rangle + {\tau_n} \, \hat
\beta_{ii}
\langle \nabla \cdot \vec Q_{n,h}^i ,w_h\rangle & = {\tau_n}\,
\hat \beta_{ii}\langle
f(t_{n,i}) ,w_h\rangle  \,,
\label{Eq:FulDis1}\\[0ex]
\langle \vec D^{-1} \vec Q_{n,h}^i,\vec v_h \rangle - \langle U_{n,h}^i ,\nabla \cdot
\vec v_h\rangle & = 0
\label{Eq:FulDis2}
\end{align}
\emph{for $i=1,\ldots, r$ and all $\{w_h,\vec v_h\}\in W_h\times \vec V_h$, where
$U_{n,h}^0 \in W_h$ and $\vec Q_{n,h}^0 \in \vec V_h$ are defined by means of the continuity constraint (cf.\ Rem.\ \ref{Rem:IniVal} ), i.e.\ }
\begin{align}
\nonumber
U_{n,h}^0 & := \sum_{j=0}^r U_{n-1,h}^j\, \varphi_{n-1,j}(t_{n-1}) \;\; \mathrm{if}\; n\geq
2 \,, &
U_{n,h}^0  & := P_h u_0 \;\; \mathrm{if}\; n=1\,,\\ 
\label{Eq:ContCond3_2}
{\color{black} \vec Q_{n,h}^0} & {\color{black} := \sum_{j=0}^r \vec Q_{n-1,h}^j\, \varphi_{n-1,j}(t_{n-1}) \;\; \mathrm{if}\; n\geq
2 \,,} & {\color{black} \vec Q_{n,h}^0} & {\color{black} := \vec P_h (-\vec D\nabla u_0) \;\; \mathrm{if}\; n=1\,,}
\end{align}
\emph{with $P_h:L^2(\Omega)\mapsto W_h$ and $\vec P_h:L^2(\Omega)\mapsto \vec V_h$ denoting the $L^2$ projections onto $W_h$ and $\vec V_h$, respectively.}

For the derivation of the algebraic formulation of the fully discrete variational problem
\eqref{Eq:FulDis1}, \eqref{Eq:FulDis2} we also refer to \cite{Bause2015,Koecher2014}. In
\cite{Bause2015,Koecher2014}, the iterative solution of the arising linear systems and
the construction of an efficient preconditioner is further addressed. For solving
the algebraic counterpart of Eqs.\ \eqref{Eq:FulDis1}, \eqref{Eq:FulDis2} we do not apply
an additional hybridization technique as it was done, for instance,
in \cite{Bause2008,Bause2013,Bause2004} and the references therein. We solve the algebraic
system by using a Schur complement technique. In \cite{Koecher2014} the efficiency of
the proposed iterative solver along with an adapted preconditioning technique is analyzed
numerically. In \cite{Bause2015,Koecher2014}, the approximation properties
of some families of space-time discretization schemes, including the
cGP($r$)--MFEM($p$) approach, in terms of convergence rates and their robustness are
studied by numerous numerical experiments. Test cases in three space dimensions
and with heterogeneous and strongly anisotropic material properties are also included.

\section{Existence and uniqueness of the semidiscrete approximation and error estimates}
\label{Sec:ExUniErr}

In this subsection we prove the existence and uniqueness of solutions to the semidiscrete
approximation scheme that is defined by \eqref{Eq:GtdP_1},
\eqref{Eq:GtdP_2} and its numerically integrated counterpart \eqref{Eq:SemiDis1},
\eqref{Eq:SemiDis2}, respectively. The time discretization error is also studied in 
this section. The spatial discretization error is analyzed in Sec.\ 
\ref{Sec:MFEMcGDis}.

\subsection{Existence and uniqueness of the semidiscrete approximation}

\begin{theorem}[Uniqueness of solutions]
\label{Thm:Uni}
Let the assumptions of Subsec.\ \ref{Sec:NotPrem_2} about $\Omega, u_0$ and $f$
be satisfied. Then the solution $\{u_\tau,\vec q_\tau\}\in
\mathcal{X}^{r}(W)\times
\mathcal{X}^{r}(\vec V)$ of the semidiscrete problem \eqref{Eq:GtdP_1},
\eqref{Eq:GtdP_2} is unique.
\end{theorem}

\begin{mproof}
Suppose that $\{u_{\tau,1},\vec q_{\tau,1}\}\in \mathcal{X}^{r}(W)\times
\mathcal{X}^{r}(\vec V)$ and $\{u_{\tau,2},\vec q_{\tau,2}\}\in
\mathcal{X}^{r}(W)\times \mathcal{X}^{r}(\vec V)$, respectively, satisfy the semidiscrete
problem \eqref{Eq:GtdP_1}, \eqref{Eq:GtdP_2} and let $u_\tau:=u_{\tau,1}-u_{\tau,2}$ and
$\vec q_\tau:= \vec q_{\tau,1}-\vec q_{\tau,2}$. Then, the tuple $\{u_\tau,\vec q_\tau\}$
satisfies \eqref{Eq:GtdP_1}, \eqref{Eq:GtdP_2} with $f\equiv 0$.
Choosing the test function $w_\tau=A^{-1} \partial_t u_\tau\in 
\mathcal{Y}^{r-1}(D(A))\subset \mathcal{Y}^{r-1}(W)$ in \eqref{Eq:GtdP_1} yields that
\begin{equation}
\label{Eq:Uni_1}
\int_0^T \langle \partial_t u_\tau, A^{-1} \partial_t u_\tau \rangle \ud t +
\int_0^T
\langle
\nabla \cdot \vec q_\tau , A^{-1} \partial_t u_\tau \rangle \ud t = 0 \,.
\end{equation}
From estimate \eqref{Eq:Prop_A} we get that
\begin{equation}
\label{Eq:Uni_2}
\int_0^T \langle \partial_t u_\tau, A^{-1} \partial_t u_\tau
\rangle \ud t \geq \frac{\alpha}{\beta^2}\,  \int_0^T
\|\partial_t u_\tau \|^2_{H^{-1}(\Omega)}\ud t \,.
\end{equation}
Using integration by parts in the second of the integrals in \eqref{Eq:Uni_1}
and recalling that $A^{-1} \partial_t u_\tau \in D(A)\subset H^1_0(\Omega)$,
Eq.\ \eqref{Eq:Uni_1} along with \eqref{Eq:Uni_2} yields that
\begin{equation}
\label{Eq:Uni_3}
0 \geq  \frac{\alpha}{\beta^2} \int_0^T\|\partial_t u_\tau \|^2_{H^{-1}(\Omega)}\ud t -
\int_0^T \langle {\color{black} \vec q_\tau} ,
\nabla A^{-1} \partial_t u_\tau \rangle \ud t \,.
\end{equation}

Next, by choosing the test function $\vec v_\tau= \vec D \nabla A^{-1}
\partial_t u_\tau \in \mathcal{Y}^{r-1}(\vec V)$ in {\color{black} Eq.\  \eqref{Eq:GtdP_2}} we find that
\begin{equation}
\label{Eq:Uni_4}
\int_0^T \langle \vec D^{-1} \vec q_\tau , \vec D \nabla A^{-1} \partial_t
u_\tau \rangle \ud t - \int_0^T
\langle u_\tau , \nabla \cdot (\vec D \nabla A^{-1} \partial_t
u_\tau) \rangle \ud t  = 0  \,.
\end{equation}
Since
\begin{equation*}
 \nabla \cdot (\vec D \nabla A^{-1} \partial_t u_\tau) = - A A^{-1}\partial_t
u_\tau = - \partial_t u_\tau
\end{equation*}
and $\vec D = \vec D^\top$ by assumption, it follows from \eqref{Eq:Uni_4} that
\begin{equation*}
\int_0^T \langle \vec q_\tau , \nabla A^{-1} \partial_t
u_\tau \rangle \ud t + \dfrac{1}{2}\int_0^T \dfrac{\ud}{\ud t} \|u_\tau \|^2 \ud t
= 0  \,.
\end{equation*}
Since $u_\tau (0) = u_{\tau,1} (0) - u_{\tau,2} (0) = 0$ it follows that
\begin{equation}
\label{Eq:Uni_5}
\int_0^T \langle \vec q_\tau , \nabla A^{-1} \partial_t
u_\tau \rangle \ud t + \dfrac{1}{2} \|u_\tau (T) \|^2 = 0  \,.
\end{equation}
Combing relations \eqref{Eq:Uni_3} and \eqref{Eq:Uni_5} shows that
\begin{equation}
\label{Eq:New1}
0 \geq c \int_0^T\|\partial_t u_\tau \|^2_{H^{-1}(\Omega)}\ud t + \dfrac{1}{2} \|u_\tau
(T)
\|^2
\end{equation}
and, therefore, $u_\tau =0$. This implies the uniqueness of solutions
$u_\tau$ to \eqref{Eq:GtdP_1}, \eqref{Eq:GtdP_2}.

To show the uniqueness of solutions $\vec q_\tau$ of \eqref{Eq:GtdP_1},
\eqref{Eq:GtdP_2}, we choose the test function $\vec v_\tau = \partial_t
\vec q_\tau \in \mathcal{Y}^{r-1}(\vec V)$. Recalling that $u_\tau=0$ by means of the 
uniqueness result \eqref{Eq:New1}
we obtain from Eq.\ \eqref{Eq:GtdP_2} that
\begin{equation*}
\label{Eq:Uni_6}
 \int_0^T \langle \vec D^{-1} \vec q_\tau, \partial_t \vec q_\tau\rangle \ud t
= 0\,.
\end{equation*}
From $\langle \vec D^{-1}\vec q_\tau , \partial_t \vec q_\tau
\rangle = \frac{1}{2} \frac{\ud}{\ud t}\langle \vec D^{-1} \vec q_\tau , \vec
q_\tau\rangle$ and $\vec q_\tau(0) = \vec 0 $ we conclude that
\begin{equation}
\label{Eq:Uni_7}
0 =\frac{1}{2} \|\vec D^{-1/2} \vec q_\tau (T)\|^2\,.
\end{equation}
Next, we choose $\vec v_\tau = \partial_t^2 \vec q_\tau\in \mathcal{Y}^{r-2}(\vec V)$, {\color{black} $\mathcal{Y}^{r-2}(\vec V)\subset \mathcal{Y}^{r-1}(\vec V)$ by definition,} in 
\eqref{Eq:GtdP_2},
recall that $u_\tau=0$ and use that
\begin{equation*}
 \dfrac{\ud}{\ud t} \langle \vec D^{-1} \vec q_\tau, \partial_t \vec q_\tau \rangle
= \langle \vec D^{-1} \partial_t \vec q_\tau, \partial_t \vec q_\tau \rangle +
\langle \vec D^{-1} \vec q_\tau, \partial_t^2 \vec q_\tau \rangle\,.
\end{equation*}
Together, this implies that
\begin{equation}
\label{Eq:Uni_8}
 0 = \int_0^T \dfrac{\ud}{\ud t} \langle \vec D^{-1} \vec q_\tau,
\partial_t \vec q_\tau \rangle \ud t - \int_0^T \|\vec D^{-1/2}
\partial_t \vec q_\tau \|^2 \ud t\,.
\end{equation}
Since $\vec q_\tau (0)=\vec q_{\tau,1}(0)-\vec q_{\tau,2}(0) = \vec 0$ and,
further, $\vec q_\tau (T)= \vec 0$ by means of \eqref{Eq:Uni_7}, it follows
from Eq.\ \eqref{Eq:Uni_8} along with property \eqref{Eq:Uni_7} that $\vec
q_\tau = 0$. The uniqueness of solutions to
the variational problem \eqref{Eq:GtdP_1}, \eqref{Eq:GtdP_2} is thus proved.
\end{mproof}

\begin{theorem}[Existence of solutions]
\label{Thm:Ex}
Let the assumptions of Subsec.\ \ref{Sec:NotPrem_2} about $\Omega, u_0, \vec D$
and $f$ be satisfied. Then the semidiscrete problem \eqref{Eq:GtdP_1},
\eqref{Eq:GtdP_2} admits a solution $\{u_\tau,\vec q_\tau\}\in
\mathcal{X}^{r}(W)\times \mathcal{X}^{r}(\vec V)$.
\end{theorem}

\begin{mproof}
To prove existence of solutions to problem \eqref{Eq:GtdP_1},
\eqref{Eq:GtdP_2}, we will use an equivalent conformal formulation, see \cite{Radu2004} for a similar approach. 

%Let $\widetilde u_\tau \in X^r(H^1_0(\Omega))$ denote the solution of the semidiscrete problem that results from the variational in time discretization of the weak formulation of \eqref{Eq:Diff}--\eqref{Eq:InVal}, i.e.\ without rewriting Eq.\\eqref{Eq:Diff} as a first order system of equations with respect to the spatial variables: 

\textit{Find $\widetilde u_\tau \in X^r (H^1_0(\Omega))$
such that $\widetilde u_\tau (0) = u_0$ and}
\begin{equation}
\label{Eq:Ex_1}
 \int_0^T \langle \partial_t \widetilde u_\tau, w_\tau \rangle \ud t + \int_0^T
a(\widetilde u_\tau,w_\tau) \ud t = \int_0^T \langle f,w_\tau\rangle \ud   t
\end{equation}
\textit{for all $w_\tau \in Y^{r-1}(H^1_0(\Omega))$.}

The existence and uniqueness of the semidiscrete approximation satisfying
\eqref{Eq:Ex_1} can be established. This is shown in the appendix of this work.
Then we define
\begin{equation}
\label{Eq:Ex_2}
u_\tau := \widetilde u_\tau \qquad \mbox{and} \qquad \vec q_\tau := - \vec D
\nabla \widetilde u_\tau\,.
\end{equation}
Obviously, it holds that $u_\tau \in X^r(W)$ since $H^1_0(\Omega)\subset W$.
Further, we have that $\partial_t \widetilde u_\tau \in L^2(I;H^1_0(\Omega))$ since on
each of the subintervals $I_n$, $n=1,\ldots, N$ the function $\widetilde u_\tau
\in X^r(H^1_0(\Omega))$ admits the representation
\begin{displaymath}
u_\tau{}_{|I_n} (t) =  \sum_{j=0}^r U_n^j  \varphi_{n,j}(t)\,, \qquad
\mbox{for}\;
t\in I_n\,,
\end{displaymath}
with coefficients $U_n^j \in H^1_0(\Omega)$ and polynomial basis functions
$\varphi_{n,j}\in \P_r(I_n;\R)$.

Next, we prove that $\vec q_\tau \in X^r(\vec V)$. Under the
assumption of Subsec.\ \ref{Sec:NotPrem_2} that  $f\in L^2(I;W)$ it
follows that
\begin{align*}
\int_0^T \langle - \vec q_\tau , \nabla w_\tau \rangle \ud t =  \int_0^T
\langle \vec D \nabla u_\tau , \nabla w_\tau \rangle \ud t & = \int_0^T \langle
f-\partial_t u_\tau ,w_\tau\rangle \ud   t =:\int_0^T \langle \widetilde f
,w_\tau\rangle \ud   t
\end{align*}
for all $w_\tau \in Y^{r-1} (C_0^\infty (\Omega))$ with $\widetilde f \in
L^2(I;L^2(\Omega))$. Thus, we have that
\begin{equation*}
 \int_0^T \langle - \vec q_\tau , \nabla w_\tau \rangle \ud t =\int_0^T \langle
\widetilde f ,w_\tau\rangle \ud   t \,.
\end{equation*}
Consequently, it holds that {\color{black} (cf.\ \cite[p.\ 18, Eq.\ (3.38)]{Brezzi1991})}
\begin{equation*}
 \int_0^T \langle \nabla \cdot \vec q_\tau, w_\tau \rangle \ud t
=\int_0^T \langle \widetilde f ,w_\tau\rangle \ud   t
\end{equation*}
for all $w_\tau \in Y^{r-1} (C_0^\infty (\Omega))$ in the sense of
distributions. Since $\widetilde f \in L^2(I;L^2(\Omega))$,  it follows that $\nabla\cdot
\vec q_\tau  \in L^2(I;L^2(\Omega))$ and, therefore, that $\vec q_\tau\in  L^2(I; \vec
V)$ is fulfilled. Finally, from the expansion in terms of polynomial basis functions
\begin{equation}
\label{Eq:Ex_3}
 \vec q_\tau(t) = - \sum_{j=0}^r \vec D \nabla U_n^j \, \varphi_{n,j}(t)\,,
\end{equation}
%with coefficients $\vec D \nabla U_n^j \in \vec H(\mathrm{div};\Omega))$
we conclude that $\vec q_\tau \in C([0,T];\vec V)$.

Eq.\ \eqref{Eq:Ex_1} then directly implies that the functions $u_\tau$ and
$\vec q_\tau$ defined in \eqref{Eq:Ex_2} satisfy the first equation of the variational
problem \eqref{Eq:GtdP_1}, \eqref{Eq:GtdP_2}. The second equation of the system
\eqref{Eq:GtdP_1}, \eqref{Eq:GtdP_2} then follows from the representation
\eqref{Eq:Ex_3} of the variable $\vec q_\tau$ by testing the identity \eqref{Eq:Ex_3} with
some function $\vec v_\tau \in \mathcal{Y}^{r-1}(\vec V)$ and applying the divergence
theorem of Gauss. Hence, the assertion of the theorem is proved.
\end{mproof}

As a corollary of the previous two theorems proving the existence of a unique
solution to the semidiscrete problem \eqref{Eq:GtdP_1}, \eqref{Eq:GtdP_2} we obtain an
inf-sup stability condition within our space-time framework. This result will play a
fundamental role in our error analyses. For this we need some further notation. Let
$\{u_\tau,\vec q_\tau\}\in \mathcal{X}^{r}(W)\times \mathcal{X}^{r}(\vec V)$ denote the
solution of the semidiscrete problem \eqref{Eq:GtdP_1}, \eqref{Eq:GtdP_2}. We split
$u_\tau$ as
\begin{equation}
\label{Eq:Split}
u_\tau(t) = u_0 + u^0_\tau(t) \qquad \text{with}\quad u^0_\tau \in
\mathcal{X}_0^{r}(W)\,.
\end{equation}
In terms of the tuple $\{u^0_\tau, \vec q_\tau\}$ of unknowns we recast the
existence and uniqueness result of Thm.\ \ref{Thm:Uni} and \ref{Thm:Ex} in the following
form.

\begin{cor}
\label{Thm:ExUniHom}
Let the assumptions of Subsec.\ \ref{Sec:NotPrem_2} about $\Omega, u_0, \vec D$
and $f$ be satisfied. Let $\{u_\tau,\vec q_\tau\}\in \mathcal{X}^{r}(W)\times
\mathcal{X}^{r}(\vec V)$ be the unique solution of the semidiscrete problem
\eqref{Eq:GtdP_1}, \eqref{Eq:GtdP_2} according to
Thm.\ \ref{Thm:Uni} and \ref{Thm:Ex}. Then, the tuple $\{u_\tau^0,\vec
q_\tau\}\in \mathcal{X}_0^{r}(W)\times \mathcal{X}^{r}(\vec V)$ with $u_\tau^0$ being
defined in \eqref{Eq:Split} is the unique solution of the following variational problem:
\emph{Find $\{u^0_\tau, \vec q_\tau\}\in \mathcal{X}_0^r (W)\times \mathcal{X}^r(\vec V)$
such that}
\begin{align}
\label{Eq:GtdP_1_0}
\int_0^T \langle \partial_t u^0_\tau  ,w_\tau  \rangle \ud t + \int_0^T \langle
\nabla \cdot \vec q_\tau, w_\tau \rangle \ud t & = \int_0^T \langle f,
w_\tau \rangle \ud t\,,\\[2ex]
\label{Eq:GtdP_2_0}
\int_0^T
\langle \vec D^{-1} \vec q_\tau , \vec v_\tau \rangle \ud t - \int_0^T
\langle u^0_\tau , \nabla \cdot \vec v_\tau \rangle \ud t & = \int_0^T \langle
u_0 , \nabla \cdot \vec v_\tau \rangle \ud t
\end{align}
\emph{for all $w_\tau \in \mathcal{Y}^{r-1}(W)$ and $\vec v \in
\mathcal{Y}^{r-1}(\vec V)$.}
\end{cor}

As a corollary we get the following inf-sup stability condition.
\begin{cor}
\label{Thm:InfSup}
Let the assumptions of Subsec.\ \ref{Sec:NotPrem_2} about $\Omega, u_0, \vec D$
and $f$ be satisfied. Then, there exists a constant $\gamma > 0$ such that
\begin{equation}
 \label{Eq:InfSupStab}
 \inf_{\{u^0_\tau,\vec q_\tau\}\in \mathcal{W}\backslash \{\vec 0\}}
\sup_{\{w_\tau,\vec v_\tau\}\in
\mathcal{V}\backslash \{\vec 0\}} \frac{a_\tau(\{u^0_\tau,\vec q_\tau \},\{w_\tau,\vec
v_\tau\})}{\|\{u^0_\tau,\vec q_\tau\}\|_{\mathcal{W}} \, \|\{w_\tau,\vec
v_\tau\}\|_{\mathcal{V}}} \geq \gamma > 0\,.
\end{equation}
\end{cor}

\begin{mproof}
The discrete problem \eqref{Eq:GtdP_1_0}, \eqref{Eq:GtdP_2_0} satisfies
the assumptions of the Banach-Ne\v{c}as-Babu\v{s}ka theorem \cite[p.\ 85]{Ern2010}.
Since the discrete problem \eqref{Eq:GtdP_1_0}, \eqref{Eq:GtdP_2_0} is
well-posed according to Corollary \ref{Thm:ExUniHom}, the
Banach-Ne\v{c}as-Babu\v{s}ka theorem implies the inf-sup stability
condition \eqref{Eq:InfSupStab}.
\end{mproof}

\subsection{Estimates for the error between the continuous and the semidiscrete solution}
\label{Sec:ErrExctFrm}

Now we shall show error estimates for the exact form \eqref{Eq:GtdP_1},
\eqref{Eq:GtdP_2} of the cGP($r$) approach applied to the mixed formulation
\eqref{Eq:IntMix_1}, \eqref{Eq:IntMix_2} of our parabolic model problem.

For this we assume that the following approximation property are satisfied.
There exist interpolation operators $I_\tau: H^1_0(I,W) \mapsto
\mathcal{X}_0^r(W)$, $J_\tau: L^2(I;\vec V)\mapsto \mathcal{X}^r(\vec V)$
such that for sufficiently smooth functions $u\in H^1(I;W)$ and $\vec
q \in L(I;\vec V)$ and all time intervals $I_n$, for $n=1,\ldots ,N$,
it holds that
\begin{align}
\label{Eq:Intpol1}
 \| u - I_\tau u\|_{L^2(I_n;W)} & \leq c\, \tau_n^{r+1}
\|\partial_t^{r+1} u\|_{L^2(I_n;W)}\,,\\[1ex]
\label{Eq:Intpol2}
 \| \partial_t(u - I_\tau u)\|_{L^2(I_n;W)} & \leq c\, \tau_n^{r}
\|\partial_t^{r+1} u\|_{L^2(I_n;W)}\,,\\[1ex]
\label{Eq:Intpol3}
 \| \vec q - \vec J_\tau \vec q\|_{L^2(I_n;\vec V)} & \leq c\, \tau_n^{r+1}
\|\partial_t^{r+1} \vec q\|_{L^2(I_n;\vec V)}
\end{align}
with some constant $c$ independent of $\tau_n$ and $\tau$. The existence of such 
approximations is obviously ensured, for instance, by using Lagrange 
interpolation \cite{Thomee2006}.

We get the following error estimates in the natural norm of the time discretization.

\begin{theorem}[Space-time error estimate for exact form of cGP($r$)]
\label{Th:ErrNatNorm}
Let the assumptions of Subsec.\ \ref{Sec:NotPrem_2} about $\Omega, u_0, \vec D$
and $f$ be satisfied. Let $\{u,\vec q\} \in H^1(I;W)\times L^2(I;\vec V)$
denote the unique solution of the mixed problem \eqref{Eq:IntMix_1},
\eqref{Eq:IntMix_2} that is supposed to be sufficiently regular. Then the solution
$\{u_\tau,\vec q_\tau\}\in \mathcal{X}^{r}(W)\times \mathcal{X}^{r}(\vec V)$ of the
semidiscrete problem \eqref{Eq:GtdP_1}, \eqref{Eq:GtdP_2} satisfies the error estimate
\begin{align*}
\|\{u - u_\tau,\vec q - \vec q_\tau\}\|_{\mathcal{W}} & \leq c
\left\{\sum_{n=1}^N \tau_n^{2r} \Big(\|\partial_t^{r+1} u\|^2_{L^2(I_n;W)} +
\|\partial_t^{r+1} \vec q \|^2_{L^2(I_n;\vec V)}\Big) \right\}^{1/2}\\[2ex]
& \leq c \tau^{r} \Big(\|\partial_t^{r+1} u\|_{L^2(I;W)}+\|\partial_t^{r+1} \vec q
\|_{L^2(I;\vec V)}\Big) \,,
\end{align*}
where the constant $c$ is independent of $\tau_n$, $\tau$ and $T$.
\end{theorem}

\begin{mproof} By splitting
\begin{equation}
 \label{Eq:Id_Err_u}
u(t) = u_0 + u^0(t) \qquad \mathrm{with}\quad u^0\in H^1_0(I;W)
\end{equation}
and recalling the semidiscrete counterpart \eqref{Eq:Split}, we get that
\[
u(t) - u_\tau (t) = u^0(t) - u_\tau^0(t)\,, \qquad \partial^r_t u(t)
= \partial_t^r u^0(t)
\]
for almost every $t\in (0,T)$, such that it is sufficient to derive the asserted error
bounds of the theorem for $u^0-u_\tau^0$ instead of estimating $u-u_\tau$. This will be
done in the following.

By  \eqref{Eq:Intpol1} to \eqref{Eq:Intpol3} it holds that
\begin{align}
 \nonumber
\|\{u^0- I_\tau u^0,\vec q- \vec J_\tau \vec q\}\|_{\mathcal{W}} & \leq c
\left\{\sum_{n=1}^N \tau_n^{2r} \Big(\|\partial_t^{r+1} u^0\|^2_{L^2(I_n;W)} +
\|\partial_t^{r+1} \vec q \|^2_{L^2(I_n;\vec V)}\Big) \right\}^{1/2}\\[2ex]
\label{Eq:Intpol}
&  \leq c \tau^{r} \Big(\|\partial_t^{r+1} u^0\|_{L^2(I;W)}+\|\partial_t^{r+1} \vec q
\|_{L^2(I;\vec V)}\Big) \,.
\end{align}
For the discrete functions $w_\tau := u^0_\tau -
I_\tau u^0\in \mathcal{X}_0^r(W)$, $\vec v_\tau = \vec
q_\tau - \vec J_\tau \vec q\in \mathcal{X}^r(\vec V)$ there exist, due to the
inf-sup stability
condition \eqref{Eq:InfSupStab}, functions $\varphi_\tau\in
\mathcal{X}_0^r(W)$, $\vec \psi_\tau \in\mathcal{X}^r(\vec V)$  such that
\begin{align}
 \gamma \| \{w_\tau,\vec v_\tau \} \|_{\mathcal{W}} \| \{\varphi_\tau {\color{black} ,}\vec
\psi_\tau\} \|_{\mathcal{V}}  & \leq a_\tau(\{w_\tau,\vec
v_\tau\},\{\varphi_\tau,\vec \psi_\tau\}) \nonumber\\[1ex]
& = a_\tau(\{u^0 - I_\tau u^0,\vec q -  \vec J_\tau \vec q\},\{\varphi_\tau,\vec
\psi_\tau\}) \nonumber\\[1ex]
& \leq c \|\{u^0 - I_\tau u^0,\vec q -  \vec J_\tau \vec q\}
\|_{\mathcal{W}}
\| \{\varphi_\tau,\vec \psi_\tau\} \|_{\mathcal{V}} \,,
\label{Eq:Err_infsup}
\end{align}
where the Galerkin orthogonalities
\begin{align*}
\int_0^T \langle \partial_t (u^0_\tau -u^0)  ,w_\tau  \rangle \ud t +
\int_0^T \langle
\nabla \cdot (\vec q_\tau - \vec q), w_\tau \rangle \ud t & = 0\,,\\[2ex]
\int_0^T \langle \vec D^{-1} (\vec q_\tau - \vec q) , \vec v_\tau \rangle
\ud t - \int_0^T
\langle u^0_\tau - u^0 , \nabla \cdot \vec v_\tau \rangle \ud t & = 0
\end{align*}
have been used. From \eqref{Eq:Err_infsup} along with \eqref{Eq:Intpol}, we find that
\begin{align}
\nonumber
\|\{u^0_\tau - I_\tau u^0,\vec q_\tau & - \vec J_\tau \vec q\}\|_{\mathcal{W}} 
\leq {c}{\gamma}^{-1}
\|\{u^0- I_\tau u^0,\vec q- \vec J_\tau \vec q\}\|_{\mathcal{W}}\\[1ex]
& \leq {c }{\gamma}^{-1} \, \tau^{r} \, \Big(\|\partial_t^{r+1}
u\|_{L^2(I;W)}+\|\partial_t^{r+1} \vec q
\|_{L^2(I;\vec V)}\Big) \,.
\label{Eq:Est_Err}
\end{align}
From inequality \eqref{Eq:Est_Err} along with the interpolation error estimate
\eqref{Eq:Intpol} we conclude the assertion of the theorem by means of the triangle
inequality.
\end{mproof}

Thm.\ \ref{Th:ErrNatNorm} yields an error estimate with respect to the
natural space-time norm of the discretization scheme. The estimate is sharp
with respect to the contribution of $\|\partial_t (u- u_\tau)\|_{L^2(I;W)}$ to
the overall norm \eqref{Eq:SpcTmNrm}. However, the
estimate is suboptimal with respect to $\| u- u_\tau\|_{L^2(I;W)}$. In the following 
theorem, we sharpen our analysis
by providing an optimal order error estimate also for $\| u-
u_\tau\|_{L^2(I;W)}$. This is done by a duality argument. For this, the following
additional regularity assumption is needed.

\textbf{Regularity condition (R${}_{\mathbf{mix}}$).}
\emph{Suppose that $g\in L^2(I;W)$. The variational problem, find $z\in
H^1(I;W)$ and $\vec p \in L^2(I;\vec V)$ with $z(T)=0$ such that}
\begin{align}
\label{Eq:DualProb01}
\int_0^T \big(\langle - \partial_t z , w \rangle + \langle \nabla \cdot \vec p , w
\rangle\big) \ud t & = \int_0^T \langle g,w\rangle \ud t\,, \\[2ex]
\label{Eq:DualProb02}
\int_0^T \big( \langle \vec D^{-1} \vec p,\vec v \rangle - \langle z, \nabla
\cdot \vec v\rangle\big) \ud t & = 0
\end{align}
\emph{for all $w\in L^2(T,0;W)$, $\vec v \in L^2(T,0;\vec V)$ admits a unique
solution $\{z,\vec p\} \in H^1(I;L^2(\Omega))$ $\times L^2(I;\vec V)$ with the improved 
regularity $\vec p \in H^1(I;\vec V)$ and the a priori estimate}
\begin{equation}
\label{Eq:AssmpRegMix}
\| \partial_t \vec p \|_{L^2(I;\vec V)} \leq c \| g\|_{L^2(I;W)}\,.
\end{equation}

Formally, the corresponding strong form of \eqref{Eq:DualProb01}, \eqref{Eq:DualProb02} is
given by
\begin{equation}
\label{Eq:DualMix}
- \partial_t z + \nabla \cdot \vec p = g \,, \qquad \vec D^{-1} \vec p + \nabla z =
0 \quad \mbox{ in } \Omega \times
(0,T) \,,
\end{equation}
with $z(T)=0$ and homogeneous Dirichlet boundary conditions, that is obtained by rewriting
the dual problem associated with \eqref{Eq:Diff}--\eqref{Eq:InVal},
\begin{equation}
\label{Eq:Dual2ndOrd}
-\partial_t z - \nabla \cdot (\vec D\nabla z) = g\; \mbox{ in } \Omega \times
(0,T)  \,, \quad z(T) = 0\; \mbox{ in } \Omega\,, \quad z = 0 \; \mbox{ on } \partial
\Omega \times (0,T)\,,
\end{equation}
as a system of first order equations. Defining the transformation $\widetilde z(t):=z(T-t)$ and $\widetilde g(t):=z(T-t)$ we recast \eqref{Eq:Dual2ndOrd} as a forward parabolic
problem in $\widetilde z$,
{\color{black}
\[
\partial_t \widetilde z - \nabla \cdot (\vec D\nabla \widetilde z) = \widetilde g\; \mbox{ in } \Omega \times
I  \,, \quad \widetilde z(0) = 0\; \mbox{ in } \Omega\,, \quad \widetilde z = 0 \; \mbox{ on } \partial
\Omega \times (0,T)\,,
\]
}
such that standard existence and stability estimates
can be applied; cf.\  \cite[p.\ 382, Theorem 5]{Evans2010}. Then, defining the variable
$\vec p$ by means of the second of the identities in \eqref{Eq:DualMix}, the thus
obtained tuple $\{z,\vec p\}$ satisfies the variational problem \eqref{Eq:DualProb01},
\eqref{Eq:DualProb02}. Moreover, for $g \in L^2(I;W)$, from \cite[p.\ 382,
Theorem 5]{Evans2010} we get the a priori estimate
\begin{equation}
\label{Eq:Stab1}
\|\partial_t z \|_{L^2(I;W)} + \| \vec p \|_{L^2(I;\vec V)} \leq c \| g
\|_{L^2(I;W)}\,.
\end{equation}
For this we note that $\vec p \in L^2(I;\vec V)$ can be shown by using the arguments
of the proof of Thm.\ \ref{Thm:Ex}. The a priori estimate of the vector variable $\vec
p$ in \eqref{Eq:Stab1} is then a direct consequence of the variational equation
\eqref{Eq:DualProb01}.

\begin{rem}
A regularity condition similar to (R${}_{\mathrm{mix}}$) is also used in \cite[p.\ 48, 
Eq.\ (6.16)]{Schieweck2010} to prove the optimal order convergence of a variational time 
discretization of second order parabolic problems in the non-mixed formulation. Currently, 
it remains an open problem how this limiting condition can be avoided in the theoretical 
analysis. The techniques that were developed recently in \cite{Ern2016} might be helpful. 
However, in our numerical convergence studies of Sec.\ \ref{Sec:NumStudies} the optimal 
convergence rate that is proved in Thm.\ \ref{Th:ErrL2NormMix} under the condition 
(R${}_{\mathrm{mix}}$) is nicely observed. 
\end{rem}

Below we also need the following auxiliary lemma. 

\begin{lemm}
Let $I_0:H^1(I,W)\mapsto \mathcal{Y}^0(W)$ and $\vec J_0:H^1(I,\vec V)\mapsto 
\mathcal{Y}^0(\vec V)$ be interpolation operators that are defined on each subinterval 
$I_n$ be means of
\[
I_0 u(t) := u(t_{n-1}) \quad and \quad \vec J_0 \vec v(t) := \vec v(t_{n-1}) 
\quad \text{for all}\;\;  t \in I_n\,.
\]
Then it holds that
\begin{align}
\label{Eq:Itpz}
\|z-I_0 z\|_{L^2(I_n,W)} & \leq \tau_n \| \partial_t z \|_{L^2(I_n;W)}\,,\\[1ex]
\label{Eq:Itpp}
\|\vec p -\vec J_0 \vec p \|_{L^2(I_n;\vec V)} & \leq \tau_n \| \partial_t \vec p
\|_{L^2(I_n;\vec V)}\, .
\end{align}
\end{lemm}

\begin{mproof}
The assertions directly follow from \cite[Lemma 6.2]{Schieweck2010}; cf.\ also 
\cite[p.\ 49]{Schieweck2010}.
\end{mproof}

\begin{theorem}[$L^2$ Error estimate for the exact form of cGP($r$)]
\label{Th:ErrL2NormMix}
Let the assumptions of Subsec.\ \ref{Sec:NotPrem_2} about $\Omega, u_0, \vec D$ and $f$ be
satisfied. Further, suppose that the regularity condition (R${}_{\mathrm{mix}}$) holds.
Let $\{u,\vec q\} \in H^1(I;W)\times L^2(I;\vec V)$ denote the unique solution of the
mixed problem \eqref{Eq:IntMix_1}, \eqref{Eq:IntMix_2} that is supposed to be sufficiently
regular. Then the solution $\{u_\tau,\vec q_\tau\}\in \mathcal{X}^{r}(W)\times
\mathcal{X}^{r}(\vec V)$ of the semidiscrete problem \eqref{Eq:GtdP_1}, \eqref{Eq:GtdP_2}
satisfies the error estimate
\begin{align*}
\| u - u_\tau \|_{L^2(I;W)} & \leq c \, \tau
\left\{\sum_{n=1}^N \tau_n^{2r} \Big(\|\partial_t^{r+1} u\|^2_{L^2(I_n;W)} +
\|\partial_t^{r+1} \vec q \|^2_{L^2(I_n;\vec V)}\Big) \right\}^{1/2}\\[2ex]
& \leq c \, \tau^{r+1} \Big(\|\partial_t^{r+1}
u\|_{L^2(I;W)}+\|\partial_t^{r+1} \vec q
\|_{L^2(I;\vec V)}\Big) \,.
\end{align*}
\end{theorem}

\begin{mproof}
We put $e_u:=u^0 -u^0_\tau\in L^2(I;W)$ and $\vec e_{\vec q}:= \vec q - \vec q_\tau\in 
L^2(I;\vec V)$ with the splitting \eqref{Eq:Id_Err_u} and \eqref{Eq:Split} of the scalar 
variable and its semidiscrete approximation, respectively. Further, let $\{z,\vec p\}\in 
H^1(0,T;W)\cap C([0,T];W) \times L^2(0,T;\vec V)$ with $z(T)=0$ denote the unique 
solution of \eqref{Eq:DualProb01}, \eqref{Eq:DualProb02} with right-hand side function 
$g=e_u$. 

Firstly, recalling that $z(T)=0$ and $e_u(0)=0$ by definition, we get that
\begin{equation}
\label{Eq:L2err_1}
\int_0^T \langle - \partial_t z , e_u \rangle \ud t = -{z(T)}e_u(T)+
z(0){e_u(0)} + \int_0^T \langle \partial_t e_u ,z\rangle \ud t= \int_0^T \langle
\partial_t e_u ,z\rangle \ud t\,.
\end{equation}
Choosing the test function $w=e_u$ in \eqref{Eq:DualProb01} and using \eqref{Eq:L2err_1},
we find that
\begin{align}
\nonumber
&\int_0^T \|e_u \|^2 \ud t  = \int_0^T \big(\langle \partial_t e_u , z\rangle +
\langle \nabla \cdot \vec p, e_u \rangle \big)\ud t \\[2ex]
& = \int_0^T \big(\langle \partial_t e_u , z\rangle +
\langle \nabla \cdot \vec e_{\vec q}, z \rangle \big)\ud t
-\int_0^T \big(\langle \nabla \cdot \vec e_q, z \rangle - \langle
\nabla \cdot \vec p, e_u \rangle \big)\ud t\,.
\label{Eq:L2err_1_1}
\end{align}
Choosing the test function $\vec v=\vec e_{\vec q}$ in \eqref{Eq:DualProb02} and recalling
that the matrix $\vec D$ is symmetric by assumption, we conclude that
\begin{equation}
 \label{Eq:L2err_2}
\int_0^T \langle \nabla \cdot \vec e_{\vec q}, z \rangle\ud t =  \int_0^T
\langle \vec D^{-1} \vec p, \vec e_{\vec q}\rangle \ud t = \int_0^T \langle \vec
D^{-1} \vec e_{\vec q}, \vec p \rangle \ud t\,.
\end{equation}
From \eqref{Eq:L2err_1_1} and \eqref{Eq:L2err_2} it then follows that
\begin{equation}
\label{Eq:L2Err1}
 \int_0^T \|e_u \|^2 \ud t  =\int_0^T \big(\langle \partial_t e_u , z\rangle +
\langle \nabla \cdot \vec e_q, z \rangle \big)\ud t  - \int_0^T \big(\langle
\vec D^{-1} \vec e_{\vec q}, \vec p \rangle - \langle e_u, \nabla \cdot \vec p
\rangle \big)\ud t\,.
\end{equation}

Secondly, by Galerkin orthogonality we find that
\begin{align}
\label{Eq:GalOrth1}
  \int_0^T \big(\langle \partial_t e_u , w_\tau \rangle +
\langle \nabla \cdot \vec e_{\vec q}, w_\tau \rangle \big)\ud t & = 0\,, \\[2ex]
\label{Eq:GalOrth2}
\int_0^T \big( \langle \vec D^{-1} \vec e_{\vec q},\vec v_\tau \rangle -
\langle e_u, \nabla
\cdot \vec v_\tau \rangle\big) \ud t & = 0
\end{align}
for all $w_\tau \in \mathcal{Y}^{r-1}(W)$ and $\vec v_\tau \in
\mathcal{Y}^{r-1}(\vec V)$. Choosing $w_\tau = I_0 z$ in \eqref{Eq:GalOrth1},
it follows that
\begin{equation}
\label{Eq:GalOrth3}
 \int_0^T \big(\langle \partial_t e_u , I_0 z \rangle +  \langle \nabla \cdot
\vec e_{\vec q}, I_0 z \rangle \big)\ud t  = 0\,.
\end{equation}
Further, choosing $\vec v_\tau = \vec J_0 \vec p$ in \eqref{Eq:GalOrth2} yields that
\begin{equation}
\label{Eq:GalOrth4}
\int_0^T \big( \langle \vec D^{-1} \vec e_{\vec q},\vec J_0 \vec p
\rangle - \langle e_u, \nabla
\cdot \vec J_0 \vec p_\tau \rangle\big) \ud t  = 0 \,.
\end{equation}

Thirdly, combining \eqref{Eq:L2Err1} with \eqref{Eq:GalOrth3} and 
\eqref{Eq:GalOrth4}, and then
using the Cauchy--Schwarz inequality as well as the interpolation error
estimates \eqref{Eq:Itpz} and \eqref{Eq:Itpp} yields that
\begin{align*}
\|e_u \|^2_{L^2(I;W)} & = \int_0^T \big(\langle \partial_t e_u , z-I_0 z\rangle +
\langle \nabla \cdot \vec e_q, z-I_0 z \rangle \big)\ud t \\[1ex]
& \quad - \int_0^T
\big(\langle
\vec D^{-1} \vec e_{\vec q}, \vec p - \vec J_0 \vec p \rangle - \langle e_u,
\nabla \cdot (\vec p - \vec J_0 \vec p)  \rangle \big)\ud t\,.\\[2ex]
& \leq  \big(\|\partial_t e_u\|_{L^2(I;W)} + \| \vec e_{\vec q}\|_{L^2(I;\vec V)}
\big)\|z-I_0 z\|_{L^2(I;W)}\\[2ex]
& \quad + \big(\|\vec D^{-1} \|_2 \|\vec e_{\vec q}\|_{_{L^2(I;\vec V)}}+ \|e_u
\|_{L^2(I;W)}\big) \| \vec
p - \vec J_0 \vec p)\|_{L^2(I;\vec V)}\\[2ex]
& \leq \tau \big(\theta_M \|\partial_t e_u\|_{L^2(I;W)} + \| \vec e_{\vec q}\|_{L^2(I;\vec
V)}\big) \|\partial_t z \|_{L^2(I;W)} \\[2ex]
& \quad +  c \tau \big(\| \vec e_{\vec q}\|_{L^2(I;\vec V)} + \|e_u
\|_{L^2(I;W)}\big)
\|\partial_t \vec p \|_{L^2(I;\vec V)}\,.
\end{align*}
Applying the a priori estimate \eqref{Eq:Stab1} and the additional regularity assumption
\eqref{Eq:AssmpRegMix} with $g=e_u$ as well as using the error estimate of Thm.\ 
\ref{Th:ErrNatNorm}, we then find that
\begin{align*}
\|e_u \|_{L^2(I;W)} & \leq c\, \tau \,
\left\{\sum_{n=1}^N \tau_n^{2r} \Big(\|\partial_t^{r+1} u\|^2_{L^2(I_n;W)} +
\|\partial_t^{r+1} \vec q \|^2_{L^2(I_n;\vec V)}\Big) \right\}^{1/2}\\[2ex]
& \leq c\, \tau^{r+1} \Big(\|\partial_t^{r+1}
u\|_{L^2(I;W)}+\|\partial_t^{r+1} \vec q
\|_{L^2(I;\vec V)}\Big) \,.
\end{align*}
This proves the assertion of the theorem.

%Since it holds that $e_u = u-u_\tau\in H^1(I;W)$ and $z(T)=0$, it follows by an
%regularity result for parabolic problems (cf.\ \cite[p.\ 382, Theorem 5]{Evans2010})
%that $z\in L^\infty(I;H^2(\Omega))$ and $z'\in L^2(I;L^2(\Omega))$, if a boundary
%$\partial \Omega$ of class $C^{1,1}$ and coefficients $\vec d_{ij} \in
%W^{1,\infty}(\Omega)$, for $i,j=1,\ldots,d$, are assumed.
\end{mproof}

%\subsection{Error estimate for the non-exakt form of cGP($r$)}
%\label{Sec:ErrNonExaCGP}
%In the previous subsections we studied the exact form of the cGP($r$) approach given
%by the system \eqref{Eq:GtdP_1}, \eqref{Eq:GtdP_2}. In

Next  we derive an error  estimate for the non-exact form \eqref{Eq:SemiDis1},
\eqref{Eq:SemiDis2} of the cGP($r$) method. The difference of the non-exact form of
cGP($r$) to \eqref{Eq:GtdP_1}, \eqref{Eq:GtdP_2} comes through the numerically
integrated right-hand side term in \eqref{Eq:SemiDis1}. Firstly, we ensure the existence
and uniqueness of the solution to the non-exact form of cGP($r$).

\begin{theorem}[Existence and uniqueness]
Let the assumptions of Subsec.\ \ref{Sec:NotPrem_2} about $\Omega, u_0, \vec D$
and $f$ be satisfied. Then the non-exact form \eqref{Eq:SemiDis1},
\eqref{Eq:SemiDis2} of the semidiscrete problem admits a unique solution $\{U_n^j,\vec
Q_n^j\}\in W\times \vec V$ for $j=1,\ldots,r$ and $n=1,\ldots,N$ defining semidiscrete
approximations $\{u_\tau,\vec q_\tau\}\in \mathcal{X}^{r}(W)\times \mathcal{X}^{r}(\vec
V)$ by means of the expansions \eqref{Eq:RepBasis} and the initial condition 
$u_\tau(0)=u_0$.
\end{theorem}

\begin{mproof}
By the definition of the Lagrange interpolation operator $\Pi_r$ given in
\eqref{Eq:DefLagIntRhs}
and the representations \eqref{Eq:RepBasis} of $u_\tau$ and $\vec q_\tau$ in terms of
basis functions we recast the non-exact form \eqref{Eq:SemiDis1}, \eqref{Eq:SemiDis2} of
the semidiscrete problem in the equivalent form
\begin{align}
\label{Eq:NonExGtdP_1}
\int_0^T \langle \partial_t u_\tau  ,w_\tau  \rangle \ud t + \int_0^T \langle
\nabla \cdot \vec q_\tau, w_\tau \rangle \ud t & = \int_0^T \langle \Pi_r f,
w_\tau \rangle \ud t\\[2ex]
\label{Eq:NonExGtdP_2}
\int_0^T
\langle \vec D^{-1} \vec q_\tau , \vec v_\tau \rangle \ud t - \int_0^T
\langle u_\tau , \nabla \cdot \vec v_\tau \rangle \ud t & = 0
\end{align}
for all $w_\tau \in \mathcal{Y}^{r-1}(W)$ and $\vec v \in
\mathcal{Y}^{r-1}(\vec V)$ with the initial condition $u_\tau (0) = u_0$.

Existence and uniqueness of the solution $\{u_\tau,\vec q_\tau\}\in
\mathcal{X}^r(W)\times \mathcal{X}^r(\vec V)$ to the system \eqref{Eq:NonExGtdP_1},
\eqref{Eq:NonExGtdP_2} then follows as in Thm.\ \ref{Thm:Uni} and \ref{Thm:Ex} with
$\Pi_r f$ replacing $f$ in the arguments of the proofs.
\end{mproof}

Next, we present the corresponding a priori error estimate.
\begin{theorem}
\label{Thm:ErrNonExNat}
Let the assumptions of Subsec.\ \ref{Sec:NotPrem_2} about $\Omega, u_0, \vec D$
and $f$ be satisfied. Suppose that $f$ is sufficiently regular with respect to the time
variable. Let $\{u,\vec q\} \in H^1(I;W)\times L^2(I;\vec V)$
denote the unique solution of the mixed problem \eqref{Eq:IntMix_1},
\eqref{Eq:IntMix_2} that is supposed to be sufficiently regular. Then the solution
$\{u_\tau,\vec q_\tau\}\in \mathcal{X}^{r}(W)\times \mathcal{X}^{r}(\vec V)$ of the
non-exact semidiscrete problem \eqref{Eq:GtdP_1}, \eqref{Eq:GtdP_2} satisfies
the error estimate
\begin{equation*}
\begin{split}
\|\{u - u_\tau,\vec q & - \vec
q_\tau\}\|_{\mathcal{W}} \leq c \Bigg(\sum_{n=1}^N
 \tau_n^{2r} \Big\{\|\partial_t^{r+1} u\|^2_{L^2(I_n;W)} + \tau_n^2 \|\partial_t^{r+1}
\vec q
\|^2_{L^2(I_n;\vec V)}\\[0ex]
& \hspace*{20ex} + \tau_n^2 \|\partial_t^{r+1} f
\|^2_{L^2(I_n;W)}\Big\}\Bigg)^{1/2}\\[1ex]
& \leq c \, \tau^r \Big\{\| \partial_t^{r+1} u|\|_{L^2(I;W)} +\tau \| \partial_t^{r+1}
\vec q|\|_{L^2(I;W)} +  \tau \|\partial_t^{r+1} f\|_{L^2(I;W)} \Big\}\,.
\end{split}
\end{equation*}
where the constants $c$ is independent of $\tau_n$, $\tau$ and $T$.
\end{theorem}

Since the proof of Thm.\ \ref{Thm:ErrNonExNat} follows from the proof of
Thm.\ \ref{Th:ErrNatNorm} by a standard estimate of the interpolation error,
we skip it here. For the sake of completeness we summarize the proof in the
appendix of this work.

\section{Existence and uniqueness of the fully discrete approximation and error estimates}
\label{Sec:MFEMcGDis}

In the first subsection of Sec.\ \ref{Sec:MFEMcGDis} we prove the existence and uniqueness
of solutions to the fully discrete approximation scheme \eqref{Eq:FulDis1},
\eqref{Eq:FulDis2}. Then, in Subsec.\ \ref{Sec:ErrMFEM} we establish an estimate for the
error between the non-exact form of the semidiscrete approximation defined by Eqs.\
\eqref{Eq:SemiDis1}, \eqref{Eq:SemiDis2} and the fully discrete solution given
by Eqs.\ \eqref{Eq:FulDis1}, \eqref{Eq:FulDis2}. Finally, in Subsec.\ \ref{Sec:ErrCGMFEM}
we combine the error estimates of the temporal discretization that are derived in
Sec.~\ref{Sec:ExUniErr} with the error estimates of Subsec.\ \ref{Sec:ErrMFEM} to get the
desired error estimates.

% Eq:FulDisAppFun

\subsection{Existence and uniqueness of the fully discrete approximation}
\label{Sec:ExUniFulDis}

Firstly we prove the existence and uniqueness of solutions to the fully discrete
cGP($r$)--MFEM($p$) scheme \eqref{Eq:FulDis1}, \eqref{Eq:FulDis2}. For this we
need the following lemma {\color{black} (cf.\ \cite[p.\ 302]{Radu2008}}.

\begin{lemm}
\label{lemma_thomas}
For given $ w_h \in W_h$ there exits a function $ \vec v_h \in \vec V_h$  satisfying
\begin{displaymath}
\nabla \cdot \vec v_h =  w_h  \qquad { \text{ and } }\qquad   \| \vec v_h \| \le
c  \| w_h \|
\end{displaymath}
for some constant $c > 0$ depending on $\Omega$ and the space dimension $d$
but not on $w_h$ or the mesh size $h$.
\end{lemm}

\begin{theorem}[Existence and uniqueness of solutions]
\label{Eq:ExUniFullDisc}
Let the assumptions of Subsec.\ \ref{Sec:NotPrem_2} about $\Omega, u_0, \vec D$
and $f$ be satisfied. Then the fully discrete problem \eqref{Eq:FulDis1}, 
\eqref{Eq:FulDis2}  admits a unique solution $\{u_{\tau,h},\vec
q_{\tau,h}\}\in \mathcal{X}^{r}(W_h)\times \mathcal{X}^{r}(\vec V_h)$.
\end{theorem}
%\eqref{Eq:locPdis_1}, \eqref{Eq:locPdis_2}, with its equivalent algebraic form \Eqref{Eq:Fuldis1}, \Eqref{Eq:Fuldis2}, 
\begin{mproof}
Since the fully discrete problem %\eqref{Eq:locPdis_1}, \eqref{Eq:locPdis_2} or
 \eqref{Eq:FulDis1}, \eqref{Eq:FulDis2} 
is finite dimensional and linear, it is sufficient to show
the uniqueness of the solution. The existence is then a direct consequence. Assume that
there exist two pairs of solutions $\{u_{\tau,h}^k,\vec q_{\tau,h}^k\}\in
\mathcal{X}^{r}(W_h)\times \mathcal{X}^{r}(\vec V_h)$, for $k=1,2$, that are represented  
in terms of basis functions by
\begin{equation*}
 \label{Eq:RepFulDis}
u_{\tau,h}^{k} (t)_{|I_n} = \sum_{j=0}^r U_{n,h}^{j,k} \varphi_{n,j}(t) \quad \mbox{and}
\quad \vec q_{\tau,h}^{k} (t)_{|I_n} =  \sum_{j=0}^r \vec
Q_{n,h}^{j,k} \varphi_{n,j}(t)\,, \quad \mbox{for}\; k = 1,2\,,
\end{equation*}
and $t \in I_n$ with coefficient functions $U_{n,h}^{j,k}\in W_h$ and $\vec Q_{n,h}^{j,k}
\in \vec V_h$. The continuity constraint that is imposed by the definition of 
$\mathcal{X}^{r}(W_h)$ and $\mathcal{X}^{r}(\vec V_h)$, respectively,
directly implies that $U_{n,h}^{0,1} = U_{n,h}^{0,2}$ and $\vec Q_{n,h}^{0,1} = \vec
Q_{n,h}^{0,2}$. Further, the pairs $\{u_{\tau,h}^{k} (t), \vec q_{\tau,h}^{k}
(t)\}$, for $k = 1,2$, both satisfy the discrete equations \eqref{Eq:FulDis1},
\eqref{Eq:FulDis2}. Therefore, it follows that
%in the case we have also the components Q^0, we have to add that they are the same
\begin{align}
\sum_{j=0}^r \hat \alpha_{ij} \langle U_{n,h}^{j,1}- U_{n,h}^{j,2} , w_h\rangle +
{\tau_n}
\, \hat  \beta_{ii} \langle \nabla \cdot (\vec Q_{n,h}^{i,1}- \vec Q_{n,h}^{i,2})
,w_h\rangle & = 0,  \label{Eq:proof_uniq_1}\\[1ex]
\langle \vec D^{-1} (\vec Q_{n,h}^{i,1}- \vec Q_{n,h}^{i,2}),\vec v_h \rangle - \langle
U_{n,h}^{i,1}-U_{n,h}^{i,2},\nabla \cdot  \vec v_h\rangle & = 0
\label{Eq:proof_uniq_2}
\end{align}
for $i=1,\ldots, r$ and all $\{w_h,\vec v_h\}\in W_h\times \vec V_h$. Now, by subtracting
the equations \eqref{Eq:proof_uniq_1} and \eqref{Eq:proof_uniq_2} from each other and
choosing the test functions $w_h= U_{n,h}^{i,1} - U_{n,h}^{i,2}$ and $\vec v_h =\tau_n
\hat \beta_{ii} (\vec Q_{n,h}^{i,1} - \vec Q_{n,h}^{i,2})$, for $i = 1,\ldots, r$ in
\eqref{Eq:proof_uniq_1} and \eqref{Eq:proof_uniq_2}, respectively,
 we get that
\begin{equation}
\label{Eq:proof_uniq_3}
\sum_{j=0}^r \hat \alpha_{ij} \langle U_{n,h}^{j,1} - U_{n,h}^{j,2}, U_{n,h}^{i,1} -
U_{n,h}^{i,2}\rangle + {\tau_n} \, \hat  \beta_{ii} \langle \vec D^{-1} (\vec
Q_{n,h}^{i,1} - \vec Q_{n,h}^{i,2}),\vec Q_{n,h}^{i,1} - \vec Q_{n,h}^{i,2} \rangle = 0\,,
\end{equation}
for $i=1,\ldots, r$. Summing up Eq.\ \eqref{Eq:proof_uniq_3} from $i = 1$ to $i=r$, using
Lemma \ref{lemma_timeintegration_basicfunc} and recalling that $U_{n,h}^{0,1} =
U_{n,h}^{0,2}$ then implies that
\begin{equation}
%\begin{split}
\label{Eq:proof_uniq_4}
 \dfrac{1}{2} \|u_{\tau,h}^{1} (t_n)  - u_{\tau,h}^{2} (t_n) \|^2 + \sum_{i = 1}^r
{\tau_n} \, \hat  \beta_{ii}
\langle \vec D^{-1} (\vec Q_{n,h}^{i,1} - \vec Q_{n,h}^{i,2}),\vec Q_{n,h}^{i,1} - \vec
Q_{n,h}^{i,2} \rangle = 0\,.
%\end{split}
\end{equation}
%The matrix $(\hat \alpha_{ij})_{i,j = 1, \dots, r}$ in Eq.\ \eqref{Eq:proof_uniq_4} is  positive definite. This is shown in \cite[p.\ 1784]{Karakashian1999}.
The symmetric matrix $\vec D^{-1}$ is positive definite by assumption \eqref{Eq:PosDefD}
and $\hat  \beta_{ii} > 0$ under the coefficient property (C); cf.\ Lemma 
\ref{Lem:CoePropC}. Therefore, Eq.\
\eqref{Eq:proof_uniq_4} immediately implies that  $\vec Q_{n,h}^{i,1} = \vec
Q_{n,h}^{i,2}$ for $i = 1, \dots, r$. By Lemma \ref{lemma_thomas} there exists some $\vec
v_h \in \vec V_h$ such that $\nabla \cdot \vec v_h = U_{n,h}^{i,1}-U_{n,h}^{i,2}$. Using
this $\vec v_h$ as test function in \eqref{Eq:proof_uniq_2} and noting that the first term
in \eqref{Eq:proof_uniq_2} now vanishes, we obtain that $U_{n,h}^{i,1} = U_{n,h}^{i,2}$,
for $i = 1, \dots, r$.  This implies the uniqueness of the solution to the fully discrete
problem \eqref{Eq:FulDis1}, \eqref{Eq:FulDis2} and proves the assertion of the theorem.
\end{mproof}

\subsection{Estimates for the error between the semidiscrete and the fully
discrete solution}
\label{Sec:ErrMFEM}

In this subsection we derive estimates for the error between the
semidiscrete approximation defined by Eqs.\ \eqref{Eq:GtdP_1}, \eqref{Eq:GtdP_2} and the
fully  discrete solution given by Eqs.\ \eqref{Eq:FulDis1}, \eqref{Eq:FulDis2}. For this
we use  the following projection operators (cf.\ \cite{Brezzi1991}, \cite{Arbogast96} and \cite[p.\
237]{Quarteroni1994}) defined in $ W$  and $\vec V$,
respectively, by
\begin{equation}
P_h: W \rightarrow W_h,
\qquad \la P_h w - w, w_h \ra = 0
\label{Eq:l2_projector}
\end{equation}
for all $w_h \in W_h$ and
\begin{align}
\label{Eq:hdiv_projector}
\vec \Pi_h: & \; \vec V \rightarrow \vec V_h, \qquad &&
\la \nabla \cdot (\vec \Pi_h \vec v - \vec v) , w_h \ra = 0, \\[1ex]
\label{Eq:hdiv_projector2}
\vec P_h: & \; \vec V \rightarrow \vec V_h \qquad && \la \vec P_h \vec v - \vec v, \vec
v_h \ra = 0,
\end{align}
for all $w_h\in W_h$ and $\vec v_h \in \vec V_h$, respectively. We point out that $\vec \Pi_h$ 
is first defined on $\vec H^1(\Omega)$ and then extended to $\vec V$ by following  \cite[p.\ 237]{Quarteroni1994}. For these operators and the family of Raviart--Thomas 
elements on quadrilateral elements for the two-dimensional case and the class of 
Raviart--Thomas--N\'{e}d\'{e}lec elements in three space dimensions there holds that
\begin{eqnarray}
\label{lemmaproiectors}
\| w - P_h w \| \leq c h^{p+1} \| w \|_{p+1},  &&\label{lemmaproiectors_1} \\[1ex]
\| \vec v - \vec \Pi_h \vec v \| \leq c h^{p+1} \| \vec v \|_{p+1}, && \| \nabla
\cdot (\vec v - \vec \Pi_h \vec v) \| \leq c h^{p+1} \| \nabla \cdot \vec v
\|_{p+1},\label{lemmaproiectors_2} \\[1ex]
\| \vec v - \vec P_h \vec v \| \leq c h^{p+1} \| \vec v \|_{p+1}, && \| \nabla \cdot
(\vec v - \vec P_h \vec v) \| \leq c h^{p+1} \| \nabla \cdot \vec v
\|_{p+1},\label{lemmaproiectors_3}
\end{eqnarray}
for any  $w \in H^{p + 1}(\Omega)$ and $\vec v \in \vec H^{p+1}(\Omega)$,  {\color{black} $\nabla
\cdot \vec v \in H^{p+1}(\Omega)$}, respectively.

For the error between the semidiscrete solution and fully discrete we use the notation
\begin{equation*}
\label{Eq:def_error_semi2fully}
\begin{array}{r@{\;}c@{\;}l@{\qquad}r@{\;}c@{\;}l}
E_{u}(t) & = & u_\tau (t) -  u_{\tau,h} (t), & \vec E_{\vec q}(t) & = & \vec q_\tau (t) -
\vec q_{\tau,h}(t)\,,\\[1ex]
E_{u,n}^i & = & E_u(t_{n,i}), & \vec E_{\vec q,n}^i & = & \vec E_{\vec q}(t_{n,i})
\end{array}
\end{equation*}
for $t\in I$ and $n \in \{1,\ldots,N\}$, $ i \in \{0,\ldots,r\}$. Representing the 
semidiscrete and fully discrete solution in terms of basis functions (cf.\ 
\eqref{Eq:RepBasis}) there holds that
\begin{equation*}
\label{Eq:FormulaEonIn}
E_{u}(t)=\sum_{i=0}^r E_{u,n}^i \varphi_{n,i}(t) \quad \mbox{ and } \quad \vec
E_{\vec q}(t)=\sum_{i=0}^r \vec E_{\vec q,n}^i \varphi_{n,i}(t),  \quad \mbox{for} \;
t\in I_n\,.
\end{equation*}

{\color{black} Next, we prove two preliminary lemmas.}
\begin{lemm}
\label{lemma_semi2fully_1}
Let the assumptions of Subsec.\ \ref{Sec:NotPrem_2} about $\Omega, u_0, \vec D$ and $f$ be
satisfied. Let the semidiscrete approximation $\{u_{\tau},\vec q_{\tau}\}\in 
\mathcal{X}^{r}(W)\times \mathcal{X}^{r}(\vec V)$ be defined by 
\eqref{Eq:RepBasis}--\eqref{Eq:SemiDis2}. Further, let  $\{u_{\tau,h},\vec q_{\tau,h}\}\in
\mathcal{X}^{r}(W_h)\times \mathcal{X}^{r}(\vec V_h)$ be the unique solution of
the fully discrete problem \eqref{Eq:FulDis1}, \eqref{Eq:FulDis2}. Then, for any $K 
=1,\ldots , N$ it holds that
\begin{equation}
\label{Eq:error_est_semi2fully_1}
\begin{array}{l}
\ds \| E_{u} (t_K) \|^2 + \sum_{n=1}^K \sum_{i=1}^r {\tau_n}\|  E_{u, n}^i  \|^2 +
\sum_{n=1}^K \sum_{i=1}^r {\tau_n}  \| \vec E_{\vec q, n}^i \|^2 \\[3ex]
\ds \hspace*{0.7cm} \le  \| u_\tau (t_K)  - P_h u_\tau (t_K) \|^2  + c \sum_{n=1}^K
\sum_{i=1}^r {\tau_n}   (\| U_{n}^i - P_h U_{n}^i \|^2 +  \| \vec Q_{n}^i -
\vec \Pi_h \vec Q_{n}^i \|^2)
\end{array}\end{equation}
with some constant $c > 0$ not depending on the discretization parameters $h$ and
$\tau$.
\end{lemm}

\begin{mproof}
By subtracting \eqref{Eq:FulDis1}, \eqref{Eq:FulDis2} from
\eqref{Eq:SemiDis1}, \eqref{Eq:SemiDis2}, respectively, it follows that
\begin{align}
\sum_{j=0}^r \hat \alpha_{ij} \langle U_{n}^j - U_{n,h}^j, w_h\rangle + {\tau_n} \, \hat
\beta_{ii} \langle \nabla \cdot (\vec Q_{n}^i -  \vec Q_{n,h}^i),w_h\rangle & = 0
\label{eq_proof_semi2fully_1}\,,\\[1ex]
\langle \vec D^{-1} (\vec Q_{n}^i -  \vec Q_{n,h}^i),\vec v_h \rangle - \langle U_{n}^i -
U_{n,h}^i,\nabla \cdot \vec v_h\rangle & = 0 \label{eq_proof_semi2fully_2}
\end{align}
for $i=1,\ldots, r$ and all $\{w_h,\vec v_h\}\in W_h\times \vec V_h$. For any
$i=1,\ldots, r$ we choose the test functions $w_h =  P_h U_{n}^i - U_{n,h}^i \in W_h$
and $\vec v_h = {\tau_n} \, \hat  \beta_{ii} (\vec \Pi_h \vec Q_{n}^i -  \vec Q_{n,h}^i)
\in \vec V_h$ in \eqref{eq_proof_semi2fully_1} and \eqref{eq_proof_semi2fully_2}, 
respectively. By adding the thus obtained equations, using the properties of the 
projection projectors $P_h$
and $\vec \Pi_h$ defined in \eqref{Eq:l2_projector} and \eqref{Eq:hdiv_projector},
respectively, and summing up from  $i = 1 $ to  $r $ we get that
\begin{equation}
\label{eq_proof_semi2fully_3}
\begin{split}
\sum_{i=1}^r \sum_{j=0}^r \hat \alpha_{ij} \langle & P_h U_{n}^j - U_{n,h}^j, P_h U_{n}^i
- U_{n,h}^i \rangle \\
& + \sum_{i=1}^r {\tau_n} \, \hat  \beta_{ii} \langle \vec
D^{-1} (\vec Q_{n}^i -  \vec Q_{n,h}^i), \vec \Pi_h \vec Q_{n}^i -  \vec Q_{n,h}^i
\rangle =
0\,.
\end{split}
\end{equation}

We note that due to Lemma \ref{lemma_timeintegration_basicfunc}, the first term
in \eqref{eq_proof_semi2fully_3} can be rewritten as
\begin{displaymath}
\begin{split}
\sum_{i=1}^r \sum_{j=0}^r \hat \alpha_{ij} \langle & P_h U_{n}^j - U_{n,h}^j, P_h U_{n}^i
- U_{n,h}^i \rangle\\
& =  \dfrac{1}{2}\| P_h E_{u,n} (t_n) \|^2 - \dfrac{1}{2}\| P_h E_{u,n-1}
(t_{n-1}) \|^2\,.
\end{split}
\end{displaymath}
Along with some further algebraic manipulations we then conclude
from \eqref{eq_proof_semi2fully_3} that
\begin{equation}\label{eq_proof_semi2fully_4}
\begin{split}
\dfrac{1}{2}\| P_h E_{u,n} (t_n) \|^2 & - \dfrac{1}{2}\| P_h E_{u,n-1} (t_{n-1}) \|^2\\
& \qquad + \sum_{i=1}^r {\tau_n} \, \hat  \beta_{ii} \langle \vec D^{-1} (\vec Q_{n}^i -
 \vec
Q_{n,h}^i), \vec Q_{n}^i -  \vec Q_{n,h}^i \rangle \\
& =  \sum_{i=1}^r {\tau_n} \, \hat  \beta_{ii} \langle \vec D^{-1} (\vec Q_{n}^i -
 \vec Q_{n,h}^i), \vec Q_{n}^i - \vec \Pi_h \vec Q_{n}^i \rangle\,.
\end{split}
\end{equation}
Recalling assumption \eqref{Eq:PosDefD} about $\vec D$ and property 
(C) \eqref{Eq:cond_on_beta} of the coefficients $\hat  \beta_{ii}$ and we obtain from 
Eq.\ \eqref{eq_proof_semi2fully_4} by applying Cauchy--Young's inequality that
\begin{equation}
\label{eq_proof_semi2fully_5}
\begin{split}
\ds \| P_h E_{u} (t_n) \|^2 - & \| P_h E_{u} (t_{n-1}) \|^2 + \sum_{i=1}^r {\tau_n}
\beta_m \theta_m \| \vec Q_{n}^i -  \vec Q_{n,h}^i\|^2 \\
& \le \frac{\beta_M^2}{\beta_m \, \theta_m} \sum_{i=1}^r {\tau_n}  \| \vec Q_{n}^i -
\vec \Pi_h \vec Q_{n}^i \|^2\,.
 \end{split}
\end{equation}

Summing up inequality \eqref{eq_proof_semi2fully_5} from $n = 1$ to $K$ and noting that $
P_h E_{u}(t_{0}) = 0 $ then shows that
\begin{equation}
\label{eq_proof_semi2fully_6}
\begin{split}
\| P_h E_{u} (t_K) \|^2 + \sum_{n=1}^K \sum_{i=1}^r & {\tau_n} \beta_m \| \vec Q_{n}^i
-  \vec Q_{n,h}^i\|^2 \\
& \le \frac{\beta_M^2}{\beta_m \, \theta_m}   \sum_{n=1}^K \sum_{i=1}^r {\tau_n}  \|
\vec Q_{n}^i - \vec \Pi_h \vec Q_{n}^i \|^2
\end{split}
\end{equation}
for any $K \in \N$ with $K \le N$. By using now Lemma \ref{lemma_thomas}, there exists
for any $i \in \{1, \ldots, r\}$ a $ \vec v_h \in \vec V_h $ such that  $\nabla \cdot \vec
v_h = P_h E_{u,n}^i $ and  $ \| \vec v_h \| \le c \| P_h E_{u,n}^i \| $. By
testing \eqref{eq_proof_semi2fully_2} with this $ \vec v_h$,  we get by using the
Cauchy--Schwarz inequality along with assumption \eqref{Eq:PosDefD} about $\vec D$ that
\begin{equation}
\label{eq_proof_semi2fully_7}
 \| P_h E_{u,n}^i \| \le c \, \theta_M \| \vec Q_{n}^i -  \vec Q_{n,h}^i \|,
\end{equation}
for  $n=1,\ldots, N$, $i=1,\ldots, r$. 

{\color{black} Combining \eqref{eq_proof_semi2fully_6} with \eqref{eq_proof_semi2fully_7} is follows that 
\begin{equation}
\label{Eq:Ad_01}
\begin{split}
\| P_h E_{u} (t_K) \|^2 + \sum_{n=1}^K \sum_{i=1}^r & {\tau_n} \beta_m \| \vec Q_{n}^i
-  \vec Q_{n,h}^i\|^2 + \sum_{n=1}^K \sum_{i=1}^r \tau_n \| P_h E_{u,n}^i \|^2 \\
& \le c \sum_{n=1}^K \sum_{i=1}^r \tau_n \| \vec Q_{n}^i -  \vec \Pi_h \vec Q_{n}^i \|^2\,.
\end{split}
\end{equation}
By $\vec E_{\vec q,n}^{i} = \vec Q_{n}^i -  \vec Q_{n,h}^i $ and the triangle inequality relation \eqref{Eq:Ad_01} implies that 
\begin{equation}
\label{Eq:Ad_02}
\begin{split}
\|& E_{u} (t_K) \|^2  + \sum_{n=1}^K \sum_{i=1}^r  {\tau_n} \beta_m \| \vec E_{\vec q,n}^{i} \|^2 + \sum_{n=1}^K \sum_{i=1}^r \tau_n \| E_{u,n}^i \|^2 \\
& \le c \sum_{n=1}^K \sum_{i=1}^r \tau_n \| \vec Q_{n}^i -  \vec \Pi_h \vec Q_{n}^i \|^2 + c \| u_{\tau}(t_k) - P_h u_{\tau}(t_k)\|^2\\
& \quad + c \sum_{n=1}^K \sum_{i=1}^r \tau_n  \| E_{u,n}^i - P_h E_{u,n}^i\|^2\,.
\end{split}
\end{equation}
Observing that $E_{u,n}^i - P_h E_{u,n}^i = U_n^i - P_h U_n^i$, inequality \eqref{Eq:Ad_02} proves \eqref{Eq:error_est_semi2fully_1}.
}
\end{mproof}

In the second lemma we restrict ourselves to the case that $\vec D = d\vec I$ with some $d>0$ is satisfied. An extension of the provided  estimates to more general matrices $\vec D(\vec x)$ still remains
an open problem. %and has to be left for our future work. We note that the error estimate of Lemma \ref{lemma_semi2fully_2} still holds under the assumption of a non-realistic permutation property of $\vec D(\vec x)$ with $\vec P_h$.

\begin{lemm}
\label{lemma_semi2fully_2}
Let the assumptions of Subsec.\ \ref{Sec:NotPrem_2} about $\Omega, u_0$ and $f$ be satisfied and $\vec D = d\vec I$ with some $d>0$. Let the semidiscrete approximation $\{u_{\tau},\vec  q_{\tau}\}\in \mathcal{X}^{r}(W)\times \mathcal{X}^{r}(\vec V)$ be defined by 
\eqref{Eq:RepBasis}--\eqref{Eq:SemiDis2}. Further, let  $\{u_{\tau,h},\vec q_{\tau,h}\}\in \mathcal{X}^{r}(W_h)\times \mathcal{X}^{r}(\vec V_h)$ be the unique solution of
the fully discrete problem \eqref{Eq:FulDis1}, \eqref{Eq:FulDis2}. Then, for any $K 
=1,\ldots , N$
it holds that
\begin{equation}\label{eq_proof_div_6}
\begin{split}
\sum_{n=1}^K \tau_n \sum_{i=1}^r \hat \beta_{ii}\| \nabla \cdot \vec \Pi_h
\vec E_{\vec q,n}^i \|^2& + \| \vec P_h \vec E_{\vec q} (t_K) \|^2 \\
& \le
\sum_{n=1}^K \sum_{i=1}^r \tau_n \, \hat  \beta_{ii}  \| \nabla \cdot \vec (\vec P_h - \vec \Pi_h)
\vec  Q_{n}^i\|^2\,.
\end{split}
\end{equation}
\end{lemm}

\begin{mproof} Introducing the projectors into the error equations
\eqref{eq_proof_semi2fully_1}--\eqref{eq_proof_semi2fully_2} yields that
\begin{align}
\sum_{j=0}^r \hat \alpha_{ij} \langle P_h E_{u,n}^j, w_h\rangle + {\tau_n} \, \hat
\beta_{ii} \langle \nabla \cdot \vec \Pi_h \vec E_{\vec q,n}^i,w_h\rangle & = 0\,,
\label{eq_proof_div_1}\\[1ex]
\langle \vec P_h  \vec E_{\vec q,n}^i, \vec v_h \rangle - \langle P_h E_{u,n}^i,\nabla
\cdot \vec v_h\rangle & = 0 \label{eq_proof_div_2}
\end{align}
for $n=1,\ldots, N$, $i=1,\ldots, r$ and all $\{w_h,\vec v_h\}\in W_h\times \vec V_h$.
Observing that for any $ n \ge 2 $ the quantities $\vec E_{\vec q,n}^0$ and $E_{u,n}^0$
are  linear combinations of $\vec E_{\vec q,n-1}^i$ and $\vec E_{u,n-1}^i$, for $i =
0, \dots, r$, respectively, and that $P_h E_{u,1}^0 = 0$ and $\vec
P_h \vec E_{\vec q,1}^0 = \vec 0$ by definition of $\{U_1^0, \vec Q_1^0\}$ and $\{U_{n,h}^1,\vec Q_{n,h}^1\}$, it 
follows that Eq.\ \eqref{eq_proof_div_2}
is also satisfied for $i = 0$ and any $n \ge 1$. Using this, we obtain by multiplying
\eqref{eq_proof_div_2} with $\hat \alpha_{ji}$ and summing up the resulting identity from
$i = 0$ to $r$ that
\begin{equation}\label{eq_proof_div_3}
\left\langle \sum_{j = 0}^r \hat \alpha_{ij} \vec P_h  \vec E_{\vec q,n}^j, \vec v_h
\right\rangle - \left\langle \sum_{j = 0}^r \hat \alpha_{ij} P_h E_{u,n}^j,\nabla \cdot
\vec v_h\right\rangle  = 0
\end{equation}
for any  $\vec v_h \in \vec V_h$. We note that we changed the notation for the indices.
By testing now \eqref{eq_proof_div_1} with $w_h = \sum_{j = 0}^r \hat \alpha_{ij} P_h
E_{u,n}^j \in W_h$  and \eqref{eq_proof_div_3} with $\vec v_h = \tau_n \, \hat
\beta_{ii} \vec P_h  \vec E_{\vec q,n}^i \in \vec V_h$, we get by summing the resulting
equations and using the inequalities of Cauchy--Schwarz and Cauchy--Young that
\begin{align*}
\Big\|  \sum_{j = 0}^r \hat \alpha_{ij} P_h E_{u,n}^j \Big\|^2 + & \; \tau_n \,
\hat  \beta_{ii} \Big\langle \sum_{j = 0}^r \hat \alpha_{ij} \vec P_h  \vec E_{\vec
q,n}^j,
\vec
P_h  \vec E_{\vec q,n}^i \Big\rangle \\[1ex]
= & \; \tau_n \, \hat  \beta_{ii} \Big\langle  \sum_{j = 0}^r \hat \alpha_{ij} P_h
E_{u,n}^j, \nabla \cdot \vec (\vec P_h - \vec \Pi_h) \vec  E_{\vec q,n}^i
\Big\rangle\\[1ex]
\le & \; \dfrac{1}{2} \Big\|  \sum_{j = 0}^r \hat \alpha_{ij} P_h E_{u,n}^j \Big\|^2 +
\dfrac{1}{2} \tau_n^2 \, \hat  \beta_{ii}^2  \| \nabla \cdot (\vec P_h - \vec \Pi_h)
\vec  E_{\vec q,n}^i\|^2
\end{align*}
for  $n=1,\ldots, N$ and $i=1,\ldots, r$. The inequality above further simplifies to
\begin{equation}\label{eq_proof_div_4}
\begin{split}
 \Big\|  \sum_{j = 0}^r \hat \alpha_{ij} P_h E_{u,n}^j \Big\|^2 + 2 \tau_n \, \hat
\beta_{ii} \Big\langle \sum_{j = 0}^r & \hat \alpha_{ij} \vec P_h  \vec E_{\vec q,n}^j,
\vec
P_h  \vec E_{\vec q,n}^i \Big\rangle\\[1ex]
& \le \ds \tau_n^2 \, \hat  \beta_{ii}^2  \| \nabla \cdot (\vec P_h - \vec \Pi_h) \vec
E_{\vec q,n}^i\|^2,
\end{split}
\end{equation}
for  $n=1,\ldots, N$, $i=1,\ldots, r$. Dividing \eqref{eq_proof_div_4} by $\tau_n \,\hat
\beta_{ii}$ (note that $\hat \beta_{ii} > 0$ for all $i =1, \ldots, r$),
summing up the resulting inequality from  $i=1,\ldots, r$ and using Lemma
\ref{lemma_timeintegration_basicfunc} gives that
\begin{equation}\label{eq_proof_div_5}
\begin{split}
\sum_{i=1}^r \dfrac{1}{ \tau_n \, \hat  \beta_{ii} }   \Big\|  \sum_{j = 0}^r \hat
\alpha_{ij} P_h E_{u,n}^j \Big\|^2  & +  \| \vec P_h \vec E_{\vec q} (t_n) \|^2 \le \ds
\|
\vec P_h \vec E_{\vec q} (t_{n-1}) \|^2 \\[0ex]
& + \sum_{i=1}^r \tau_n \, \hat  \beta_{ii}  \| \nabla \cdot (\vec P_h - \vec \Pi_h)
\vec  E_{\vec q,n}^i\|^2
\end{split}
\end{equation}
for  $n=1,\ldots, N$. By summing up \eqref{eq_proof_div_5} from $n=1,\ldots, K$ and noting that $\vec P_h \vec E_{\vec q} (t_{0}) = \vec P_h \vec E_{\vec q,1}^0=\vec 0$ for the choices of the semidiscrete and fully discrete coefficient functions  $\vec Q_{1}^0$ and $\vec Q_{1,h}^0$ (cf.\ their definition below \eqref{Eq:SemiDis1}, \eqref{Eq:SemiDis2} and Eq.\ \eqref{Eq:ContCond3_2}) we get that
\begin{equation}
\label{eq_proof_div_6_1}
\begin{split}
\sum_{n=1}^K \sum_{i=1}^r \dfrac{1}{ \tau_n \, \hat  \beta_{ii} }   \Big\|  \sum_{j =
0}^r \hat \alpha_{ij} & P_h E_{u,n}^j \Big\|^2 +   \| \vec P_h \vec E_{\vec q} (t_K) \|^2
\\
& \le
\sum_{n=1}^K \sum_{i=1}^r \tau_n \, \hat  \beta_{ii}  \| \nabla \cdot \vec (\vec P_h - \vec \Pi_h)
\vec  E_{\vec q,n}^i\|^2.
\end{split}
\end{equation}

We now estimate the divergence of the flux. By testing \eqref{eq_proof_div_1} with
$w_h = \nabla \cdot \vec \Pi_h  \vec E_{\vec q,n}^i \in W_h$, and using the inequalities of Cauchy--Schwarz and Cauchy--Young ($\hat \beta_{ii} > 0$ for
all $i=1, \ldots, r$) we get that
\begin{align*}\label{eq_proof_div_6}
 {\tau_n} \, \hat  \beta_{ii} \| \nabla \cdot \vec \Pi_h \vec E_{\vec q,n}^i \|^2 & =  -
\Big\langle \sum_{j=0}^r \hat \alpha_{ij} P_h E_{u,n}^j,\nabla \cdot \vec \Pi_h \vec
E_{\vec q,n}^i \Big\rangle \nonumber \\
  & \le  \dfrac{1}{2 {\tau_n} \, \hat  \beta_{ii}} \Big\| \sum_{j=0}^r \hat \alpha_{ij}
P_h E_{u,n}^j\Big\|^2 +  \dfrac{{\tau_n} \, \hat  \beta_{ii}}{2} \Big\| \nabla \cdot
\vec \Pi_h \vec E_{\vec q,n}^i \Big\|^2
\end{align*}
for  $n=1,\ldots, N$, $i=1,\ldots, r$. Summing up the previous inequality from $n=1,\ldots, K$ as well as from $i=1,\ldots, r$, using {\color{black} \eqref{eq_proof_div_6_1}} along with $(\vec P_h - \vec \Pi_h) \vec E_{\vec q,n}^i = (\vec P_h - \vec \Pi_h) \vec Q_{n}^i$ by definition of the projectors $\vec P_h$ and $\vec \Pi_h$ we
obtain that
\begin{displaymath}
\sum_{n=1}^K \sum_{i=1}^r  {\tau_n} \, \hat  \beta_{ii} \| \nabla \cdot \vec \Pi_h \vec
E_{\vec q,n}^i \|^2 \le  \sum_{n=1}^K \sum_{i=1}^r \tau_n \, \hat  \beta_{ii}  \| \nabla
\cdot (\vec P_h - \vec \Pi_h) \vec  Q_n^i\|^2\,,
\end{displaymath}
which proves the assertions of the lemma.
\end{mproof}

Now we combine the inequalities of the previous lemmas to estimate the error between the semidiscrete and the
fully discrete solutions in the norms of $L^2(I; W)$ and $L^2(I; \vec V)$.

\begin{theorem} \label{Thm:L2norm_error_semi2fully}
Let the assumptions of Subsec.\ \ref{Sec:NotPrem_2} about $\Omega, u_0, \vec D$ and $f$ be satisfied. Let the sufficienctly regular semidiscrete approximation $\{u_{\tau},\vec q_{\tau}\}\in \mathcal{X}^{r}(W)\times \mathcal{X}^{r}(\vec V)$ be defined by 
\eqref{Eq:RepBasis}--\eqref{Eq:SemiDis2}. Further, let  $\{u_{\tau,h},\vec q_{\tau,h}\}\in
\mathcal{X}^{r}(W_h)\times \mathcal{X}^{r}(\vec V_h)$ be the solution of the fully discrete problem \eqref{Eq:FulDis1}, \eqref{Eq:FulDis2}. 
For the scalar variable $u_\tau$ it holds that
\begin{align}
\| u_\tau - u_{\tau, h} \|_{L^2(I; W)} &\le c h^{p+1} \,.
\label{L2norm_error_semi2fully_1}
\end{align}
For the vectorial variable $\vec q_\tau$ it holds that 
\begin{align}
{\color{black} \left(\sum_{n=1}^N \tau_n \sum_{i=1}^r \| \vec q_\tau(t_{n,i}) - \vec q_{\tau,h}(t_{n,i}) \|^2\right)^{1/2} \leq c h^{p+1}\,.}
\label{L2norm_error_semi2fully_3}
\end{align}
Further, for $\vec D(\vec x) = d \vec I$, for some $d>0$, it holds that
\begin{align}
\| \vec q_\tau - \vec q_{\tau, h} \|_{L^2(I; \vec L^2(\Omega))} &\le c h^{p+1}
\label{L2norm_error_semi2fully_2}
\end{align}
and
\begin{align}
{\color{black} \left(\sum_{n=1}^N \tau_n \sum_{i=1}^r \| \vec q_\tau(t_{n,i}) - \vec q_{\tau,h}(t_{n,i}) \|_{\vec V}^2\right)^{1/2} \leq c h^{p+1}\,.}
\label{L2norm_error_semi2fully_4}
\end{align}
The constant $c$ does not depend on the discretization parameters 
$h$ and $\tau$.
\end{theorem}

\begin{mproof}
By using \eqref{result_lemma_L2_I_X_error} and recalling that $E_{u,n}^0 = E_u(t_{n-1})$ we find that 
\begin{align}
\| u_\tau - u_{\tau, h}\|_{L^2(I; W)}^2 \leq c \bigg(\sum_{n=1}^N \sum_{i=1}^r \tau_n  \|E_{u,n}^i \|^2 + \sum_{n=1}^N \tau_n \|E_{u}(t_{n-1})\|^2 \bigg)\,.
\label{Eq:SemiSpace_00}
\end{align}
By inequality \eqref{Eq:error_est_semi2fully_1} we get for the first term on the right-hand side of \eqref{Eq:SemiSpace_00} that 
\begin{equation}
\label{Eq:SemiSpace_01}
\begin{split}
\sum_{n=1}^N   \sum_{i=1}^r \tau_n   \|  
E_{u,n}^i \|^2 & {\color{black} + \sum_{n=1}^N \sum_{i=1}^r \tau_n \|
\vec E_{\vec q,n}^i \|^2 } \le  c \bigg(\sum_{n=1}^N \sum_{i=1}^r \tau_n \| U_n^i - P_h U_n^i \|^2\\
&  +  \tau_N \| u_{\tau}(t_N) - P_h u_{\tau}(t_N) \|^2 + \sum_{n=1}^N \sum_{i=1}^r \tau_n \| \vec Q_n^i - \vec \Pi_h \vec Q_n^i \|^2\bigg)\,.
\end{split}
\end{equation}
Using \eqref{Eq:error_est_semi2fully_1} again, we find for the second term on the right-hand side of \eqref{Eq:SemiSpace_00} that 
\begin{equation}
\label{Eq:SemiSpace_02}
\begin{aligned}
\| E_{u}(t_K)\|^2 & \leq \| u_\tau (t_K)-P_h u_\tau(t_k)\|^2 \\
& \quad + c \sum_{n=1}^K \sum_{i=1}^r (\|U_n^i - P_h U_n^i \|^2 + \| \vec Q_n^i - \vec \Pi_h \vec Q_n^i\|^2) 
\end{aligned}
\end{equation}
for $K=1,\ldots, N$. 
Combining now \eqref{Eq:SemiSpace_00} with \eqref{Eq:SemiSpace_01} and \eqref{Eq:SemiSpace_02} and using the approximation properties \eqref{lemmaproiectors_1}--\eqref{lemmaproiectors_3} of the projection operators we then get that
\begin{equation}
\begin{split}
\| u_\tau  - u_{\tau, h} & \|_{L^2(I; W)}^2\le  c \bigg(\sum_{n=1}^N \tau_n  h^{2(p+1)} \sum_{i=0}^r  \| U_n^i \|_{p+1}^2 \\ 
&  +  h^{2(p+1)} \max_{K=0,\ldots, N} \| u_{\tau}(t_K) \|_{p+1}^2  + \sum_{n=1}^N \tau_n h^{2(p+1)} \sum_{i=0}^r \| \vec Q_n^i 
\|_{p+1}^2 \bigg) \,,
\end{split}
\label{Eq:EstSpaceScal00}
\end{equation}
where the arising constant does not depend on the discretization parameters $h$ and $\tau$. The result \eqref{L2norm_error_semi2fully_1} directly follows from 
\eqref{Eq:EstSpaceScal00} under the assumption
of the theorem of sufficiently regular coefficient functions $\{U_n^j,\vec Q_n^j\}\in W\times \vec V$. {\color{black}
From \eqref{Eq:SemiSpace_01} along with \eqref{Eq:SemiSpace_02} and \eqref{lemmaproiectors_1}--\eqref{lemmaproiectors_3} we further conclude that 
\begin{equation}
\label{Eq:l2V_01}
\sum_{n=1}^N \tau_n \sum_{i=1}^r \| \vec E_{\vec q}(t_{n,i}) \|^2 \leq c h^{2(p+1)}\,.
\end{equation}
This proves \eqref{L2norm_error_semi2fully_3}.
}

{\color{black}
By using \eqref{result_lemma_L2_I_X_error}, recalling that $\vec E_{\vec q,n}^0 = \vec E_{\vec q}(t_{n-1})$ and applying the boundedness of the $\vec L^2$ projection operator $\vec P_h$  we find that 
\begin{align}
\nonumber
& \| \vec q_\tau - \vec q_{\tau, h}  \|_{L^2(I; \vec L^2(\Omega))}^2 
  \\ \nonumber
&\le c \Big(\| \vec q_\tau - \vec P_h \vec q_{\tau} \|_{L^2(I; \vec L^2(\Omega))}^2 +  \| \vec P_h \vec q_\tau - \vec q_{\tau, h}  \|_{L^2(I; \vec L^2(\Omega))}^2\Big) \\ \nonumber
& \le  c \bigg(\sum_{n=1}^N \sum_{i=1}^r \tau_n \|\vec E_{\vec q,n}^i \|^2 + \sum_{n=1}^N \tau_n \|\vec P_h \vec E_{\vec q}(t_{n-1})\|^2\\
\label{Eq:SemiSpace_10}
& \qquad  + \| \vec q_\tau - \vec P_h\vec q_{\tau}  \|_{L^2(I; \vec L^2(\Omega))}^2 \bigg)\,.
\end{align}
For $\vec D(\vec x) = d \vec I$ the second term on the right-hand side of \eqref{Eq:SemiSpace_10} can be bounded from above by means of the inequality \eqref{eq_proof_div_6} along with the observation that $(\vec P_h - \vec \Pi_h) \vec E_{\vec q,n}^i = (\vec P_h - \vec \Pi_h) \vec Q_{n}^i$ by definition of the projectors $\vec P_h$ and $\vec \Pi_h$. Recalling further the boundedness of $\hat \beta_{ii}$ (cf.\ Lem.\ \ref{Lem:CoePropC}) we conclude that 
\begin{equation}
\label{Eq:SemiSpace_03}
\| \vec P_h \vec E_{\vec q}(t_K)\|^2 \leq c \sum_{n=1}^K \sum_{i=1}^r \tau_n \|
(\vec P_h - \vec \Pi_h) \vec Q_{n}^i \|^2
\end{equation}
for $K=1,\ldots, N$. Finally, combining \eqref{Eq:SemiSpace_10} with \eqref{Eq:SemiSpace_01} and \eqref{Eq:SemiSpace_03} and using the approximation properties \eqref{lemmaproiectors_1}--\eqref{lemmaproiectors_3} of the projection operators we then get that
\begin{align}
\nonumber
 \| \vec q_\tau & - \vec q_{\tau, h}  \|_{L^2(I; \vec L^2(\Omega))}^2\\
\label{Eq:EstSpaceScal2}
& \le\; c \bigg(\sum_{n=1}^N \tau_n 
h^{2(p+1)} \sum_{i=0}^r  \| U_n^i \|_{p+1}^2 +  h^{2(p+1)} \max_{K=0,\ldots, N} \| 
u_{\tau}(t_K) \|_{p+1}^2  \\ \nonumber 
& \qquad + \sum_{n=1}^N \tau_n h^{2(p+1)} \sum_{i=0}^r \| \vec Q_n^i 
\|_{p+1}^2 {\color{black} + h^{2(p+1)} \|\vec q_\tau \|^2_{L^2(I;\vec H^{p+1}(\Omega))}}\bigg) \,,
\end{align}
where the arising constant does not depend on the discretization parameters $h$ and $\tau$. The result \eqref{L2norm_error_semi2fully_2} directly follows from 
\eqref{Eq:EstSpaceScal2} under the assumption of sufficiently regular coefficient functions $\{U_n^j,\vec Q_n^j\}\in W\times \vec V$. 

To estimate the divergence part of the error in  \eqref{L2norm_error_semi2fully_4}, we use that by definition of the projection operators it holds that 
\begin{equation}
\label{Eq:l2V_02}
\begin{split}
& \sum_{n=1}^N \tau_n \sum_{i=1}^r \| \nabla \cdot \vec E_{\vec q}(t_{n,i}) \|^2 \\
& \leq  \sum_{n=1}^N \tau_n \sum_{i=1}^r \| \nabla \cdot (\vec Q_{n}^i- \vec \Pi_h \vec Q_{n}^i) \|^2
+ \sum_{n=1}^N \tau_n \sum_{i=1}^r \| \nabla \cdot \vec \Pi_h \vec E_{\vec q,n}^i \|^2 \,.
\end{split}
\end{equation}
The assertion \eqref{L2norm_error_semi2fully_4} then follows from \eqref{Eq:l2V_02} combined with \eqref{eq_proof_div_6} and the approximation properties  \eqref{lemmaproiectors_1}--\eqref{lemmaproiectors_3}.
}
\end{mproof}

{\color{black}
We remark that the inequalities \eqref{L2norm_error_semi2fully_2} and \eqref{L2norm_error_semi2fully_4} provide an error control for the spatial discretization in the Gaussian quadrature points or temporal degrees of freedom of the subintervals $I_n$ with respect to the norm of $\vec L^2(\Omega)$ and $\vec V$, respectively. For an error control with respect to the norm of $L^2(I;\vec V)$ or  $L^2(I;\vec V)$ a further estimate of $E_{\vec q,n}^{0}$ is required which remains an open problem.}

\subsection{Error estimates for the error between the continuous and the fully discrete solution}
\label{Sec:ErrCGMFEM}

In this section we combine the results of Thm.\ \ref{Th:ErrL2NormMix} and  Thm.\ \ref{Thm:ErrNonExNat}
with the estimates of Thm.\ \ref{Thm:L2norm_error_semi2fully} to prove the convergence of the fully discrete
scheme.

\begin{theorem} \label{theorem_error_cont2fully}
Let the assumptions of Subsec.\ \ref{Sec:NotPrem_2} about $\Omega, u_0, \vec D$  and $f$ be 
satisfied. Let $\{u,\vec q\} \in H^1(I;W)\times L^2(I;\vec V)$  denote the unique 
solution of \eqref{Eq:IntMix_1},  \eqref{Eq:IntMix_2} that is supposed 
to be sufficiently regular. Further, let $\{u_{\tau,h},\vec q_{\tau,h}\}\in 
\mathcal{X}^{r}(W_h)\times \mathcal{X}^{r}(\vec V_h)$ be the uniquely
defined solution of the fully discrete problem \eqref{Eq:FulDis1}, \eqref{Eq:FulDis2}, 
respectively. Suppose that the semidiscrete problem \eqref{Eq:SemiDis1}, 
\eqref{Eq:SemiDis2} admits a sufficiently regular solution $\{u_{\tau},\vec 
q_{\tau}\}\in \mathcal{X}^{r}(W)\times \mathcal{X}^{r}(\vec V)$. Then, there holds that
\begin{align}
\| u - u_{\tau, h} \|_{L^2(I; L^2(\Omega))} &\le c (\tau^{r} + h^{p+1})\,. 
\label{theorem_error_cont2fully_1} 
\end{align}
For homogeneous diffusion coefficients $\vec D(\vec x) = d \vec I$, with some constant $d>0$, there holds that
\begin{equation}
\| \vec q - \vec q_{\tau, h} \|_{L^2(I; \vec L^2(\Omega))}  \le  c (\tau^{r} + h^{p+1})\,. 
\label{theorem_error_cont2fully_3}
\end{equation}
Under the regularity condition (R${}_{\mathrm{mix}}$) given in \eqref{Eq:AssmpRegMix} and for interpolated
right-hand side functions \eqref{Eq:DefLagIntRhs} there holds that
\begin{equation}
\label{theorem_error_cont2fully_4}
\| u - u_{\tau, h} \|_{L^2(I; L^2(\Omega))} \le c (\tau^{r+1} + h^{p+1})\,. 
\end{equation}
The constant $c$ in 
\eqref{theorem_error_cont2fully_1}--\eqref{theorem_error_cont2fully_4}, respectively,
does not depend on the discretization parameters $h$ and $\tau$.
\end{theorem}

\begin{mproof} By using the triangle inequality, Thm.\ \ref{Thm:ErrNonExNat}  and Thm.\  \ref{Thm:L2norm_error_semi2fully} it follows that
\begin{align*}
\| u - u_{\tau, h} \|_{L^2(I; L^2(\Omega))}^2 &\le 2 \| u - u_{\tau} \|_{L^2(I; L^2(\Omega))}^2  + 2  \| u_{\tau} - u_{\tau, h} \|_{L^2(I; L^2(\Omega))}^2 \nonumber \\[1ex]
&\le c \left(\tau^{2r} +  h^{2(p+1)}\right)\,,
\end{align*}
where sufficient regularity of the continuous and semidiscrete solution with
appropriate upper bounds for the solutions (cf.\ Thm.\ \ref{Thm:ErrNonExNat} and Thm.\
\ref{Thm:L2norm_error_semi2fully}) is assumed. The
inequality \eqref{theorem_error_cont2fully_3} is obtained similarly. The estimate \eqref{theorem_error_cont2fully_4} can be concluded
in the same way by using now the result of Thm.\ \ref{Th:ErrL2NormMix} instead of Thm.\ 
\ref{Thm:ErrNonExNat}.
\end{mproof}

\begin{rem}
\label{Rem:OptOrd}
\begin{itemize}
\item The error estimate \eqref{theorem_error_cont2fully_4} is optimal in time and space. The assumption of an interpolated
right-hand side function \eqref{Eq:DefLagIntRhs} can still be dropped even though this is not explicitly done in this work.
It requires to estimate the error between the exact form of cGP($r$) defined in \eqref{Eq:GtdP_1}, \eqref{Eq:GtdP_2} and the fully
discrete solution. In this case the arguments used to prove Thm.\  \ref{Thm:L2norm_error_semi2fully} have to be
augmented by an estimate of the interpolation error for the right-hand side function, similarly to the proof of Thm.\
\ref{Thm:ErrNonExNat}.
\item The error estimate \eqref{theorem_error_cont2fully_3} is  suboptimal in time. It remains an open problem to analyze if the estimates can still be sharpened to order $r+1$.
In our numerical study presented in Sec.\ \ref{Sec:NumStudies}  convergence of order $r+1$ will be observed for the temporal discretization of the scalar and the flux variable. Moreover, this is even observed in the (spatially) stronger norm of $L^2(0,T;\vec V)$ instead of $L^2(0,T;\vec L^2(\Omega))$ for the flux variable.
\end{itemize}
\end{rem}

\section{Numerical studies}
\label{Sec:NumStudies}

In this section we present numerical studies in order to illustrate the error
estimate given in Thm.\ \ref{theorem_error_cont2fully} for the fully discrete
scheme \eqref{Eq:FulDis1}, \eqref{Eq:FulDis2} combining a variational time
discretization with the MFEM. Moreover, we analyze the robustness 
of the convergence behaviour with respect to random perturbations of the meshes. Thereby 
we mimic mesh distributions of applications that are of practical interest.    
Additional convergence studies for variational space-time discretizations
of the proposed type as well as for discontinuous time discretizations can be
found in \cite{Bause2015,Koecher2015} for parabolic problems and in
\cite{Koecher2014,Koecher2015} for variational space-time discretizations of
wave equations.
In \cite{Koecher2015,Bause2015} the efficient iterative solution of the resulting
algebraic system of equations \eqref{Eq:FulDis1}, \eqref{Eq:FulDis2} along with
the construction of appropriate preconditioning techniques is carefully addressed.
In the literature, further computational studies of variational time discretization
schemes are presented also for different kind of flow and transport problems in, e.g.,
\cite{Ahmed2012,Ahmed2015,Ahmed2012_2,Schieweck2010,Matthies2011,Hussain2011,
Hussain2013,Hussain2012}.

In order to determine the space-time convergence behavior we consider in our
numerical study the cGP($2$)--MFEM($2$) approach.
That is \eqref{Eq:GtdP_1}--\eqref{Eq:GtdP_2} with $r=2$ combined with the
mixed finite element method MFEM($2$) based on the choice $p=2$ in the definition
\eqref{Eq:DefWh} and \eqref{Eq:DefVh} of the tuple of MFE spaces. We 
prescribe the solution
\begin{displaymath}
% \label{eq:E1:ExactPrimal:ST}
u_{\textnormal{E}}(\vec x, t) :=
\sin(\omega t) \sin(\pi x_1) \sin(\pi x_2)\,, \quad
\textnormal{in}\quad \Omega \times (0,T)\,,
\end{displaymath}
with $\Omega = (0,1)^2$, $\omega = 10\pi$ of problem \eqref{Eq:Diff}--\eqref{Eq:InVal}. 
The corresponding flux function is then given by $\boldsymbol q_{\textnormal{E}} = 
-\boldsymbol D \nabla u_{\textnormal{E}}$ for $\boldsymbol D = \vec I$. We choose the 
final time $T=1$. On the coarsest level (level 0) the temporal mesh is uniformly refined 
into $N=10$ time subintervals and the corresponding spatial mesh consists of a single 
cell. To determine the experimental orders of convergence the space-time mesh is refined 
uniformly by a factor of two in each of the space dimensions and in the time dimension.
We use the abbreviation 
\begin{equation*}
\label{eq:E4:ErrorRepresentation:cG}
e_{u}^{\textnormal{cGP}(2)}(t) :=
u_{\textnormal{E}}(t) -
u_{\tau,h}(t)\quad
\textnormal{and}\quad
e_{{\boldsymbol{q}}}^{\textnormal{cGP}(2)}(t) :=
{\boldsymbol{q}}_{\textnormal{E}}(t) -
{\boldsymbol{q}}_{\tau,h}(t)\,,
\end{equation*}
where we denote by $u_{\tau,h}$ and by ${\boldsymbol{q}}_{\tau,h}$
the fully discrete cGP($2$)--MFEM($2$) approximation of the primal variable and
the flux field.
The discretization errors
for $e_u^{\textnormal{cGP(2)}}$ are measured in the $L^2(I; L^2(\Omega))$-norm and
for $\vec e_{\boldsymbol{q}}^{\textnormal{cGP(2)}}$ in the
$L^2(I; \vec V)$-norm.
As usual, the integral over the spatial domain $\Omega$ and the integral
over the time domain $I=(0,T)$ in the error norms are evaluated elementwise
in space and time by appropriate quadrature rules of sufficiently high order of
accuracy.

\begin{table}[t]
\centering
\begin{tabular}{c||cc|ccc}
Level
& $N$ & $\tau_n$
& $|\mathcal T_h|$ & $h$ & $N_{\textnormal{DoF}}$\\
\hline
\hline
0 &  10 & 1.000e-01 &    1 & 1.4142e-00 &    33 \\
1 &  20 & 5.000e-02 &    4 & 7.0711e-01 &   120 \\
2 &  40 & 2.500e-02 &   16 & 3.5355e-01 &   456 \\
3 &  80 & 1.250e-02 &   64 & 1.7678e-01 &  1776 \\
4 & 160 & 6.250e-03 &  256 & 8.8388e-02 &  7008 \\
5 & 320 & 3.125e-03 & 1024 & 4.4194e-02 & 27840
\end{tabular}
\caption{Space-time mesh with number of time subintervals $N$,
global time discretization parameter $\tau_n$, number of cells $|\mathcal T_h|$,
global space discretization parameter $h$ and
degrees of freedom $N_{\textnormal{DoF}}$ per degree of freedom in time.}
\label{table:T1:ConvSpaceTime:MeshData}
\end{table}

\begin{table}
\centering
\begin{tabular}{c||cc|cc}
Level &%
$ \big\| e_u^{\textnormal{cGP($2$)}} \big\|_{L^2(I; L^2(\Omega))} $ & EOC &%
$ \big\| e_{{\boldsymbol{q}}}^{\textnormal{cGP($2$)}}
\big\|_{L^2(I; \vec V)} $ & EOC \\
\hline
\hline
0 & 4.0298e-02 & ---  & 8.2000e-01 & ---  \\
1 & 1.1316e-02 & 1.83 & 2.2827e-01 & 1.84 \\
2 & 1.4371e-03 & 2.98 & 2.8876e-02 & 2.98 \\
3 & 1.8037e-04 & 2.99 & 3.6208e-03 & 3.00 \\
4 & 2.2569e-05 & 3.00 & 4.5295e-04 & 3.00 \\
5 & 2.8219e-06 & 3.00 & 5.6631e-05 & 3.00
\end{tabular}
\caption{Norm values and corresponding experimental orders of convergence in
space-time for cGP($2$)--MFEM($2$) on the refinement levels as given in Tab.\ 
\ref{table:T1:ConvSpaceTime:MeshData}.}
\label{table:T2:ConvSpaceTime:MFEM2:cG2}
\end{table}

\begin{figure}
\centering
\includegraphics[width=\textwidth]{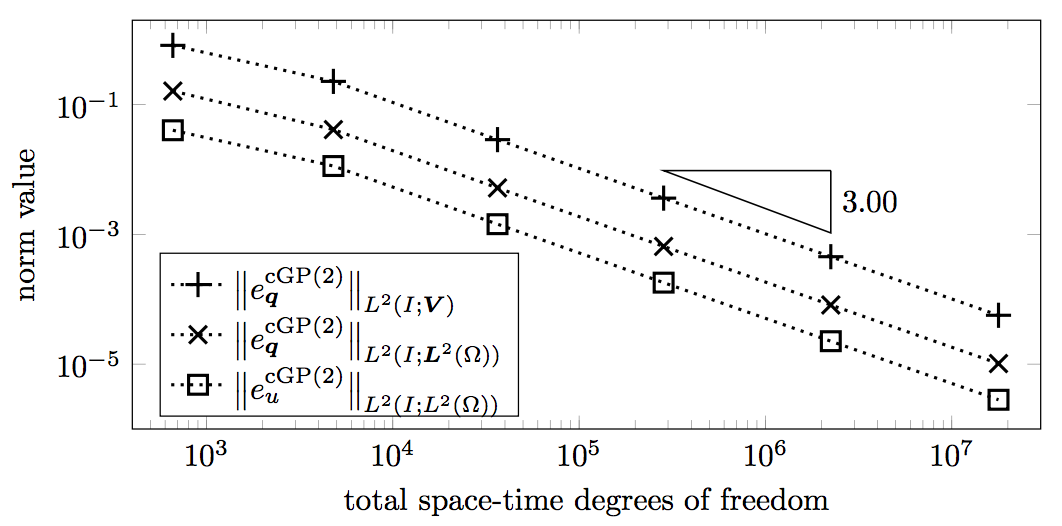}

\caption{Calculated errors and experimental orders of convergence in space-time for 
cGP($2$)--MFEM($2$).
% 
% $N$ denotes the total number of time subintervals and
% $N_{\textnormal{DoF}}$ denotes the total number of degrees of freedom in space.
}
\label{fig:F1:EOC:L2_L2_SpaceTime}
\end{figure}

\subsection{Uniform meshes}
\label{Sec:NumStudies:uniform}

We summarize the calculated errors and their experimental order of convergence
(EOC) for the proposed space-time discretization in Tab.\
\ref{table:T2:ConvSpaceTime:MFEM2:cG2} and further illustrate them in Fig.\
\ref{fig:F1:EOC:L2_L2_SpaceTime}.
The numerical results confirm the expected third order rate of convergence
established in Thm.\ \ref{theorem_error_cont2fully} (cf.\ also Rem.\  \ref{Rem:OptOrd})
for the discretization in the space-time domain with polynomial order
$r=2$ and $p=2$, respectively, in the definition of the underlying finite element
spaces.
We note that the optimal order convergence in time and space is obtained for the
primal and the flux variable. Thus, the estimates \eqref{L2norm_error_semi2fully_2}
and \eqref{L2norm_error_semi2fully_3} might be suboptimal with respect to
the time discretization; cf.\ Rem.\ \ref{Rem:OptOrd}. The estimate 
\eqref{theorem_error_cont2fully_4} is nicely confirmed by the 
presented numerical results. 
Further, we note that the optimal rate of convergence is obtained for the spatial
discretization of the flux field in the norm of $\vec V$.
In this point the family of Raviart--Thomas pairs of mixed finite elements is
superior to the family of Brezzi--Douglas--Marini pairs of mixed finite elements
(cf.\ \cite{Brezzi1991}) for that the optimal order of convergence of the flux
variable can be obtained only in the norm of $\vec L^2(\Omega)$.

%% distorted mesh %%%%%%%%%%%%%%%%%%%%%%%%%%%%%%%%%%%%%%%%%%%%%%%%%%%%%%%%%%%%%%
\subsection{Distorted meshes}
\label{Sec:NumStudies:distorted}
\begin{figure}
\centering

\subfloat[$5\%$]{
\begin{minipage}[c][1\width]{.315\linewidth}
\centering
\resizebox{1.\linewidth}{!}{
\includegraphics[width=\linewidth]
{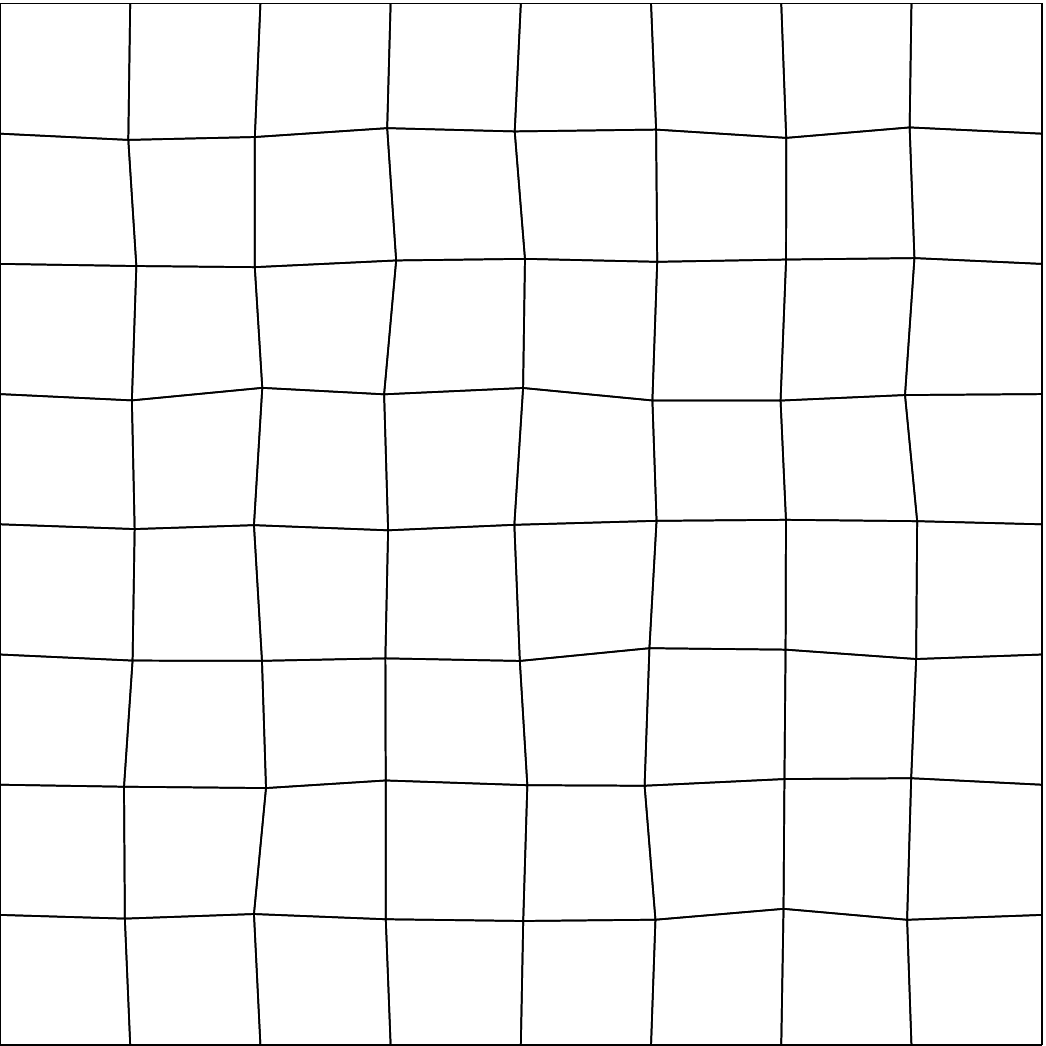}
}
\end{minipage}}
~
\subfloat[$10\%$]{
\begin{minipage}[c][1\width]{.315\linewidth}
\centering
\resizebox{1.\linewidth}{!}{
\includegraphics[width=\linewidth]
{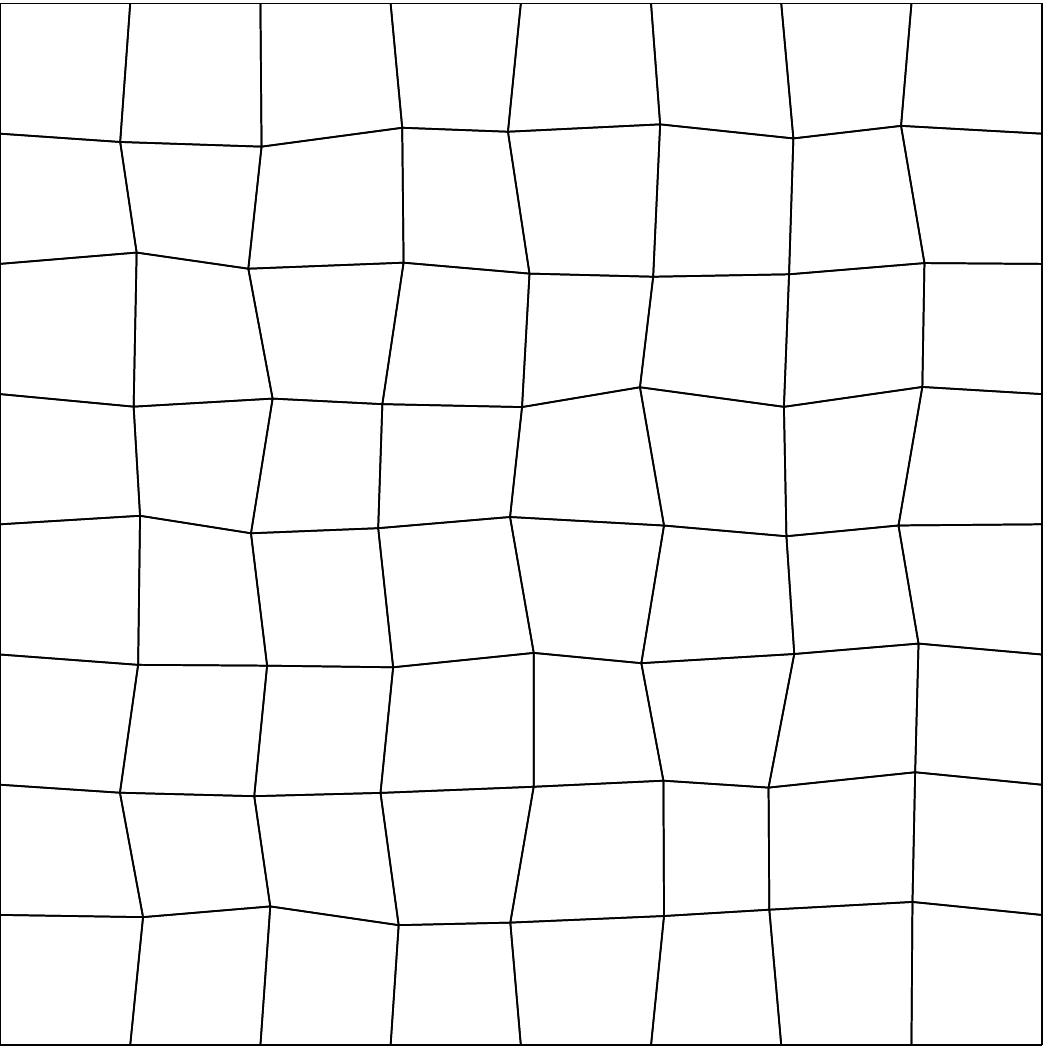}
}
\end{minipage}}
~
\subfloat[$25\%$]{
\begin{minipage}[c][1\width]{.315\linewidth}
\centering
\resizebox{1.\linewidth}{!}{
\includegraphics[width=\linewidth]
{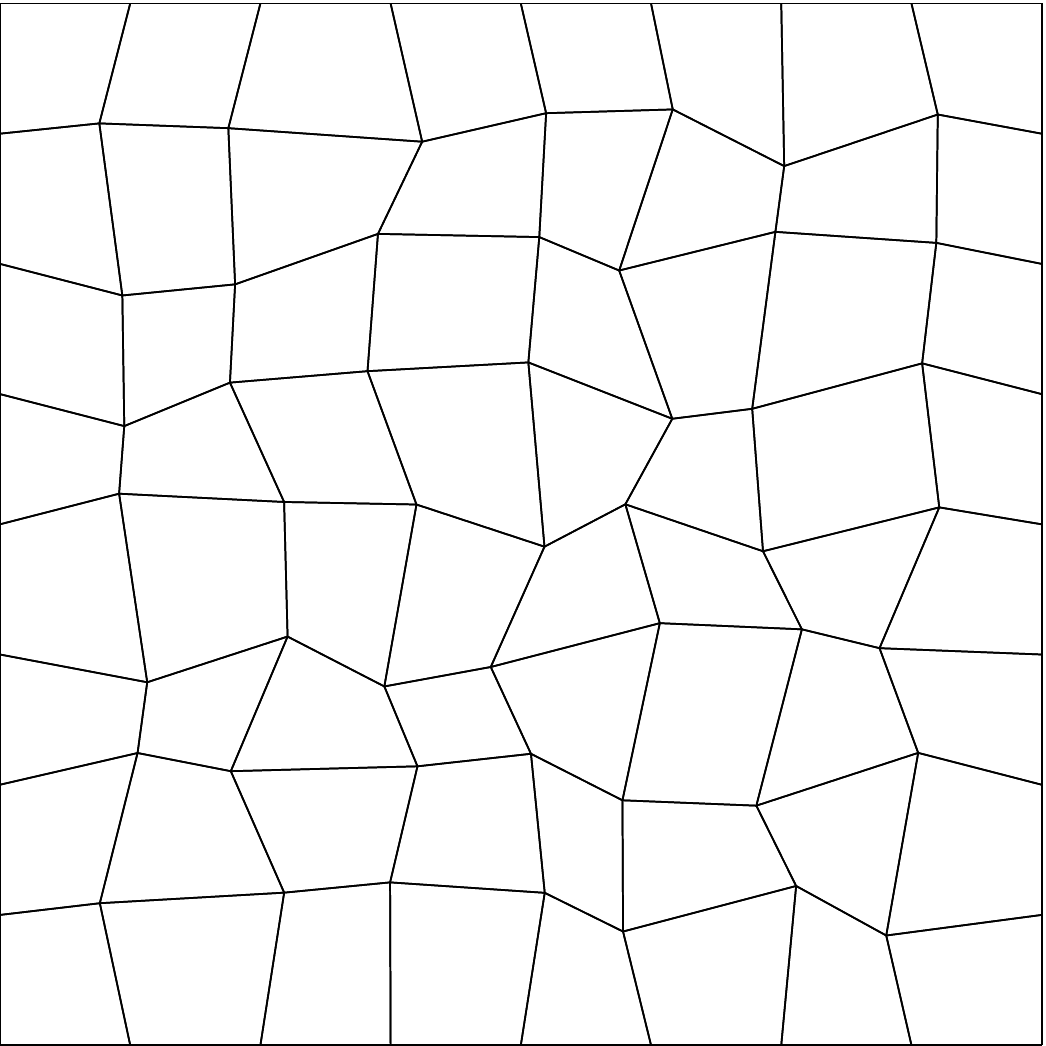}
}
\end{minipage}}

\caption{Distorted spatial meshes for $5\%$ (a), $10\%$ (b) and $25\%$ (c)
random vertex movement for refinement level \MeshRefinement.}

\label{fig:F2:distorted_meshes_fixed_level}
\end{figure}

In the second part of the numerical convergence studies we approximate the same
analytic solution as before but we use spatial meshes with
randomly distorted interior vertices. Precisely, each of the interior vertices is 
distorted by a randomly chosen vector. The magnitude of the distortion vector is chosen
randomly up to a given factor of relative length to the corresponding edge length.
The characteristic numbers of the refinement levels are summarized in Tab.\ 
\ref{table:T3:distorted_mesh_specification}. The resulting distorted meshes 
are illustrated in Fig.\ \ref{fig:F2:distorted_meshes_fixed_level} for the refinement 
level \MeshRefinement.
The temporal mesh is chosen in the same way as in the first numerical experiment;
cf.\ Tab.\ \ref{table:T1:ConvSpaceTime:MeshData}.

\begin{table}
\centering
\begin{tabular}{c||cc|cc|cc|cc}
\multicolumn{1}{c||}{~}
& \multicolumn{2}{c|}{$0 \%$} & \multicolumn{2}{c|}{$5 \%$}
& \multicolumn{2}{c|}{$10 \%$} & \multicolumn{2}{c}{$25 \%$}\\
Level
 & $h_\mathrm{max}$  & $h_\mathrm{black}$ & $h_\mathrm{max}$  & $h_\mathrm{black}$
 & $h_\mathrm{max}$  & $h_\mathrm{black}$ & $h_\mathrm{max}$  & $h_\mathrm{black}$\\
\hline
\hline
0 & 1.4142 &  --- & 1.4142 &  --- & 1.4142 &  --- & 1.4142 &  ---\\
1 & 0.7071 & 2.00 & 0.7258 & 1.95 & 0.7504 & 1.88 & 0.8127 & 1.74\\
2 & 0.3536 & 2.00 & 0.3661 & 1.98 & 0.3974 & 1.89 & 0.4694 & 1.73\\
3 & 0.1768 & 2.00 & 0.1888 & 1.94 & 0.2006 & 1.98 & 0.2393 & 1.96\\
4 & 0.0884 & 2.00 & 0.0946 & 2.00 & 0.1008 & 1.99 & 0.1195 & 2.00
\end{tabular}
\caption{Distorted spatial mesh: $h_\mathrm{max}$ largest cell diameter and
$h_\mathrm{black}$ cell diameter reduction factor
for $0 \%$, $5 \%$, $10 \%$ and $25 \%$ random vertex movement.}
\label{table:T3:distorted_mesh_specification}
\end{table}

\begin{table}
\centering
\begin{tabular}{c||cc|cc|cc}
Level &%
$ 5\% $ & EOC &%
$ 10\% $ & EOC &%
$ 25\% $ & EOC \\
\hline
\hline
0 & 4.0298e-02 & ---  & 4.0298e-02 & ---  & 4.0298e-02 & ---  \\
1 & 1.1353e-02 & 1.83 & 1.1463e-02 & 1.81 & 1.2155e-02 & 1.73 \\
2 & 1.4381e-03 & 2.98 & 1.4783e-03 & 2.95 & 1.8690e-03 & 2.70 \\
3 & 1.8290e-04 & 2.97 & 1.9117e-04 & 2.95 & 2.5232e-04 & 2.89 \\
4 & 2.3024e-05 & 2.99 & 2.4402e-05 & 2.97 & 3.4463e-05 & 2.87
\end{tabular}
\caption{Calculated errors and corresponding experimental order of convergence
for $\big\| e_u^{\textnormal{cGP($2$)}} \big\|_{L^2(I; L^2(\Omega))}$
on distorted meshes given in Tab. \ref{table:T3:distorted_mesh_specification}.}
\label{table:T4:ConvSpaceTime:MFEM2:cG2:distorted:u}
\end{table}

\begin{table}
\centering
\begin{tabular}{c||cc|cc|cc}
Level &%
$ 5\% $ & EOC &%
$ 10\% $ & EOC &%
$ 25\% $ & EOC \\
\hline
\hline
0 & 8.2000e-01 & ---  & 8.2000e-01 & ---  & 8.2000e-01 & ---  \\
1 & 2.2964e-01 & 1.84 & 2.3376e-01 & 1.81 & 2.6241e-01 & 1.64 \\
2 & 3.0172e-02 & 2.93 & 3.5136e-02 & 2.73 & 7.5321e-02 & 1.80 \\
3 & 3.9244e-03 & 2.94 & 4.7838e-03 & 2.88 & 1.0391e-02 & 2.86 \\
4 & 6.0990e-04 & 2.69 & 9.6250e-04 & 2.31 & 2.7408e-03 & 1.92
\end{tabular}
\caption{Calculated errors and corresponding experimental order of convergence
for $\big\| e_{{\boldsymbol{q}}}^{\textnormal{cGP($2$)}} \big\|_{L^2(I; \vec V)}$
on distorted meshes given in Tab. \ref{table:T3:distorted_mesh_specification}.}
\label{table:T5:ConvSpaceTime:MFEM2:cG2:distorted:q:Hdiv}
\end{table}

\begin{figure}[t!]
\centering
%%%%%%%%%%%%%%%%%%%%%%%%%%%%%%%%%%%%%%%%%%%%%%%%%%%%%%%%%%%%%%%%%%%%%%%%%%%%%%%%
\includegraphics[width=\textwidth]{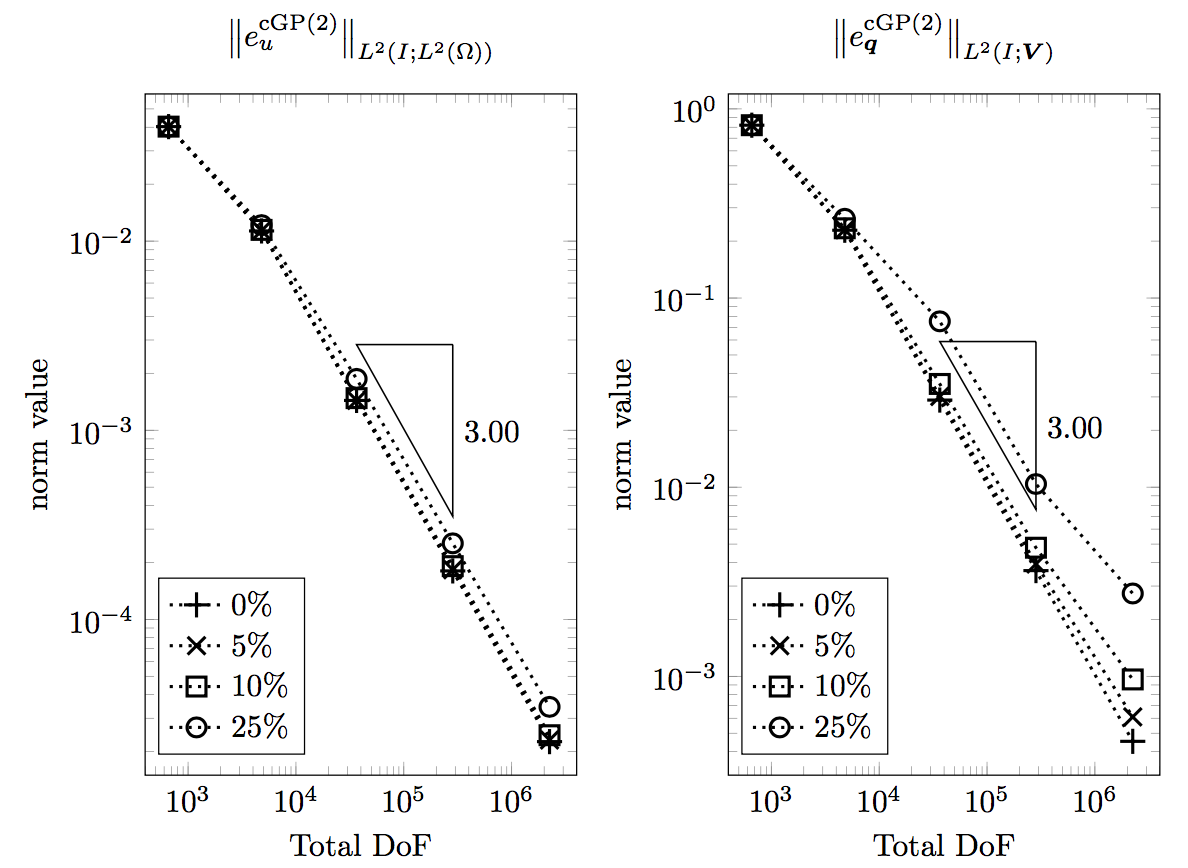}
\caption{Calculated errors and corresponding experimental order of convergence on 
distorted meshes given in Tab. \ref{table:T3:distorted_mesh_specification}.}
% total space-time DoF: ($2 \cdot N \cdot N_{\textnormal{DoF}}$)
\label{fig:F3:EOC:distorted_meshes}
\end{figure}

We summarize the calculated errors and the corresponding experimental order of 
convergence (EOC) for the proposed space-time discretization on the distorted spatial 
meshes in Tab.\ \ref{table:T4:ConvSpaceTime:MFEM2:cG2:distorted:u} for the scalar-valued
primal variable and in Tab.\ \ref{table:T5:ConvSpaceTime:MFEM2:cG2:distorted:q:Hdiv}
for the vector-valued flux variable and further illustrate them in
Fig.\ \ref{fig:F3:EOC:distorted_meshes}.
The expected experimental order of convergence in space and time, for the primal variable measured in the $L^2(I;L^2(\Omega))$-norm and for the flux variable in the $L^2(I; \vec V)$-norm, is largely confirmed even for the strongly perturbed meshes with a distortion factor of 25\,\%. This nicely demonstrates the  robustness of the numerical scheme. We note that the space-time convergence studies on the 
distorted spatial meshes were done with exactly the same numerical solver settings as for the above-given studies on uniform meshes.

\section{Conclusions}
\label{Sec:Conclusions}

In this work a numerical analysis of a family of variational space approximation schemes
that combine continuous finite elements in time with the MFEM in space was presented for 
a parabolic prototype model of flow in porous media. The existence and uniqueness of the 
temporally semidiscrete and the fully discrete approximations were
proved. Error estimates with explicit rates of convergence, including an optimal order 
error estimate, in natural norms of the scheme were established. The error estimates were 
illustrated and confirmed by numerical convergence studies. We believe that our analyses 
and techniques can be extended and applied to more sophisticated flow and transport 
processes in porous media or to incompressible viscous free flow. This will be our work 
for the future.

\section*{Acknowledgements}

This work was supported by the German Academic Exchange Service (DAAD) under
the grant IDs 56435737 and 57238185 and by the Research Council of Norway under the grant DAADppp225267 and DAADppp255715.

The authors wish to thank the anonymous reviewers for their help to improve the 
presentation of this paper.

\appendix
\numberwithin{equation}{section}
\section{Supplementary proofs}

In the sequel we introduce a variational semidiscretization in time of the weak
formulation of the second order problem \eqref{Eq:Diff}--\eqref{Eq:InVal}, i.e.\ without
rewriting Eq.\ \eqref{Eq:Diff} as a first order system of equations as it is done in
Subsec.\ \ref{Sec:SemiTimeCont}. Then we prove the existence and uniqueness of solutions
to the resulting semidiscrete variational problem. This results is used to establish the
existence of the semidiscrete approximation in mixed form defined by the variational
problem \eqref{Eq:GtdP_1}, \eqref{Eq:GtdP_2} in Sec.\ \ref{Sec:ExUniErr}.
Here, we present a different technique of proof than in \cite{Schieweck2010} since one of
the arguments that is used \cite[Lemma 6.1]{Schieweck2010} does not hold in the
applied form from our point of view. Thereby, we aim to keep our work self-contained.
Further, we summarize the proof of Thm.\ \ref{Thm:ErrNonExNat}.

\subsection{Variational time discretization of the second order problem}

In the following we use the notation that is introduced in Subsec.\ \ref{Sec:NotPrem_1}
and \ref{Sec:SemiTimeCont}, respectively. Moreover we use the splitting (cf.\ Eq.\
\eqref{Eq:Split})
\begin{equation*}
\label{Eq:Split2nd}
u_\tau(t) = u_0 + u^0_\tau(t) \qquad \text{with}\quad u^0_\tau \in
\mathcal{X}_0^{r}(H^1_0(\Omega))\,.
\end{equation*}
Further, we put
\[
f_0 (t) = f(t) - A u_0
\]
for $u_0 \in H^1_0(\Omega)$ such that by definition $Au_0 \in H^{-1}(\Omega)$, cf.\ Sec.\
\ref{Sec:NotPrem_1}. Under the additional regularity condition that
\begin{equation*}
\label{Eq:RegCondH2}
D(A)=H^2(\Omega)\cap H^1_0(\Omega)
\end{equation*}
it even holds that $f_0\in \braum$. The semidiscrete variational
approximation of the system \eqref{Eq:Diff}--\eqref{Eq:InVal} is now defined by:
\textit{Find $u_\tau^0 \in X^r_0 (H^1_0(\Omega))$ such that}
\begin{equation}
\label{Eq:SemiDis2nd}
 \int_0^T \langle \partial_t u_\tau^0, w_\tau \rangle \ud t + \int_0^T
a(u_\tau^0,w_\tau) \ud t = \int_0^T \langle f_0, w_\tau\rangle \ud   t
\end{equation}
\textit{for all $w_\tau \in Y^{r-1}(H^1_0(\Omega))$.}

Firstly, we show the uniqueness of solutions to \eqref{Eq:SemiDis2nd}. In the sequel we 
denote by $\varphi_{n,j}=\varphi_{n,j}(t)$ for $j=0,\ldots,r$ the Lagrange basis 
functions in $I_n=(t_{n-1},t_n]$ with respect to $r+1$
quadrature points $t_{n,l}$, $l=0,\ldots,r$. Here, we choose the Gauss-Lobatto 
quadrature rule that is exact for polynomials of maximum degree $2r-1$. In particular, 
for the quadrature nodes in $\overline{I_n}$ it holds that $t_{n,0}=t_{n-1}$ and 
$t_{n,r}=t_n$. Then, any function $u_\tau^0 \in X^r_0
(H^1_0(\Omega))$ and its time derivative admit the representation
\begin{equation}
\label{Eq:RepSemiDis2nd}
u_\tau^0 (t) =  \sum_{j=0}^r U_n^j  \varphi_{n,j}(t)\,,\qquad
\partial_t u_\tau^0 (t) =  \sum_{j=0}^r U_n^j  \varphi_{n,j}^\prime (t)
\end{equation}
for all $t\in I_n$ with coefficient functions $U_n^j \in H^1_0(\Omega)$ for
$j=0,\ldots ,r$.

\begin{theorem}[Uniqueness of solutions to \eqref{Eq:SemiDis2nd}]
\label{Thm:Uni2nd}
Let the assumptions of Subsec.\ \ref{Sec:NotPrem_2} about $\Omega, u_0$ and $f$ be 
satisfied. Then the solution $u_\tau^0\in \mathcal{X}^{r}_0(H^1_0(\Omega))$ of
the semidiscrete problem \eqref{Eq:SemiDis2nd} is unique.
\end{theorem}

\begin{mproof}
Let $u_{\tau,1}^{0}$, $u_{\tau,2}^{0}\in \mathcal{X}^r_0{(H^1_0(\Omega))}$ denote two 
solutions of the semidiscrete variational problem \eqref{Eq:SemiDis2nd}. We put 
$u_\tau^0(t):= u_{\tau,1}^{0} - u_{\tau,2}^{0}$.
We choose the test function $w_\tau := A^{-1} \partial_t u_\tau^0 + \mu \partial_t
u_\tau^0$ for some fixed parameter $\mu \geq 0$. By means of \eqref{Eq:RepSemiDis2nd}, it
holds that $w_\tau \in \mathcal{Y}^{r-1}{(H^1_0(\Omega))}$. For this choice of $w_\tau$
it follows that
\begin{equation}
 \label{Eq:A_Uni_1}
 I:= \int_0^T \langle \partial_t u_\tau^0(t), A^{-1} \partial_t
u_\tau^0  + \mu \partial_t u_\tau^0 \rangle \ud t +\int_0^T
\langle A u_\tau^0(t), A^{-1} \partial_t u_\tau^0 + \mu \partial_t u_\tau^0
\rangle \ud t  = 0 \,.
\end{equation}
By the symmetry of $a(\cdot,\cdot)$ we have that $a(u_\tau^0,\partial_t
u_\tau^0) = \frac{1}{2\,}\frac{d}{dt} a(u_\tau^0 , u_\tau^0 )$. Further, we have that
$\langle u_\tau^0, \partial_t u_\tau^0 \rangle = \frac{1}{2\,}\frac{d}{dt}
\| u_\tau^0 \|^2_{L^2(\Omega)}$.  Recalling \eqref{Eq:Coerc_1}--\eqref{Eq:Prop_A}
and noting that $u_\tau^0(0)=0$ and $\partial_t u_\tau^0\in
\mathcal{Y}^{r-1}{(H^1_0(\Omega))}$, it follows from \eqref{Eq:A_Uni_1} that
\begin{align*}
 0 = I & =  \int_0^T \langle \partial_t u_\tau^0(t), A^{-1} \partial_t
u_\tau^0 \rangle \ud t + \int_0^T \langle \partial_t
u_\tau^0(t),  \mu \partial_t u_\tau^0 \rangle \ud t\\[1ex]
& \qquad + \int_0^T
\langle A u_\tau(t), A^{-1} \partial_t u_\tau^0 \rangle \ud t + \int_0^T
\langle A u_\tau^0 (t),  \mu \partial_t u_\tau^0 \rangle \ud t\\[2ex]
& \geq c \int_0^T \| \partial_t u_\tau^0 \|_{H^{-1}( \Omega)}^2 \ud t + \mu \int_0^T
\| \partial_t u_\tau^0 \|_{L^2(\Omega)}^2 \ud t \\[1ex]
& \qquad + \int_0^T \frac{1}{2}\frac{d}{dt} \| u_\tau^0 \|_{L^2(\Omega)}^2 \ud t + \mu
\int_0^T \frac{1}{2}\frac{d}{dt} a(u_\tau^0,u_\tau^0) \ud t \\[2ex]
& \geq c \int_0^T \| \partial_t u_\tau^0 \|_{H^{-1}(\Omega)}^2 \ud t + \mu \int_0^T
\| \partial_t u_\tau^0 \|_{L^2(\Omega)}^2 \ud t\\[1ex]
& \qquad + \frac{1}{2}\, \|u_\tau^0
(T)\|_{L^2(\Omega)}^2 +
\frac{\alpha\mu }{2} \| u_\tau^0(T)\|^2_{H^1_0(\Omega)}\,.
\end{align*}
This implies that $u_\tau^0 =0$ and, consequently, that $u_{\tau,1}^{0} =
u_{\tau,2}^{0}$. The uniqueness of solutions to \eqref{Eq:SemiDis2nd} is thus
established.
\end{mproof}

We remark that testing Eq.\ \eqref{Eq:SemiDis2nd} with $v_\tau= A^{-1}\partial_t
u_\tau^0$ or $v_\tau = \partial_\tau u_\tau^0$ would already be sufficient for proving
the uniqueness result. Further, the symmetry of $a(\cdot,\cdot)$ is essential in the
previous proof. A generalization of the arguments to problems with nonsymmetric
bilinearforms, for instance to convection-diffusion equations, still remains an open
problem.

The existence of a solution to the semidiscrete problem \eqref{Eq:SemiDis2nd} follows
from the uniqueness of the solutions. Using the eigenspaces of $A$, problem
\eqref{Eq:SemiDis2nd} can be reduced to a set of finite dimensional problems, for each of
which obviously uniqueness implies existence. For this we recall the following result
from \cite[Appendix D.6]{Evans2010}.

\begin{lemm}
\label{Lem:Eigfunc}
Let $H$ be a separable Hilbert space, and suppose that $S: H\mapsto H$ is a compact and
symmetric operator. Then there exists a countable orthonormal basis of $H$ consisting of
eigenfunctions of $S$.
\end{lemm}

\begin{theorem}[Existence of solutions to \eqref{Eq:SemiDis2nd}]
\label{Thm:Ex2nd}
Let the assumptions of Subsec.\ \ref{Sec:NotPrem_2} about $\Omega, u_0$ and $f$ as the
be satisfied. Then the semidiscrete problem \eqref{Eq:SemiDis2nd} admits a solution
$u_\tau^0\in \mathcal{X}^{r}_0(H^1_0(\Omega))$.
\end{theorem}

\begin{mproof}
The operator $S:=A^{-1}: L^2(\Omega)\mapsto L^2(\Omega)$ with $A$ being defined in
\eqref{Eq:Def_A} is a bounded, linear compact operator mapping $L^2(\Omega)$
into itself. By means of Lemma \ref{Lem:Eigfunc} there exists a set of appropriately
scaled eigenfunctions $\{w_k\}_{k=1}^\infty\subset L^2(\Omega)$ with $w_k\in
H^1_0(\Omega)$ such that $ \{w_k\}_{k=1}^\infty$ is an orthogonal basis of $H^1_0(\Omega)$
and an orthonormal basis of $L^2(\Omega)$.
%; cf.\ \cite[Sec.\ 6.6 and 7.1]{Evans2010}.

In terms of these eigenfunctions $\{w_k\}_{k=1}^\infty\subset H^1_0(\Omega)$
the solution $u_\tau^0$ of problem \eqref{Eq:SemiDis2nd} can be represented as
\begin{displaymath}
 u_\tau(x,t) = \sum_{j=0}^r U_{n}^{(j)}(x) \varphi_{n}^{(j)}(t) =
\sum_{j=0}^r\sum_{k=1}^\infty d_{n,k}^{(j)} w_k(x) \varphi_{n}^{(j)}(t)\,,
\quad \mbox{for } t\in \overline{I}_{n}\,.
\end{displaymath}
with coefficients $d_{n,k}^{(j)}\in \R$ for $k=1,\ldots ,\infty$ and each $j=0,\ldots ,r$
and $n=0,\ldots , N$. We choose test functions $v_\tau 
\in \mathcal{Y}^{r-1}(H^1_0(\Omega))$ being defined by
\begin{displaymath}
 v_\tau = \left\{
 \begin{array}{@{}ll}
  w_k \psi_{n}^{(i)}\,, & \mbox{for } t \in \overline{I}_n\,,\\[2ex]
  0 \,, & \mbox{for } t \in I \backslash \overline{I}_n
 \end{array}
 \right.
\end{displaymath}
for $i=1,\ldots , r$, $k=1,\ldots, \infty$ and $n=1,\ldots,N$. Then, for each
$k=1,\ldots, \infty$, we get the finite dimensional problem
\begin{equation}
\label{Eq:FinDimProb}
\begin{split}
& \sum_{j=0}^r d_{n,k}^{(j)} \underbrace{\int_{I_n} d_t \varphi_{n}^{(j)}(t)
\psi_{n}^{(i)}(t) \ud t}_{:=\alpha_{ij}} +\sum_{j=0}^r d_{n,k}^{(j)}
\underbrace{a(w_k,w_k)}_{:=\gamma_k} \underbrace{\int_{I_n}  \varphi_{n}^{(j)}(t)
\psi_{n}^{(i)}(t) \ud t}_{=:\beta_{ij}} \\[1ex]
& = \underbrace{\int_0^T \langle f_0(t),w_k \rangle
\psi_n^{(i)}(t) \ud t}_{=: b_{k,i}}
\end{split}
\end{equation}
for $i=1,\ldots r$ and $n=1,\ldots,N$. Due to the continuity of functions $u_\tau \in
\mathcal{X}^r_0{(H^1_0(\Omega))}$ and the choice of the Gauss-Lobatto quadrature rule it 
holds that
\begin{displaymath}
 U_n^{(0)} = U_{n-1}^{(r)} \qquad \mbox{or} \qquad  d_{n,k}^{(0)} =
d_{n-1,k}^{(r)}\,,
\end{displaymath}
respectively. Therefore, we recast the finite dimensional problem \eqref{Eq:FinDimProb} as
\begin{equation}
\label{Eq:FinDimProb2}
\sum_{j=1}^r \left(\alpha_{ij} + \gamma_k \beta_{ij}\right) d_{n,k}^{(j)} =
b_{k,i} + \left(\alpha_{i0} + \gamma_k \beta_{i0}\right) d_{n,k}^{(0)}
\end{equation}
for $i=1,\ldots r$ and each $k=1,\ldots ,\infty$. For the finite dimensional problem
\eqref{Eq:FinDimProb2} the uniqueness of the solution established in Thm.\
\ref{Thm:Uni2nd} then implies the existence of a solution.
\end{mproof}

Similarly to Corollary \ref{Thm:InfSup}, the existence and uniqueness of the semidiscrete
solution implies that an inf-sup stability condition in the underlying space-time
framework is satisfied \cite[p.\ 85, Thm.\ 2.6]{Ern2010}.

\begin{cor}
\label{Thm:InfSup2nd}
Let the assumptions of Subsec.\ \ref{Sec:NotPrem_2} about $\Omega, u_0, \vec D$
and $f$ be satisfied. Let $u_\tau^0\in \mathcal{X}^{r}_0(H^1_0(\Omega))$ be the unique
solution of the semidiscrete problem \eqref{Eq:SemiDis2nd} according to Thm.\
\ref{Thm:Uni2nd} and \ref{Thm:Ex2nd}. Then, there exists a constant $c > 0$ such that
\begin{equation*}
 \inf_{u_\tau \in \mathcal{X}^r_0{(H^1_0(\Omega))}\backslash \{0\}} \sup _{v_\tau \in
\mathcal{Y}^{r-1}{(H^1_0(\Omega))}\backslash \{0\}}
\dfrac{B(u_\tau,v_\tau)}{\|u_\tau\|_{\mathcal{X}}\|v_\tau\|_{\mathcal{Y}}} \geq
c  > 0
\end{equation*}
with
\begin{align}
\nonumber
 B(u_\tau,v_\tau) & = \int_0^T \langle \partial_t u_\tau + A
u_\tau ,v_\tau\rangle \ud t\,, \\[1ex]
\|u_\tau\|_{\mathcal{X}} & = \left(\|u_\tau\|^2_{L^2(I,H^1_0(\Omega))} +
\|\partial_t u_\tau\|^2_{L^2(I,H^{-1}(\Omega))}\right)^{1/2}\,,
\label{Eq:NatNorm2nd}\\[2ex]
\nonumber
\|v_\tau\|_{\mathcal{Y}} & = \|u_\tau\|_{L^2(I,H^1_0(\Omega))}\,.
\end{align}
\end{cor}

\begin{rem}
By the arguments of \cite[Thm.\ 6.2]{Schieweck2010} the inf-sup stability condition
implies an error estimate for the semidiscretization \eqref{Eq:SemiDis2nd} where
the error is measured in the corresponding natural norm \eqref{Eq:NatNorm2nd} of
the scheme.
\end{rem}

\subsection{Proof of Thm.\ \ref{Thm:ErrNonExNat}}

\begin{mproof}
We let $w_\tau := I_\tau u - u_\tau$,
$\vec v_\tau = \vec J_\tau \vec q - \vec q_\tau $ with the interpolation operators 
$I_\tau$
and $\vec J_\tau$ of Subsec.\ \ref{Sec:ErrExctFrm}. Using the inf-sup stability
condition along with problem \eqref{Eq:IntMix_1}, \eqref{Eq:IntMix_2} and the non-exact 
semidiscrete problem \eqref{Eq:NonExGtdP_1}, \eqref{Eq:NonExGtdP_2}, applying the 
inequality of Cauchy--Schwarz and the continuity \eqref{Eq:ContA} of 
$a_\tau(\cdot,\cdot)$ we conclude that
\begin{align*}
& \alpha \| \{w_\tau,\vec v_\tau \} \|_{\mathcal{W}} \| \{\varphi_\tau \vec
\psi_\tau\} \|_{\mathcal{V}}  \leq a_\tau(\{w_\tau,\vec
v_\tau\},\{\varphi_\tau,\vec \psi_\tau\})\\[1ex]
& = a_\tau(\{u^0,\vec q\},\{\varphi_\tau,\vec \psi_\tau\})- a_\tau(\{u^0_\tau,\vec
q_\tau\},\{\varphi_\tau,\vec \psi_\tau\})\\[1ex]
&  \quad - a_\tau(\{u^0- I_\tau u^0_\tau,\vec q -\Vec J_\tau \vec q\},\{\varphi_\tau,\vec
\psi_\tau\})\\[1ex]
& = \int_0^T \langle f(t)- \Pi_r f(t),\varphi_\tau \rangle \ud t -a_\tau(\{u^0-
I_\tau u^0_\tau,\vec q -\vec J_\tau \vec q\},\{\varphi_\tau,\vec
\psi_\tau\}) \\[1ex]
& \leq \|f- \Pi_r f\|_{L^2(I;W)} \| \{\varphi_\tau,\vec
\psi_\tau\}\|_{\mathcal{V}} + c \|\{u^0-I_\tau u^0,\vec q - \vec J_\tau \vec q
\}\|_{\mathcal{W}} \|
\{\varphi_\tau,\vec \psi_\tau\}\|_{\mathcal{V}}
\\[2ex]
& \leq \Bigg(\sum_{n=1}^N \Big\{\|f- \Pi_r f\|_{L^2(I;W)}^2 + c \|u^0 -
I_\tau u^0 \|^2_{L^2(I_n;W)}\\[0ex]
&  \qquad + c \|\partial_t (u^0 - I_\tau u^0) \|^2_{L^2(I_n;W)} + c \|\vec q -
\vec
J_\tau \vec
q\|^2_{L^2(I_n;\vec V)} \Big\}\Bigg)^{1/2} \cdot \| \{\varphi_\tau,\vec
\psi_\tau\}\|_{\mathcal{V}}\,.
\end{align*}
By means of the approximation properties \eqref{Eq:Intpol1} to
\eqref{Eq:Intpol3} we then get that
\begin{equation*}
\begin{split}
\| \{w_\tau,\vec v_\tau \} \|_{\mathcal{W}} \leq c \Bigg(\sum_{n=1}^N
 & \tau_n^{2r} \Big\{\|\partial_t^{r+1} u\|^2_{L^2(I_n;W)} \\[-1ex]
 & + \tau_n^2 \|\partial_t^{r+1}
\vec q \|^2_{L^2(I_n;\vec V)} + \tau_n^2 \|\partial_t^{r+1} f
\|^2_{L^2(I_n;W)}\Big\}\Bigg)^{1/2}\,.
\end{split}
\end{equation*}
Combining this estimate with the triangle inequality yields the assertion of Thm.\
\ref{Thm:ErrNonExNat}. \linebreak[4] \mbox{}
\end{mproof}

\section{Summary of notation}

\subsection*{Function spaces and norms}
\begin{center}
\small
\begin{tabular}{|p{5cm}|p{7cm}|}
\hline
$I$, $\Omega$ & time and space domain, $I=(0,T]$\\
$L^2(\Omega)$, $H^p(\Omega)$, $H_0^1(\Omega)$, $W^{1,\infty}(\Omega)$ &  standard Sobolev spaces\\ 
$H^{-1}(\Omega)$ & dual space of $H_0^1(\Omega)$\\
$\| \cdot \|$, $\| \cdot \|_p$, $\langle \cdot, \cdot \rangle$ & Norms in $L^2(\Omega)$, $H^p(\Omega)$, inner product in $L^2(\Omega)$ \\ 
$W$, $\vec V$ & $W=L^2(\Omega)$, $\vec V = \vec H(\text{div};\Omega)$\\
$\| \cdot \|_{\vec V}$ & $\| \cdot \|_{\vec V} = (\| \cdot \|^2+ \| \nabla \cdot (\cdot) \|^2)^{1/2}$\\
$C(\overline I;X)$, $L^2(I;X)$, $H^1(I;X)$ & Bogner spaces with values in the Banach space $X$\\
$\vec D$, $u_0$, $f$ & data of the problem\\
$u$, $\vec q$ & solution and flux of the continuous problem \\
$a(u,v)$ & bilinear form $a(u,v)= \langle \vec D u, v\rangle$\\
\hline 
\end{tabular}
\end{center}

\subsection*{Time discretization} 
\begin{center}
\small
\begin{tabular}{|p{5cm}|p{7cm}|}
\hline
$\mathcal{X}^{r}{(W)}$, $\mathcal{X}^{r}{(\vec V)}$ & trial spaces of semidiscretization in time\\
$\mathcal{X}_0^{r}{(W)}$ & trial space of semidiscretization with  $u_\tau(0)=0$\\
$\mathcal{Y}^{r-1}{(W)}$, $\mathcal{Y}^{r}{(\vec V)}$ & test spaces  of semidiscretization in time\\
$\P_r(J;\,X)$ & $\P_r(J;\,X) = \{p(t) = \sum_{j=0}^{r}{\xi_n^j\, t^j} \mid \xi_n^j \in X\}$\\
$\mathcal{W}$, $\mathcal{V}$ & $\mathcal{W} = X_0^r (W) \times X^r (\vec V)$, $ \mathcal{V} = Y^{r-1} (W) \times Y^{r-1} (\vec V)$\\
$I_n$, $\tau_n$ & subinterval $I_n=(t_{n-1},t_n]$, $\tau_n = t_n-t_{n-1}$\\
$t_{n,0}$ & $t_{n,0}=t_{n-1}$\\
$t_{n,1},\ldots, t_{n,r}$ & $r$-point Gaussian  quadrature nodes\\
$ \varphi_{n,j}(t)$ & Lagrange polynomial on $I_n$ w.r.t $t_{n,0}, \ldots, t_{n,r}$\\
$ \psi_{n,i}(t)$ & Lagrange polynomial on $I_n$ w.r.t $t_{n,1}, \ldots, t_{n,r}$\\
$u_\tau$  & $u_\tau{}_{|\overline{I}_n} (t ) = \sum_{j=0}^r U_n^j\,  \varphi_{n,j}(t) $\\ & semidiscrete solution\\
$\vec q_\tau$ & $\vec q_\tau{}_{|\overline{I}_n} (t ) = \sum_{j=0}^r \vec Q_n^j \, \varphi_{n,j}(t)$ \\ & semidiscrete flux\\
$\hat I$ & $\hat I = [0,1]$ reference interval\\
$\hat t_{1},\ldots , \hat t_{r}$ & Gaussian quadrature nodes on $\hat I$\\
$\hat \omega_{1},\ldots , \hat \omega_{r}$ & Gaussian quadrature weights on $\hat I$\\
$\hat \varphi_j(\hat t)$ & transformed Lagrange polynomial on $\hat I$ \\
$\hat \alpha_{ij}$, $\hat \beta_{ij}$ & $\hat \alpha_{ij} = \hat \omega_i$, $\hat \beta_{ij} = \hat \omega_i \delta_{i,j}$ with Kronecker symbol $\delta_{i,j}$;\\
& from transformation of time integrals to $\hat I$ and\\
& application of quadrature on $\hat I$\\
$\Pi_r$ & temporal Lagrange interpolant\\
&  w.r.t.\ to $t_{n,0},\ldots, t_{n,r}$\\
$ I_\tau$ & temporal interpolation operator for variable $u$\\
$ \vec J_\tau$ & temporal interpolation operator for variable $\vec q$\\
$a_{\tau}(\{\cdot,\cdot\},\{\cdot,\cdot\})$ & space-time bilinear form\\
$u^0(t)$, $u_\tau^0(t)$ & $u^0(t) = u(t)-u_0$, $u_\tau^0(t) = u_\tau(t)-u_0$  \\ 
$z$, $\vec p$ & solution of dual problem\\
$ I_0$ & temporal interpolation operator for $z$\\
$ \vec J_0$ & temporal interpolation operator for $\vec p$\\ \hline
\end{tabular}
\end{center}

\subsection*{Space discretization and error analysis} 
\begin{center}
\small
\begin{tabular}{|p{5cm}|p{7cm}|}
\hline
$W_h$, $\vec V_h$ & inf-sup stable pair of finite element spaces\\
& Raviart--Thomas(--N\'ed\'elec) elements\\
$\mathcal{X}^{r}{(W_h)}$, $\mathcal{X}^{r}{(\vec V_h)}$ & trial spaces of space-time discretization \\
$\mathcal{Y}^{r-1}{(W_h)}$, $\mathcal{Y}^{r-1}{(\vec V_h)}$ & test spaces of space-time discretization\\
$u_{\tau,h}$ & $u_{\tau,h} (t)_{|I_n} = \sum_{j=0}^r U_{n,h}^{j} \varphi_{n,j}(t) $ \\ & fully discrete solution\\
$\vec q_{\tau,h}$ &  $\vec q_{\tau,h} (t)_{|I_n} =  \sum_{j=0}^r \vec
Q_{n,h}^{j} \varphi_{n,j}(t)$ \\ & fully discrete flux\\
$P_h$, $\vec P_h$ & $L^2$ projection onto $W_h$ and $\vec V_h$\\
$\vec \Pi_h$ & projection onto $\vec V_h$:\\ & $\la \nabla \cdot (\vec \Pi_h \vec v - \vec v) , w_h \ra = 0$ for all $w_h\in W_h$ \\  
$E_{u}(t)$ & $E_{u}(t)  =  u_\tau (t) -  u_{\tau,h} (t)$\\ 
& $E_{u}(t)=\sum_{j=0}^r E_{u,n}^j \varphi_{n,j}(t) $\\ & error between semidiscrete and fully discrete\\& approximation of scalar variable\\
$\vec E_{\vec q}(t)$ & $\vec E_{\vec q}(t) =\vec q_\tau (t) -
\vec q_{\tau,h}(t)$ \\ 
& $\vec E_{\vec q}(t) = \sum_{j=0}^r \vec E_{\vec q,n}^j \varphi_{n,j}(t)$ \\
& error between semidiscrete and fully discrete\\ &  approximation of flux variable \\
$E_{u,n}^i$ & $E_{u,n}^i = E_u(t_{n,i})$: error in node $t_{n,i}$, $i=0,\ldots ,r$\\  
$\vec E_{\vec q,n}^i$ & $\vec E_{\vec q,n}^i = \vec E_{\vec q}(t_{n,i})$: error in node $t_{n,i}$, $i=0,\ldots ,r$\\ \hline 

\end{tabular}
\end{center}

\end{document}